\documentclass[10pt]{article}

\usepackage{latexsym,amsmath,amscd,amssymb,graphics}
\usepackage{enumerate}

\newcommand{\bfi}{\bfseries\itshape}

\makeatletter

\@addtoreset{figure}{section}
\def\thefigure{\thesection.\@arabic\c@figure}
\def\fps@figure{h, t}
\@addtoreset{table}{bsection}
\def\thetable{\thesection.\@arabic\c@table}
\def\fps@table{h, t}
\@addtoreset{equation}{section}

\makeatother

\begin{document}

\newtheorem{theorem}{Theorem}[section]
\newtheorem{definition}[theorem]{Definition}
\newtheorem{lemma}[theorem]{Lemma}
\newtheorem{remark}[theorem]{Remark}
\newtheorem{proposition}[theorem]{Proposition}
\newtheorem{corollary}[theorem]{Corollary}
\newtheorem{example}[theorem]{Example}
      
\makeatletter

\title{The Banach Poisson geometry of 
     multi-diagonal  Toda lattices}
\author{Anatol Odzijewicz$^{1}$ and Tudor
S. Ratiu$^{2}$} 
\addtocounter{footnote}{1} 
\footnotetext{Institute
of Mathematics, University of Bialystok, Lipowa 41,
PL-15424 Bialystok, Poland.
\texttt{aodzijew@uwb.edu.pl}}
\addtocounter{footnote}{1} 
\footnotetext{Section de Math\'ematiques and Bernoulli Center,
\'Ecole Polytechnique F\'ed\'erale de Lausanne. CH--1015 Lausanne.
Switzerland. \texttt{tudor.ratiu@epfl.ch} }
     
\date{}
     
\maketitle 

\makeatother

\maketitle



\begin{abstract}

The Banach Poisson geometry of multi-diagonal Hamiltonian
systems having infinitely many integrals in involution is
studied. 
It is shown that these systems can be considered as
generalizing the semi-infinite Toda lattice which is an
example of a bidiagonal system, a case to which special
attention is given. The generic coadjoint orbits of the Banach
Lie group of bidiagonal bounded operators are studied. It
is shown that the infinite dimensional generalization of the
Flaschka map is a momentum map. Action-angle variables for
the Toda system are constructed.
\end{abstract}

\tableofcontents
 
     
\section{Introduction}
\label{section: Introduction}

Many important conservative systems have a 
Hamiltonian formulation in terms of Lie-Poisson
brackets. With few notable
exceptions, such as the Euler, Poisson-Vlasov, KdV, or
sine-Gordon equations, for example, for infinite
dimensional systems this Lie-Poisson bracket formulation
is mostly formal. It is our belief that these
formal approaches can be given a solid functional analytic
underpinning. The present paper formulates such an
approach for various generalizations of the semi-infinite Toda
lattice. It raises fundamental issues about the nature of
coadjoint orbits for the  Banach Lie
groups having only a finite number of non-zero upper diagonals
and it poses questions about the integrability of certain
generalizations of the Toda lattice in infinite
dimensions by providing a functional analytic framework in
which these problems can be rigorously formulated. The
background of the present work is
\cite{OR} where the theory of Banach Lie-Poisson spaces was
developed.

The paper is organized as follows. The first two sections develop
the theoretical background for the constructions carried out later.
Section
\ref{section: induced and coinduced} presents the general theory
of induced and coinduced Banach Lie-Poisson
structures and derives the analogue of the classical
Adler-Kostant-Symes involution theorem \cite{A,K,S} from this point
of view in the infinite dimensional context. Section 
\ref{section: symplectic induction} introduces the notion of
a momentum map for weak symplectic Banach manifolds and presents the
abstract symplectic induction method in infinite dimensions.

The next two sections concentrate of the Banach Lie-Poisson
geometry of several spaces of trace class operators on a
real separable Hilbert space. The general constructions of Section 2
are implemented explicitly on these spaces in Section
\ref{section: induced and coinduced from ell one}. The
multi-diagonal Banach Lie group, its Lie algebra, and its dual are
introduced and studied. The naturally induced and coinduced Poisson
structures on the preduals of their Banach Lie algebras are
presented.  Section \ref{section: dynamics generated by Casimirs of
ell one} formulates the equations  of motion induced by the Casimir
functions of the Banach Lie-Poisson space of trace class operators
relative to the various induced and coinduced Poisson brackets
discussed previously.

Starting  with Section \ref{sect: bidiagonal case} the
emphasis is  on the important particular case of bidiagonal 
operators, that is, operators having all entries equal to zero
with the possible exception of those on the main and upper 
$k$ diagonal. The Banach Lie group of upper
bidiagonal bounded operators is  studied in detail and the
topological and symplectic structure of the generic coadjoint
orbit is presented. The Banach space analogue of the Flaschka
map is analyzed and its relationship to the coadjoint orbits
is pointed out. There are new, typical infinite dimensional,
phenomena that appear in this context. For example, the
Banach space of trace zero lower bidiagonal trace class
operators does not form a single coadjoint orbit and there
are non-algebraic invariants for the coadjoint orbits. 

Section \ref{sect: generalized Flaschka} uses the
method  of symplectic induction developed Section 
\ref{section: symplectic induction} to derive explicit formulas
that are used for the concrete case of the bidiagonal Banach Lie
group. A generalization of the Flaschka map introduced in the
previous section is presented. This is a Poisson map whose range is
the weak symplectic manifold
$(\ell ^{\infty})^{k-1}\times (\ell^1)^{k-1}$, endowed with a
non-canonical weak symplectic form. Systems with an infinite number
of integrals in involution are also introduced in this section. As an
example of the theory, the semi-infinite Toda lattice is solved
in Section \ref{toda example} using the method of orthogonal
polynomials first introduced, to our knowledge, in \cite{be}. The explicit solution of this system is
obtained, both in action-angle as well as in the original
variables, thereby extending the formulas in \cite{M} from the finite to the semi-infinite Toda lattice.

\medskip

\noindent \textbf{Conventions.} In this paper all Banach
manifolds and Lie groups are real. The
definition of the notion of a Banach Lie subgroup follows
Bourbaki \cite{Bou1972}, that is, a subgroup $H$ of a Banach
Lie group $G $ is necessarily a submanifold (and not just
injectively immersed). In particular, Banach Lie subgroups are
necessarily closed.

\section{Induced and coinduced Banach Lie-Poisson spaces}
\label{section: induced and coinduced}

In this section  we quickly review some material from
\cite{OR} and present some constructions that are necessary
for the development of the ideas in the rest of the paper.

\paragraph{Preliminaries.} Let us recall how a given
Banach Lie-Poisson structure induces and coinduces similar
structures on other Banach spaces. All the proofs of the
statements below can be found in \cite{OR}. Throughout
this paper, unless specified otherwise, all objects are
over
$\mathbb{R}$.

A {\bfi Banach Lie algebra} $(\mathfrak{g}, [\cdot,
\cdot])$ is a Banach space $\mathfrak{g}$ that is also a
Lie algebra such that the Lie bracket is a bilinear
continuous map $\mathfrak{g} \times
\mathfrak{g} \rightarrow \mathfrak{g}$. Thus the adjoint and
coadjoint maps $\operatorname{ad}_x :\mathfrak{g} \rightarrow
\mathfrak{g}$, $\operatorname{ad}_x y:= [x,y]$, and
$\operatorname{ad}_x^\ast : \mathfrak{g}^\ast \rightarrow
\mathfrak{g}^\ast$ are also continuous for each $x \in
\mathfrak{g}$. Here $\mathfrak{g}^\ast$ denotes the
dual of $\mathfrak{g}$, that is, the
Banach space of all linear continuous functionals on
$\mathfrak{g}$.

A {\bfi Banach Lie-Poisson space} $(\mathfrak{b},\{\cdot,
\cdot\})$ is defined to be a real Poisson manifold such
that $\mathfrak{b}$ is a Banach space and the dual
$\mathfrak{b}^\ast \subset C^\infty(\mathfrak{b})$ is a Banach Lie
algebra under the Poisson bracket operation. We need to
explain what does it mean for $\mathfrak{b}$ to be a
Banach Poisson manifold. The Poisson bracket induces the
derivation $h \mapsto \{\cdot , h\}$ on $C
^{\infty}(\mathfrak{b})$ which defines a map $X_h :
\mathfrak{b} \rightarrow \mathfrak{b}^{**}$ by
$\langle X_h (b), Df(b) \rangle = \{f, h \}(b)$
for any $b \in \mathfrak{b}$ and $f$ a smooth real
valued function defined in an open subset of $\mathfrak{b}$
containing $b$. Thus, $X_h(b) \in  \mathfrak{b}^{\ast\ast}
\cong T_b^{\ast\ast}\mathfrak{b}$ and therefore $X_h(b)$ is
not a tangent vector to $\mathfrak{b}$ at $b $. The
requirement that $\mathfrak{b}$ be a Banach
Poisson manifold is that $X_h(b) \in \mathfrak{b} \cong
T_b \mathfrak{b}$ for all $b \in \mathfrak{b}$. 

Denote by $[\cdot, \cdot]$ the restriction of the Poisson 
bracket $\{\cdot, \cdot \}$ from $C^\infty(\mathfrak{b})$
to the Lie subalgebra $\mathfrak{b}^\ast$. The following
criterion characterizes the Banach Lie-Poisson structure.
The Banach space $\mathfrak{b}$ is a Banach Lie-Poisson
space $(\mathfrak{b}, \{\cdot,\cdot\})$ if and only if its
dual $\mathfrak{b}^\ast$ is a Banach Lie algebra
$(\mathfrak{b}^\ast, [\cdot, \cdot ])$ satisfying
$\operatorname{ad}_x^\ast \mathfrak{b} \subset \mathfrak{b}
\subset \mathfrak{b}^{\ast \ast}$ for all $x \in
\mathfrak{b}^\ast$. Moreover, the Poisson bracket of $f, h
\in C^\infty(\mathfrak{b})$ is given by
\begin{equation}\label{general LP}
\{f, h\}(b) = \langle [Df(b), Dh(b)], b
\rangle,
\end{equation}
where $b \in \mathfrak{b}$ and $Df(b) \in
\mathfrak{b}^\ast$ denotes the Fr\'echet derivative of $f
$ at the point $b$. If $h$ is a smooth function on
$\mathfrak{b}$, the associated {\bfi Hamiltonian vector
field\/} is given by
\begin{equation}
\label{general Hamiltonian vector field}
X_h(b) = - \operatorname{ad}^\ast _{Dh(b)} b \in
\mathfrak{b}
\end{equation}
for any $b \in \mathfrak{b}$. Therefore Hamilton's
equations are
\begin{equation}
\label{general Hamilton equations}
\frac{d}{dt}b(t) = - \operatorname{ad}^\ast _{Dh(b(t))} b(t).
\end{equation}

Given two Banach Lie-Poisson spaces $(\mathfrak{b}_1,
\{\,,\}_1)$ and 
$(\mathfrak{b}_2, \{\,,\}_2)$, a smooth map 
$\varphi: \mathfrak{b}_1 \rightarrow \mathfrak{b}_2$ is
said to be {\bfi canonical\/} or a {\bfi Poisson map\/} if 
\begin{equation}
\label{Poisson map definition}
\{f, h\}_2 \circ \varphi = \{f \circ \varphi, h \circ
\varphi \}_1
\end{equation}
for any two smooth locally defined functions $f$ and $h$ on
$\mathfrak{b}_2$. Like in the finite dimensional case,
\eqref{Poisson map definition} is equivalent  to
\begin{equation}
\label{Poisson map condition}
X^2_h \circ  \varphi = T\varphi \circ X^1_{h \circ \varphi}
\end{equation}
for any smooth locally defined function $h$ on
$\mathfrak{b}_2$. Therefore, the  flow of a
Hamiltonian  vector field is a Poisson map and Hamilton's
equations $\dot{f} = \{f, h \}$ in Poisson bracket 
formulation are valid. If the Poisson map $\varphi$ is, in
addition, linear, then it is called a {\bfi linear Poisson
map}.

Given the Banach Lie-Poisson spaces $(\mathfrak{b}_1,
\{\,,\}_1)$ and  $(\mathfrak{b}_2, \{\,,\}_2)$ there is a
unique Banach Poisson structure 
$\{\,,\}$ on the product space $\mathfrak{b}_1 \times
\mathfrak{b}_2$ such that:
\begin{enumerate}
\item[{\rm (i)}] the canonical projections 
$\pi_1: \mathfrak{b}_1 \times \mathfrak{b}_2
\rightarrow \mathfrak{b}_1$  and $\pi_2: \mathfrak{b}_1
\times \mathfrak{b}_2 \rightarrow \mathfrak{b}_2$ are
Poisson maps;
\item[{\rm (ii)}] $\pi_1^\ast (C^\infty(\mathfrak{b}_1))$
and $\pi_2^\ast (C^\infty(\mathfrak{b}_2))$ are Poisson
commuting subalgebras of $C^\infty(\mathfrak{b}_1 \times
\mathfrak{b}_2)$.
\end{enumerate}
This unique Poisson structure on $\mathfrak{b}_1 \times
\mathfrak{b}_2$ is called the  {\bfi product Poisson
structure\/} and its bracket is given by the formula
\begin{equation}
\label{product bracket}
\{f, g\}(b_1, b_2) = \{f_{b_2}, g_{b_2}\}_1(b_1) + 
\{f_{b_1}, g_{b_1}\}_2(b_2),
\end{equation}
where $f_{b_1}, g_{b_1}  \in C^\infty(\mathfrak{b}_2)$ and 
$f_{b_2}, g_{b_2} \in C^\infty(\mathfrak{b}_1)$ are the
partial functions given by $f_{b_1}(b_2) := f_{b_2}(b_1) :=
f(b_1, b_2)$ and $g_{b_1}(b_2) := g_{b_2}(b_1) :=
g(b_1, b_2)$. In addition, this formula shows that this
unique Banach Poisson structure is Lie-Poisson and that
the inclusions $\iota_1: \mathfrak{b}_1 \hookrightarrow
\mathfrak{b}_1 \times \mathfrak{b}_2$, $\iota_2:
\mathfrak{b}_2 \hookrightarrow \mathfrak{b}_1 \times
\mathfrak{b}_2$ given by $\iota_1(b_1): = (b_1, 0)$ and
$\iota_2 (b_2) : = (0,b_2)$, respectively, are also linear
Poisson maps.

\paragraph{Induced Structures.}
Let $\mathfrak{b}_1$ be a Banach space,
$(\mathfrak{b}, \{\cdot, \cdot \})$ a Banach
Lie-Poisson space, and $\iota : \mathfrak{b}_1
\hookrightarrow \mathfrak{b}$ an injective continuous
linear map with closed range. Then $\ker \iota^\ast$ is an
ideal in the Banach Lie algebra $(\mathfrak{b}^\ast,
[\cdot,
\cdot])$ if and only if $\mathfrak{b}_1$ carries a unique 
Banach Lie-Poisson bracket $\{\cdot, \cdot \}_1^{\rm ind}$
such that
\begin{equation}
\label{iota poisson}
\{F \circ \iota, G \circ \iota\}_1^{\rm ind} = \{F, G\}
\circ \iota
\end{equation}
for any $F,G \in C ^{\infty}(\mathfrak{b})$; see
Proposition 4.10 in \cite{OR}. This Poisson structure on
$\mathfrak{b}_1$ is said to be {\bfi induced\/} by the
mapping $\iota$ and it is given by
\begin{equation}
\label{induced one bracket}
\{f, g\}_1^{\rm ind}(b_1) = \left\langle [
\iota^\ast]\left(\left[ [\iota^\ast]^{-1} (Df(b_1)),
[\iota^\ast]^{-1} (Dg(b_1))
\right]_1\right), b_1 \right\rangle
\end{equation}
for any $f, g \in C ^{\infty}(\mathfrak{b}_1) $ and $b_1
\in \mathfrak{b}_1$, where $[\iota ^\ast]:
\mathfrak{b}^\ast/\ker \iota^\ast
\rightarrow \mathfrak{b}_1 ^\ast$ is the Banach space
isomorphism induced by $\iota^\ast: \mathfrak{b}^\ast
\rightarrow \mathfrak{b}_1^\ast$ and $[\cdot, \cdot]_1$
denotes the Lie bracket on the
quotient Lie algebra $\mathfrak{b}^\ast/\ker
\iota^\ast$.

Let us assume now that the range $\iota(\mathfrak{b}_1)$
is a closed split subspace of $\mathfrak{b}$, that is,
there exists a projector $R= R^2 :
\mathfrak{b}\rightarrow \mathfrak{b}$ such that
$\iota(\mathfrak{b}_1) = R(\mathfrak{b})$. Taking in
\eqref{iota poisson} $F: = f \circ \iota^{-1}\circ R, G: =
g \circ \iota^{-1}\circ R \in C^{\infty} (\mathfrak{b})$ for
$f,g \in C ^{\infty}( \mathfrak{b}_1)$ and noting that 
$\iota^{-1}
\circ R
\circ
\iota =
\operatorname{id}_{\mathfrak{b}_1}$, we get
\begin{align}
\label{induced bracket third form}
\{f, g\}_1^{\rm ind}(b_1) 
& = \{f \circ \iota^{-1}\circ R, g \circ \iota^{-1}\circ
R\}(\iota(b_1)) \nonumber \\
&= \left\langle \left[D(f \circ
\iota^{-1}\circ R)(\iota(b_1)), D(g \circ
\iota^{-1}\circ R)(\iota(b_1)) \right],
\iota(b_1) \right\rangle.
\end{align}
We shall make use of this formula in 
\S\ref{section: induced and coinduced from ell one}.

\medskip

We return now to the general case, that is, we consider an
arbitrary quasi-immersion $\iota: \mathfrak{b}_1 \hookrightarrow
\mathfrak{b}$ of Banach spaces which means that the 
range $\iota(\mathfrak{b}_1)$ is closed but does not 
necessarily possess a closed complement.

\begin{proposition}
\label{induction proposition}
Let $\iota: \mathfrak{b}_1 \hookrightarrow \mathfrak{b}$
be a quasi-immersion of Banach Lie-Poisson spaces {\rm (}so
$\operatorname{range}\iota $ is a closed subspace of
$\mathfrak{b}$ and $\ker \iota^\ast$ is an ideal in
the Banach Lie algebra $\mathfrak{b}^\ast${\rm )}. Assume
that there is a connected Banach Lie group $G$ with
Banach Lie algebra $\mathfrak{g}: =\mathfrak{b}^\ast$.
Then the $G $-coadjoint orbit $\mathcal{O}_{\iota(b_1)}: =
\operatorname{Ad}^\ast_G
\iota(b_1)$ is contained in $\iota(\mathfrak{b}_1)$ for any $b_1 \in
\mathfrak{b}_1$. In addition, if $N \subset G$ is a
closed connected normal Lie subgroup of $G$ whose Lie
algebra is $\ker \iota ^\ast$, then the $N$-coadjoint
action restricted to $\iota(\mathfrak{b}_1)$ is trivial.
Therefore the Banach Lie group $G/N: = \{[g] : =
gN \mid g \in G\}$ naturally acts on
$\iota(\mathfrak{b}_1)$ and the orbit of $\iota(b_1) $
under this action coincides with $\mathcal{O}_{\iota(b_1)}$
for any $b_1 \in \mathfrak{b}_1 $.
\end{proposition}

\noindent \textbf{Proof.} Since $\ker \iota^\ast$ is an
ideal in $\mathfrak{g} = \mathfrak{b}^\ast$, it follows
that $[x, y] \in \ker \iota^\ast$ for all $x \in
\mathfrak{g}$ and $y \in \ker \iota^\ast$. Therefore,
since $\ker \iota^\ast$ is closed in
$\mathfrak{g}$, it follows that $\operatorname{Ad}
_{\exp x} y = e^{\operatorname{ad}_x} y \in \ker
\iota^\ast$ for any $x \in \mathfrak{g}$ and $y \in \ker
\iota^\ast$. This shows that for any $g \in G$ in an open
neighborhood of the identity element of $G$ we have 
$\operatorname{Ad}_{g} \ker \iota^\ast \subset \ker
\iota^\ast$. Since $G $ is connected, it is generated by
a neighborhood of the  identity and  we conclude that 
$\operatorname{Ad}_{g} \ker \iota^\ast \subset \ker
\iota^\ast$ for any $g \in G $.

The upper index $^\circ$ on a set denotes the annihilator of
that set relative to a duality pairing; the annihilator of
a set is always a vector subspace. Let $b_1 \in
\mathfrak{b}_1$ and $g \in G$. Since $\ker \iota^\ast =
\iota(\mathfrak{b}_1) ^\circ $, closedness of
$\iota(\mathfrak{b}_1)$ in $\mathfrak{b}$ implies that $(\ker
\iota ^\ast) ^\circ = \iota(\mathfrak{b}_1)$. Thus, for any
$g \in G $ and $x \in \ker \iota^\ast$, we have
\[
\langle \operatorname{Ad}^\ast_g \iota(b_1), x \rangle = 
\langle \iota(b_1), \operatorname{Ad}_g x \rangle = 0
\]
which proves that $\operatorname{Ad}^\ast_G\iota(b_1)
\subset \iota(\mathfrak{b}_1)$.

Now let $N \subset G$ be a closed connected normal Lie
subgroup of $G$ with Banach Lie algebra $\ker \iota ^\ast
\subset \mathfrak{g}$. For any $b_1 \in \mathfrak{b}_1$, $x
 \in \mathfrak{g} = \mathfrak{b}^\ast$, $y \in \ker
\iota^\ast$, we have
\[
\langle\operatorname{ad}^\ast_y \iota(b_1), x \rangle 
= \langle \iota(b_1), [y, x]\rangle = 0
\] 
since $\ker \iota^\ast$ is an ideal in $\mathfrak{g}$ and
$\ker \iota^\ast = \iota(\mathfrak{b}_1) ^\circ$. Since
this is valid for all $x  \in \mathfrak{g}$, it follows
that $\operatorname{ad}^\ast_y \iota(b_1) = 0 $ for all $y
\in \ker \iota^\ast$ and $b_1 \in \mathfrak{b}_1$. Using 
the exponential map, this shows that
$\operatorname{Ad}^\ast_n \iota(b_1) = \iota(b_1)$ for any
$n $ in a neighborhood of the identity in $N $. Since $N $
is connected, it is generated by a neighborhood of the
identity and we conclude that $\operatorname{Ad}^\ast_n
\iota(b_1) = \iota(b_1)$ for all $n \in N $.

The quotient $G/N : = \{[g] : = gN \mid g \in G \}$ is a
Banach Lie group and the projection $G \rightarrow G/N $ is
a smooth surjective submersive Banach Lie group
homomorphism (\cite{Bou1972}, Chapter III, \S 1.6). Since
the coadjoint action of $N$ on $\iota(\mathfrak{b}_1)$ is
trivial, the Banach Lie group $G/N $ acts smoothly on
$\iota(\mathfrak{b}_1)$ by $[g] \cdot \iota(b_1): =
\operatorname{Ad}^\ast_{g ^{-1}} \iota(b_1)$. The orbit of
a fixed element
$\iota(b_1) \in \iota(\mathfrak{b}_1) $ by this group
action is obviously equal to the $G$-orbit
$\mathcal{O}_{\iota(b_1)}$.
\quad $\blacksquare$

\paragraph{Coinduced Structures.} Let $(\mathfrak{b}, 
\{\,, \})$ be a Banach Lie-Poisson space 
and $\pi : \mathfrak{b} \rightarrow \mathfrak{b}_1$ 
a continuous linear surjective map onto the Banach space
$\mathfrak{b}_1$. Then $\mathfrak{b}_1$ carries a
unique Banach Lie-Poisson bracket $\{\,,\}_1^{\rm coind} $
making $\pi$ into a linear Poisson map, that is,
\begin{equation}
\label{pi poisson}
\{f \circ \pi, g \circ \pi\} = \{f, g\}_1^{\rm coind} \circ
\pi
\end{equation}
for any $f, g \in C ^{\infty}(\mathfrak{b}_1) $ if and only
if $\pi^\ast(\mathfrak{b}_1 ^\ast) \subset
\mathfrak{b}^\ast$ is closed under the Lie bracket
$[\cdot\,, \cdot]$ of $\mathfrak{b}^\ast$; see Proposition
4.8 of \cite {OR}. This unique Poisson structure on
$\mathfrak{b}_1$ is said to be {\bfi coinduced\/} by the
Banach Lie-Poisson structure on
$\mathfrak{b}$ and the linear continuous map
$\pi$. It should be noted that $\operatorname{im}
\pi^\ast$ is a closed subspace of $\mathfrak{b}^\ast$
since $\operatorname{im} \pi^\ast = (\ker \pi)^\circ$. 
To determine the coinduced bracket on $\mathfrak{b}_1$ note
that $\pi^\ast: \mathfrak{b}_1 ^\ast \rightarrow
\mathfrak{b}^\ast$ is an injective linear continuous map
whose closed range is a Banach Lie subalgebra of
$\mathfrak{b}^\ast$. Thus, on $\operatorname{im}\pi^\ast$
we can invert $\pi^\ast$. The coinduced bracket on
$\mathfrak{b}_1$ has then the form
\begin{equation}
\label{coinduced one bracket}
\{f, g\}_1^{\rm coind}(b_1) = \left\langle
(\pi^\ast)^{-1}\left[
\pi^\ast(Df(b_1)), \pi^\ast(Dg(b_1))\right], b_1
\right\rangle
\end{equation}
for any $f, g \in C ^{\infty}(\mathfrak{b}_1) $ and $b_1
\in \mathfrak{b}_1 $.
\medskip

Let us assume that $\ker \pi$ admits a closed complement.
This is equivalent to the existence of a linear continuous
injective map  $\iota: \mathfrak{b}_1
\hookrightarrow \mathfrak{b}$ with closed range such that
$\pi\circ \iota = \operatorname{id}_{\mathfrak{b}_1}$.
Thus \eqref{pi poisson} implies that 
\begin{equation}
\label{coinduced bracket}
\{f, g\}_1^{\rm coind}= \{f \circ \pi, g \circ \pi\}
\circ \iota.
\end{equation}
for any $f, g \in C ^{\infty}(\mathfrak{b}_1) $.
\medskip

Assume now
that the Banach Lie-Poisson space
$\mathfrak{b}$ splits into a direct sum $\mathfrak{b} =
\mathfrak{b}_1 \oplus \mathfrak{b}_2 $ of closed Banach
subspaces. Denote by $R_j: \mathfrak{b} \rightarrow
\mathfrak{b}$ the projection onto $\mathfrak{b}_j$, for
$j = 1,2$. So we have the following relations: 
$R_1 + R_2 = \operatorname{id}_{\mathfrak{b}},\; R_1^2
= R_1,\; R_2^2 = R_2,\; R_2 R_1 = R_1R_2 = 0,\; 
\mathfrak{b}_1 := \operatorname{im}R_1$, and
$\mathfrak{b}_2 := \operatorname{im}R_2$.
Dualizing we get the projectors $R_1^\ast , R_2^\ast:
\mathfrak{b}^\ast \rightarrow \mathfrak{b}^\ast$
satisfying $R_1^\ast + R_2^\ast =
\operatorname{id}_{\mathfrak{b}^\ast},\; (R_1^\ast)^2
= R_1^\ast,\; (R_2^\ast)^2 = R_2^\ast,\; R_2^\ast R_1^\ast =
R_1^\ast R_2^\ast = 0$. The relationship between these
spaces is given by
\begin{align}
\label{b spaces relation}
\ker R_1 = \operatorname{im}R_2 = \mathfrak{b}_2
\qquad &\text{and} \qquad 
\ker R_2 = \operatorname{im}R_1 = \mathfrak{b}_1 \\
\label{b star spaces relation}
\ker R_1^\ast = \operatorname{im}R_2^\ast =
(\operatorname{im}R_1)^\circ \cong \mathfrak{b}_2^\ast 
\qquad &\text{and} \qquad
\ker R_2^\ast = \operatorname{im}R_1^\ast =
(\operatorname{im}R_2)^\circ \cong \mathfrak{b}_1^\ast \\
\label{sum relations}
\mathfrak{b} = \mathfrak{b}_1 \oplus \mathfrak{b}_2
\qquad &\text{and} \qquad
\mathfrak{b}^\ast = \mathfrak{b}_2^\circ
\oplus\mathfrak{b}_1^\circ \cong \mathfrak{b}_1 ^\ast
\oplus \mathfrak{b}_2 ^\ast. 
\end{align}
Let $\iota_j : \mathfrak{b}_j \hookrightarrow \mathfrak{b}$
be the inclusion determined by the splitting $\mathfrak{b}
= \mathfrak{b}_1 \oplus \mathfrak{b}_2$ for $j = 1,2$.
Denote by $\pi_j : \mathfrak{b} \rightarrow
\mathfrak{b}_j$ the projection determined by the
projector $R_j: \mathfrak{b} \rightarrow \mathfrak{b}$,
that is, $\iota_j \circ \pi_j = R_j$ and note that $\pi_j
\circ \iota_j= \operatorname{id}_{\mathfrak{b}_j}$. We
summarize these notations in the following diagram.

     \unitlength=5mm
     \begin{center}
     \begin{picture}(9.6,5.5)
     \put(4.8,5){\makebox(0,0){$\mathfrak{b}$}}
     \put(1,0.5){\makebox(0,0){$\mathfrak{b}_1$}}
     \put(9,0.5){\makebox(0,0){$\mathfrak{b}_2$}}
     \put(5.5,4.5){\vector(1,-1){3}}
     \put(8,1.5){\vector(-1,1){3}}
     \put(4,4.5){\vector(-1,-1){3}}
     \put(1.5,1.5){\vector(1,1){3}}
     \put(2.1,3.5){\makebox(0,0){$\pi_1$}}
     \put(3.4,2.5){\makebox(0,0){$\iota_1$}}
     \put(7.3,3.5){\makebox(0,0){$\pi_2$}}
     \put(6.4,2.5){\makebox(0,0){$\iota_2$}}
     \end{picture}
     \end{center}

From \eqref{coinduced bracket} we get
\begin{equation}
\label{coinduced bracket with projector}
\{f,g\}_j^{\rm coind} = \{f \circ \pi_j, g \circ \pi_j\}
\circ
\iota_j
\end{equation}
or, explicitly
\begin{equation}
\label{coinduced bracket with projector evaluated}
\{f,g\}_j^{\rm coind} (b_j)= \left\langle \left[D(f \circ
\pi_j)(\iota_j(b_j)), D(g \circ \pi_j)(\iota_j(b_j))
\right], \iota_j(b_j) \right\rangle, \; \text{where} \; b_j \in \mathfrak{b}_j.
\end{equation}
The following proposition presents some properties of the
induced and coinduced structures on $\mathfrak{b}_1 $ and
$\mathfrak{b}_2 $.

\begin{proposition}
\label{decomposition proposition}
Assume that $\operatorname{im}R_1^\ast$ and
$\operatorname{im}R_2 ^\ast$ are Banach Lie subalgebras of
$\mathfrak{b}^\ast$. Then:
\begin{itemize}
\item[{\rm \textbf{(i)}}] $\mathfrak{b}_j$  has a Banach
Lie-Poisson structure coinduced by $\pi_j$ and  the
expression of the coinduced bracket
$\{\,,\}_{j}^{\rm coind}$ on $\mathfrak{b}_j$ is given by
\eqref{coinduced bracket with projector}.
The  Hamiltonian vector field of
$h \in C ^{\infty}(\mathfrak{b}_j)$ at $b_j \in
\mathfrak{b}_j $ is given by
\begin{equation}
\label{coinduced vector field}
 X_h(b_j)  = -\pi_j\left(
\operatorname{ad}^\ast_{\pi_j^\ast
Dh(b_j)} \iota_j(b_j) \right), \qquad
j = 1,2,
\end{equation}
where $Dh(b_j) \in \mathfrak{b}_j^\ast$ and
$\operatorname{ad}_x $ is the adjoint action of $x \in
\mathfrak{b}^\ast$ on $\mathfrak{b}^\ast$.

\item[{\rm \textbf{(ii)}}] The Banach space isomorphism
$R:=
\frac{1}{2}(R_1 - R_2) :
\mathfrak{b} \rightarrow \mathfrak{b}$ defines a new
Banach Lie-Poisson structure 
\begin{equation}
\label{R Poisson bracket}
\{f, g\}_R(b): = \langle [R^\ast Df(b),
Dg(b) ] + [ Df(b), R^\ast Dg (b)], b \rangle
\end{equation}
on $\mathfrak{b}$, $f, g \in C^{\infty}
(\mathfrak{b})$, that coincides with the product structure
on $\mathfrak{b}_1
\times \overline{\mathfrak{b}}_2 $, where
$\mathfrak{b}_1$ carries the coinduced bracket
$\{\,,\}_{1}^{\rm coind}$ and 
$\overline{\mathfrak{b}}_2 $ denotes $\mathfrak{b}_2 $
endowed with the Lie-Poisson bracket $- \{\,,\}_{2}^{\rm
coind}$.

\item[{\rm \textbf{(iii)}}] The inclusion maps $\iota_1:
(\mathfrak{b}_1, \{\,,\}^{\rm coind}_{1}) \hookrightarrow
(\mathfrak{b},
\{\,,\}_R)$ and
$\iota_2: (\overline{\mathfrak{b}}_2,
\{\,,\}_{2}^{\rm coind})
\hookrightarrow (\mathfrak{b}, \{\,,\}_R)$ are  linear
injective Poisson maps with closed range.

\item[{\rm \textbf{(iv)}}]  The map $\iota_j$ induces from 
$(\mathfrak{b}, \{\,,\}_R)$ a Banach Lie-Poisson structure 
on $\mathfrak{b}_j $ which coincides with the coinduced
structure described in {\rm \textbf{(i)}},  for $j = 1,2$. 
\end{itemize}
\end{proposition} 

\noindent \textbf{Proof.} \textbf{(i)} By hypothesis, the
range $\operatorname{im} R_j ^\ast$ of the map $R_j^\ast:
\mathfrak{b}^\ast
\rightarrow
\mathfrak{b}^\ast$ is a Banach Lie subalgebra of
$\mathfrak{b}^\ast$. Thus $\pi_j$ coinduces a Banach
Lie-Poisson structure on $\mathfrak{b}_j^\ast$. Let $h \in
C ^{\infty}(\mathfrak{b}_j)$ and note that for any
function $f \in C^{\infty}(\mathfrak{b}_j)$ and $b_j \in
\mathfrak{b}_j$ we have
\begin{align*}
\left\langle Df(b_j), X_h(b_j) \right\rangle 
& = \{f, h\}_j^{\rm coind}(b_j) 
= \left\langle \left[D(f \circ \pi_j)(\iota_j(b_j)),
D(h \circ \pi_j)(\iota_j(b_j)) 
\right], \iota_j(b_j) \right\rangle \\
& = \left\langle \left[\pi_j ^\ast Df(b_j), \pi_j ^\ast
Dh(b_j) \right], \iota_j(b_j) \right\rangle \\
& = \left\langle \pi_j ^\ast Df(b_j),
-\operatorname{ad}^\ast_{\pi_j ^\ast Dh(b_j)} \iota_j(b_j)
\right\rangle \\
& = \left\langle Df(b_j),
-\pi_j \operatorname{ad}^\ast_{\pi_j ^\ast Dh(b_j)}
\iota_j(b_j) \right\rangle,
\end{align*}
which proves formula \eqref{coinduced vector field}.

\textbf{(ii)} Let $b = b_1 + b_2 \in \mathfrak{b}_1 \oplus
\mathfrak{b}_2 $. Then $R_j(b) = b_j$, for $j = 1,2$. A
direct verification shows then that
\begin{align*}
\{f, g\}_R(b)& = \langle [R^\ast Df(b),
Dg(b) ] + [ D(b), R^\ast Dg(b)], b
\rangle\\
& = \frac{1}{2}\langle [R_1^\ast Df(b) -
R_2^\ast Df(b), R_1^\ast Dg(b) + R_2^\ast Dg (b)], b
\rangle\\ 
&\qquad   +\frac{1}{2} \langle [R_1^\ast Df(b) +
R_2^\ast Df(b), R_1^\ast Dg(b) -
R_2^\ast Dg(b)], b \rangle\\
& = \langle  [R_1^\ast Df(b), 
R_1^\ast Dg(b)], R_1b + R_2b \rangle -
\langle [R_2^\ast Df(b), R_2^\ast Dg(b)],
R_1b + R_2b \rangle \\
&= \langle [R_1^\ast Df(b), 
R_1^\ast Dg(b)], R_1b \rangle -
\langle [R_2^\ast Df(b), 
R_2^\ast Dg(b)],  R_2b \rangle \\
& = \{f_{b_2}, g_{b_2}\}_{1}^{\rm coind}(b_1) - 
\{f_{b_1}, g_{b_1}\}_{2}^{\rm coind}(b_2),
\end{align*}
where in the third equality we have used the fact 
that $[R_1^\ast Df(b),  R_1^\ast Dg(b)] \in
\operatorname{im}R_1 ^\ast = (\operatorname{im}R_2)^\circ$
and $[R_2^\ast Df(b),  R_2^\ast Dg(b)] \in
\operatorname{im}R_2 ^\ast =
(\operatorname{im}R_1)^\circ$ and $b = b_1 + b_2 $ with
$b_j \in \mathfrak{b}_j $. To prove the last equality above
it suffices to note that 
\[
D_1 f_{b_2}(b_1) \cdot \delta b_1  = Df(b) \cdot \delta b_1
= Df(b) \cdot R_1\delta b_1 \;\text{and } \;
D_2 f_{b_1}(b_2) \cdot \delta b_2  = Df(b) \cdot \delta b_2
= Df(b) \cdot R_2\delta b_2
\]
for any $\delta b_j \in \mathfrak{b}_j$, where $D_j$ is
the Fr\'echet derivative on $\mathfrak{b}_j $, for $j =
1,2 $. The last expression is that of the product Banach
Lie-Poisson structure on $\mathfrak{b}_1 \times
\overline{\mathfrak{b}}_2$ (see \eqref{product bracket}).

\textbf{(iii)} This is an immediate consequence of
\textbf{(ii)} and the general fact, recalled earlier for
products of Banach Lie-Poisson spaces, that these
inclusions are Poisson maps with closed range.

{\rm \textbf{(iv)}} Let $\{\,,\}_j^{\operatorname{ind}}$ and
$\{\,,\}_j^{\operatorname{coind}}$ be the induced and
coinduced brackets on $\mathfrak{b}_j$ from $( \mathfrak{b},
\{\cdot  , \cdot \}_R)$ and $( \mathfrak{b},
\{\cdot  , \cdot \})$, respectively. Therefore, 
\begin{equation}
\label{induced first relation}
\{F, G\}_R \circ  \iota_j = \{F \circ \iota_j , G \circ
\iota_j\}_j^{\operatorname{ind}}
\end{equation} 
for any $F, G \in C ^{\infty}(\mathfrak{b}) $ and, by  
\eqref{coinduced bracket with projector},
\begin{equation}
\label{coinduced second  relation}
\{f, g \}_j^{\operatorname{coind}} = (-1)^{j-1}\{f \circ
\pi_j, g \circ \pi_j \} \circ \iota_j
\end{equation}
for any $f, g \in C ^{\infty}(\mathfrak{b}_j)$. Apply
relation \eqref{induced first relation} to the functions
$F: = f \circ \pi_j, G: = g \circ \pi_j $ and use $\pi_j
\circ \iota_j = \operatorname{id}_{\mathfrak{b}_j} $,
$\pi_j \circ R = \frac{1}{2}(-1)^{j-1}\pi_j $, and
\eqref{coinduced second  relation} to get for any $b_j
\in
\mathfrak{b}_j
$
\begin{align*}
\{f, g\}_j^{\operatorname{ind}}(b_j) 
&= \{f \circ \pi_j , g \circ \pi_j \}_R\left(\iota_j(b_j)
\right) \\ 
& = \langle [R^\ast D(f\circ
\pi_j)(\iota_j(b_j)),  D(g \circ
\pi_j)(\iota_j(b_j))],\iota_j(b_j) \rangle\\ 
&\qquad +
\langle [D(f \circ \pi_j)(\iota_j(b_j)),  R^\ast D(g \circ
\pi_j)(\iota_j(b_j))],
\iota_j(b_j) \rangle\\  
& = \langle [R^\ast \pi_j^\ast Df(b_j), 
\pi_j^\ast Dg(b_j)], \iota_j(b_j) \rangle\\ 
&\qquad + \langle [\pi_j^\ast
Df(b_j), R^\ast \pi_j^\ast Dg(b_j)],
\iota_j(b_j) \rangle\\ 
& = (-1)^{j-1}\langle [\pi_j^\ast Df(b_j)), 
\pi_j^\ast Dg(b_j)], \iota_j(b_j) \rangle\\ 
& = (-1)^{j-1}\langle [D(f\circ \pi_j)(\iota_j(b_j)),
D(g\circ \pi_j)(\iota_j(b_j))], \iota_j(b_j) \rangle\\ 
&= (-1)^{j-1}\{f\circ \pi_j, g \circ
\pi_j\}(\iota_j(b_j))\\ & = \{f, g
\}_j^{\operatorname{coind}}(b_j).
\quad \blacksquare
\end{align*}

This proposition implies the following involution theorem.

\begin{corollary}
\label{involution corollary}
In the notations and hypotheses of Proposition
\ref{decomposition proposition} we have:
\begin{itemize}

\item[{\rm \textbf{(i)}}] The Casimir functions on
$(\mathfrak{b}, \{ \cdot  , \cdot  \})$ are in involution on $( \mathfrak{b}, \{ \cdot ,\cdot \}_R)$ and restrict to functions in involution on
$\mathfrak{b}_j$, $j = 1,2 $.

\item[{\rm \textbf{(ii)}}] If $H$ is a Casimir function
on $\mathfrak{b}$, then its restriction $H \circ \iota_j$
to $\mathfrak{b}_j $ has the Hamiltonian vector field
\begin{equation}
\label{j Casimir vector field}
\begin{array}{ll}
&X_{H \circ \iota_1}(b_1) =
\pi_1 \big(\operatorname{ad}^\ast_{R_2^\ast
DH(\iota_1(b_1))} \iota_1(b_1) \big)\\
&X_{H \circ \iota_2}(b_2)  =
\pi_2 \big(\operatorname{ad}^\ast_{R_1^\ast
DH(\iota_2(b_2))} \iota_2(b_2)\big)
\end{array}
\end{equation}
for any $b_1 \in \mathfrak{b}_1 $ and $b_2 \in
\mathfrak{b}_2 $, where $\iota_j: \mathfrak{b}_j
\hookrightarrow \mathfrak{b}$ is the inclusion,  $j = 1,2$.

\end{itemize}

\end{corollary}

\noindent \textbf{Proof.} {\rm \textbf{(i)}} Let $F, H
\in  C ^{\infty}(\mathfrak{b})$ be Casimir functions for
the Lie-Poisson bracket $\{\,,\}$, that is,
$\operatorname{ad}^\ast_{DF(b)}b = 0$ and 
$\operatorname{ad}^\ast_{DH(b)}b = 0$ for any $b \in
\mathfrak{b}$. Therefore
\begin{align*}
\{F, H\}_R(b) &= \left\langle [R^\ast DF(b),
DH(b)] + [DF(b), R^\ast DH (b)], b \right\rangle \\
& = - \left\langle  R^\ast DF(b),
\operatorname{ad}^\ast_{DH(b)}b \right\rangle +
\left\langle R^\ast DH(b),
\operatorname{ad}^\ast_{DF(b)}b \right\rangle = 0
\end{align*}
which shows that $F$ and $H$ are in involution relative
to $\{\,,\}_R $. Then statements (iii) and (iv)
of Proposition \ref{decomposition proposition} show that
$F \circ \iota_j , H \circ \iota_j $ are also in involution
on $\mathfrak{b}_j$, $j = 1,2$. 

\textbf{(ii)} Since $H$ is a Casimir function on
$\mathfrak{b}$, we have
$\operatorname{ad}^\ast_{DH(b)} b = 0$ for any $b
\in \mathfrak{b}$. Therefore, since $R_1 ^\ast + R_2 ^\ast 
= \operatorname{id}_{\mathfrak{b}^\ast}$, we get for any
$b_1 \in \mathfrak{b}_1 $
\[
0=\operatorname{ad}^\ast_{DH(\iota_1(b_1))} \iota_1(b_1)
= \operatorname{ad}^\ast_{R_1 ^\ast DH(\iota_1(b_1))}
\iota_1(b_1) + 
\operatorname{ad}^\ast_{R_2 ^\ast DH(\iota_1(b_1))}
\iota_1(b_1).
\]
A similar relation holds for any $b_2 \in \mathfrak{b}_2$.
So, we have
\begin{equation}
\label{casimir relation on projectors}
-\operatorname{ad}^\ast_{R_j ^\ast DH(\iota_j(b_j))}
 = \operatorname{ad}^\ast_{R_{j+1} ^\ast DH(\iota_j(b_j))},
\end{equation}
where $j$ is taken modulo $2 $.

Since  $\iota_j \circ  \pi_j = R_j $ we get
\begin{align*}
\pi_j ^\ast D(H \circ \iota_j)(b_j) 
&= D(H \circ \iota_j)(b_j) \circ \pi_j 
= DH(\iota_j(b_j)) \circ \iota_j \circ  \pi_j \\
&= DH(\iota_j(b_j)) \circ R_j  
= R_j^\ast DH(\iota_j(b_j)),   
\end{align*}
so \eqref{coinduced vector field} and 
\eqref{casimir relation on projectors} yield
\begin{align}
\label{intermediate AKS computation}
\iota_j \left(X_{H \circ \iota_j}(b_j) \right) 
&= -( \iota_j \circ \pi_j) \left(
\operatorname{ad}^\ast_{\pi_j^\ast
D(H \circ \iota_j)(b_j)} \iota_j(b_j) \right)
= - R_j \left( \operatorname{ad}^\ast_{R_j^\ast
DH(\iota_j(b_j))} \iota_j(b_j) \right) \nonumber \\
& = R_j \left( \operatorname{ad}^\ast_{R_{j+1}^\ast
DH(\iota_j(b_j))} \iota_j(b_j) \right)
= \operatorname{ad}^\ast_{R_{j+1}^\ast
DH(\iota_j(b_j))} \iota_j(b_j).
\end{align}
The last equality follows from the fact that
$\operatorname{ad}^\ast_{R_{j+1}^\ast x} \iota_j(b_j)
\in \operatorname{im}R_j = \operatorname{im}\iota_j$ for
any $x
\in
\mathfrak{b}^\ast$ and $b_j \in \mathfrak{b}_j$. Indeed,
for  any $y \in (\operatorname{im}R_j)^\circ =
\operatorname{im} R_{j+1}^\ast$ we have
\[
\left\langle \operatorname{ad}^\ast_{R_{j+1}^\ast x}
\iota_j(b_j), y \right\rangle = 
\left\langle \iota_j(b_j), \left[R_{j+1}^\ast x, y  \right]
\right\rangle = 0
\]
because $\left[R_{j+1}^\ast x, y  \right] \in
\operatorname{im} R_{j+1} ^\ast = (\operatorname{im}R_j) ^
\circ$ by hypothesis (the image of
$R_{j+1}^\ast$ is a Banach Lie subalgebra of
$\mathfrak{b}^\ast$) and $\iota_j(b_j) \in
\operatorname{im}R_j$. Therefore, 
$\operatorname{ad}^\ast_{R_{j+1}^\ast x} \iota_j(b_j) \in
(\operatorname{im}R_j) ^{\circ \circ } =
\overline{\operatorname{im}R_j} = \operatorname{im}R_j$.

Finally, applying $\pi_j $ to \eqref{intermediate AKS
computation} yields \eqref{j Casimir vector field}.
\quad $\blacksquare$
\medskip

Taken together, Proposition \ref{decomposition proposition}
and Corollary \ref{involution corollary} give a version of
the Adler-Kostant-Symes Theorem (see
\cite{A, K, S, R}) formulated with the necessary additional
hypotheses in the context of Banach Lie-Poisson spaces.

\begin{proposition}
\label{double diagram}
Let $(\mathfrak{b}, \{\,,\})$ be a Banach Lie-Poisson
space and let $R_1, R_3: \mathfrak{b} \rightarrow
\mathfrak{b}$ be projectors. Assume that
$\operatorname{im}R_{21} = \operatorname{im}R_{23} = :
\mathfrak{b}_2$, where
$R_{21}:= \operatorname{id}_{\mathfrak{b}} - R_1 $, $R_{23}
: = \operatorname{id}_{\mathfrak{b}} - R_3$, and denote
$\mathfrak{b}_1: =
\operatorname{im}R_1$, $\mathfrak{b}_3 : =
\operatorname{im}R_3$.
 We summarize
this situation  in the diagram

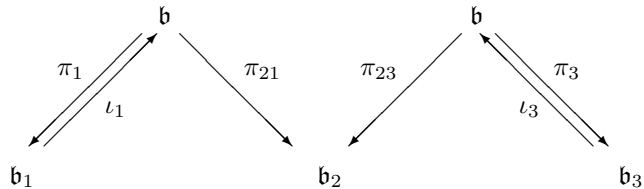
\unitlength=5mm
\begin{center}
\begin{picture}(18,5.5)
\put(4.6,5){\makebox(0,0){$\mathfrak{b}$}}
\put(.8,0.7){\makebox(0,0){$\mathfrak{b}_1$}}
\put(9,0.7){\makebox(0,0){$\mathfrak{b}_2$}}
\put(5,4.5){\vector(1,-1){3}}
\put(4,4.5){\vector(-1,-1){3}}
\put(1.4,1.5){\vector(1,1){3}}
\put(2.1,3.5){\makebox(0,0){$\pi_1$}}
\put(3.3,2.5){\makebox(0,0){$\iota_1$}}
\put(7.2,3.5){\makebox(0,0){$\pi_{21}$}}
\put(12.9,5){\makebox(0,0){$\mathfrak{b}$}}
\put(17,0.7){\makebox(0,0){$\mathfrak{b}_3$}}
\put(12.5,4.5){\vector(-1,-1){3}}
\put(13.4,4.5){\vector(1,-1){3}}
\put(16,1.5){\vector(-1,1){3}}
\put(10.3,3.5){\makebox(0,0){$\pi_{23}$}}
\put(15.3,3.5){\makebox(0,0){$\pi_3$}}
\put(14.3,2.5){\makebox(0,0){$\iota_3$}}
\end{picture}
\end{center}
where $\pi_1, \pi_{21}, \pi_{23}, \pi_3 $ are the
projections onto the ranges of $R_1, R_{21}, R_{23}$, and
$R_3$ respectively, according to the splittings
$\mathfrak{b} =
\mathfrak{b}_1 \oplus \mathfrak{b}_2 =
\mathfrak{b}_2\oplus\mathfrak{b}_3$, and $\iota_1:
\mathfrak{b}_1 \hookrightarrow \mathfrak{b}$, $\iota_3:
\mathfrak{b}_3 \hookrightarrow \mathfrak{b}$ are the
inclusions.

Then one has:
\begin{itemize}
\item[{\rm \textbf{(i)}}] If $\mathfrak{b}_2^\circ$ is a
Banach Lie subalgebra of $\mathfrak{b}^\ast$, then 
$\Phi_{31}: = \pi_3 \circ \iota_1 : (\mathfrak{b}_1,
\{\,,\}^{\rm coind}_{1}) \rightarrow (\mathfrak{b}_3,
\{\,,\}^{\rm coind}_{3})$ and
$\Phi_{13} : = \pi_1 \circ \iota_3 : (\mathfrak{b}_3,
\{\,,\}^{\rm coind}_{3}) \rightarrow (\mathfrak{b}_1,
\{\,,\}^{\rm coind}_{1})$ are mutually inverse linear
Poisson isomorphisms. 
\item[{\rm \textbf{(ii)}}] If $\mathfrak{b}_1^\circ $
and $\mathfrak{b}_3^\circ$ are Banach Lie subalgebras
of $\mathfrak{b}^\ast$, then $\mathfrak{b}_2 $ has two
coinduced Banach Lie-Poisson brackets 
$\{\,,\}^{\rm coind}_{21}$ and
$\{\,,\}^{\rm coind}_{23}$ which are not isomorphic in
general.
\end{itemize}
\end{proposition}

\noindent \textbf{Proof.} {\rm \textbf{(i)}} Since
$\mathfrak{b}_2^\circ = (\operatorname{im}R_{21})^\circ =
\operatorname{im}R_1 ^\ast$ (see \eqref{b star spaces
relation}) is a Banach Lie subalgebra of
$\mathfrak{b}^\ast$ it follows that $R_1$ coinduces a
Banach Lie-Poisson bracket $\{\,,\}^{\rm coind}_{1}$ on
$\mathfrak{b}_1$. Similarly, the relation $\mathfrak{b}_2
^ \circ = (\operatorname{im}R_{23})^ \circ =
\operatorname{im}R_3^\ast$ implies that $R_3$ coinduces a
Banach Lie-Poisson bracket $\{\,,\}^{\rm coind}_{3}$ on
$\mathfrak{b}_3$. 

Let us notice that  $\Phi_{31} \circ \Phi_{13} = \pi_3
\circ \iota_1
\circ \pi_1 \circ \iota_3 = \pi_3 \circ R_1 \circ \iota_3 
= \pi_3 \circ (\operatorname{id}_{\mathfrak{b}}  - R_{21})
\circ \iota_3 = \pi_3 \circ \iota_3  - \pi_3 \circ R_{21}
\circ \iota_3  = \operatorname{id}_{\mathfrak{b}_3}$ since
$\pi_3 \circ R_{21} = 0 $. One proves similarly  that 
$\Phi_{13} \circ \Phi_{31} =
\operatorname{id}_{\mathfrak{b}_1}$.
 
From $\ker \pi_1 = \ker \pi_3  = \mathfrak{b}_2 $ and $b -
(\iota_3 \circ \pi_3)(b) \in \ker \pi_3 $ for any $b \in
\mathfrak{b}$, it follows that $\pi_1 \circ \iota_3 \circ
\pi_ 3 = \pi_1 $. Therefore, if $f,g  \in 
C^{\infty}(\mathfrak{b}_1)$ we get from \eqref{coinduced
bracket with projector} and the fact that $\pi_1 :
\mathfrak{b}\rightarrow \mathfrak{b}_1 $ is a Poisson map
\begin{align*}
\{f \circ \Phi_{13}, g \circ \Phi_{13} \}^{\rm coind}_3 
&= \{f \circ \pi_1 \circ \iota_3, g \circ \pi_1 \circ
\iota_3 \}^{\rm coind}_3 \\
&= \{f \circ \pi_1 \circ \iota_3 \circ \pi_3, g \circ \pi_1
\circ \iota_3 \circ \pi_3\} \circ \iota_3\\
&= \{f \circ \pi_1, g \circ \pi_1 \} \circ \iota_3
= \{f, g\}^{\rm coind}_1 \circ \pi_1 \circ \iota_3 
= \{f, g\}^{\rm coind}_1 \circ \Phi_{13}.
\end{align*}
It is shown in a similar way that  $\Phi_{31}:
\mathfrak{b}_1 \rightarrow \mathfrak{b}_3 $ is a Poisson
map.

\textbf{(ii)} By \eqref{b spaces relation} we have 
$\mathfrak{b}_1 ^\circ = \operatorname{im}R_{21} ^\ast $
and $\mathfrak{b}_3 ^ \circ  = \operatorname{im}R_{23}
^\ast$ which, by hypothesis, are Banach Lie subalgebras of
$\mathfrak{b}^\ast$. Therefore, $\pi_{21} $ and $\pi_{23}$
coinduce Poisson brackets $\{\,,\}^{\rm coind}_{21}$ and
$\{\,,\}^{\rm coind}_{23}$ on $\mathfrak{b}_2 $. \quad
$\blacksquare$

\section{Symplectic induction}
\label{section: symplectic induction}

The goal of this section is to present the theory of symplectic 
induction on weak symplectic Banach manifolds. In the process we
shall define the momentum map in this setting, establish some of
its elementary properties, and give examples relevant to the
subsequent developments in this paper. 

\paragraph{Weak symplectic manifolds.} In infinite dimensions there
are two possible generalizations of the notion of a symplectic
manifold.

\begin{definition}
Let $P $ be a Banach manifold and $\omega$ a two-form. Then  $\omega$
is said to be {\bfi weakly nondegenerate\/} if for every $p \in P $
the map $v_p \in T_p P \mapsto \omega(p)(v_p, \cdot ) \in T ^\ast_pP
$ is injective. If, in addition, this map is also surjective, then
the form $\omega$ is called {\bfi strongly nondegenerate}. The form
$\omega$ is called a {\bfi weak} or {\bfi  strong symplectic
form\/} if, in addition, $\mathbf{d}\omega = 0$, where $\mathbf{d}$
denotes  the exterior differential on forms. The pair $(P, \omega)$
is called a {\bfi weak\/} or {\bfi strong symplectic manifold\/},
respectively.
\end{definition}

If $P $ is finite dimensional this distinction does not occur since
every linear  injective map is also surjective. The typical example
of an infinite dimensional strongly nondegenerate Banach manifold
is a complex Hilbert space endowed with the symplectic form equal to
the imaginary part of the Hermitian inner product. Any strong
symplectic form is locally constant but weak symplectic forms are
not, in general. The usual Hamiltonian formalism extends to the
strong  symplectic case without any difficulties.

On the other hand, if $(P , \omega)$ is a weak symplectic manifold,
the equation $\mathbf{d}h = \omega(X_h, \cdot ) $ that would define 
the Hamiltonian vector field $X_h$ associated to the function $h \in
C ^{\infty}(P) $ cannot always be  solved for $X_h$. But if $X_h $
exists, it is necessarily unique. Denote by $C ^{\infty}_{
\omega}(P)$ the vector subspace of smooth functions that admit
Hamiltonian vector fields. If $f, h \in C ^{\infty}_{ \omega}(P)$
their {\bfi Poisson bracket\/} is defined by
\begin{equation}
\label{standard PB}
\{f, h\}_{ \omega} : = \omega(X_f, X_h).
\end{equation}
In general, it is not true that $C^{\infty}_{\omega}(P)$ is a
Poisson algebra since $f, h \in C^{\infty}_{\omega}(P)$ does not
necessarily imply that $\{f, h\} \in C^{\infty}_{\omega}(P)$.
However, if $f, g, h \in C^{\infty}_{\omega}(P)$ and we assume, in
addition, that $\{f,g\},\{g,h\},\{h, f\} \in C^{\infty}_{\omega}(P)$, 
$\mathbf{d} \omega = 0 $ the same proof as in finite dimensions
implies the Jacobi identity.

Note that if $f,g \in C^{\infty}_{\omega}(P)$
then the product $fg \in C^{\infty}_{\omega}(P)$. Indeed, the
Hamiltonian vector field defined by $fg $ exists because $X_{fg} =
fX_g + gX_f$ as an easy computation shows. Another useful property
is that the Poisson bracket $\{f, g\}(p)$ for $f, g \in
C^{\infty}_{\omega}(P)$ is completely determined by
$\mathbf{d}f(p) $ and $\mathbf{d}g(p)$. Indeed, this follows from 
the fact that if $\mathbf{d} f(p) = \mathbf{d}g(p) $ then $\omega(p)
(X_f(p), \cdot ) = \mathbf{d} f(p) = \mathbf{d}g(p) =
\omega(p)(X_g(p), \cdot )$ and weak nondegeneracy of $\omega$
implies then that  $X_f(p) = X_g(p)$. Using this remark one can
recover several standard statements about Hamiltonian vector fields
in the weak symplectic case. 

\paragraph{The weak symplectic manifold $(\ell^{\infty} \times
\ell^1, \omega)$.}
In this paper we shall often work with the weak symplectic manifold
$(\ell ^{\infty} \times \ell ^1, \omega) $, where $\ell^\infty $ is
the  Banach space of bounded real sequences whose norm is given by
\[
\|\mathbf{q}\|_\infty: = \sup_{k=0,1,\dots} |q_k|, \qquad
\mathbf{q}: =\{q_k\}_{k=0} ^{\infty} \in \ell ^{\infty},
\]
$\ell ^1 $ is the Banach space of absolutely convergent real
sequences whose norm is given by
\[
\| \mathbf{p}\|_1 : = \sum_{k=0} ^\infty |p_k|, \qquad \mathbf{p}:
=\{p_k\}_{k=0} ^{\infty} \in \ell ^1,
\]
the strongly nondegenerate duality pairing
\begin{equation}
\label{qp_pairing}
\langle \mathbf{q}, \mathbf{p} \rangle = \sum_{k=0}^\infty q_k p_k,
\quad \text{for} \quad \mathbf{q} \in \ell ^{\infty}, \quad  
\mathbf{p}\in
\ell ^1,
\end{equation}
establishes the Banach space isomorphism  $(\ell ^1) ^\ast = 
\ell^\infty $, and the weak symplectic form $\omega$ has the
expression
\begin{equation}
\label{symplectic form}
\omega((\mathbf{q}, \mathbf{p}), (\mathbf{q}',
\mathbf{p}')) = \langle \mathbf{q},\mathbf{p}' \rangle -
\langle \mathbf{q}', \mathbf{p}\rangle, \quad \text{for} \quad
\mathbf{q}, \mathbf{q}' \in \ell ^\infty, \quad 
\mathbf{p}, \mathbf{p}' \in \ell ^1.
\end{equation}
for $\mathbf{q}, \mathbf{q}' \in \ell ^\infty $ and
$\mathbf{p}, \mathbf{p}' \in \ell ^1$.

The differential form $\omega$ is conveniently written as
\begin{equation}
\label{symplectic form in coordinates}
\omega = \sum_{k=0}^\infty \mathbf{d}q_k \wedge \mathbf{d}p_k.
\end{equation}
in the coordinates $q_k, p_k$. Let us elaborate on the notation used
in  \eqref{symplectic form in coordinates}. If $\mathbf{p} =
\{p _k\}_{k=0} ^{\infty} \in \ell ^1$, denote by  $\{\partial/
\partial p_k\}_{k=0} ^{\infty}$ the basis of the tangent space
$T_{\mathbf{p}}\ell ^1$ corresponding to the standard Schauder basis
$\{|k\rangle\}_{k=0} ^{\infty}$ of $\ell ^1$. The same basis in
$\ell ^{\infty}$ has a different meaning: every element $\mathbf{a}:
= \{a _k\}_{k=0} ^{\infty} \in
\ell^{\infty}$ can be uniquely written as a \textit{weakly} convergent
series $\mathbf{a} = \sum_{k=0} ^{\infty} a _k |k \rangle$. With this
notion of basis in $\ell ^{\infty}$, given $\mathbf{q} \in \ell
^{\infty}$, the sequence $\{ \partial/ \partial q _k\}_{k=0}
^{\infty} $ denotes the basis of the tangent space $T_{ \mathbf{q}}
\ell ^{\infty}$ corresponding to $\{|k \rangle \}_{k=0} ^{\infty} $.
Thus, any smooth vector field $X $ on $\ell ^{\infty} \times \ell ^1$ 
is written as
\[
X( \mathbf{q}, \mathbf{p}) = \sum_{k=0} ^{\infty} \left(A
_k(\mathbf{q},
\mathbf{p}) 
\frac{\partial}{\partial q _k} + B _k(\mathbf{q}, \mathbf{p})
\frac{ \partial}{ \partial p _k}\right),
\]
where $\{A _k(\mathbf{q}, \mathbf{p})\}_{k=0} ^{\infty} \in
\ell^{\infty}$ and $\{B _k(\mathbf{q}, \mathbf{p})\}_{k=0} ^{\infty}
\in \ell ^1$. If $Y $ is another vector field whose coefficients are
$\{C _k(\mathbf{q}, \mathbf{p})\}_{k=0} ^{\infty} \in \ell^{\infty}$,
$\{ D_k(\mathbf{q}, \mathbf{p})\}_{k=0} ^{\infty} \in \ell ^1$,
employing the usual conventions for the exterior derivatives of
coordinate functions to represent elements in the corresponding dual
spaces, formula 
\eqref{symplectic form in coordinates} gives
\begin{align*}
&\left(\sum_{k=0}^\infty \mathbf{d}q_k \wedge \mathbf{d}p_k \right)
\left(X,Y \right)(\mathbf{q}, \mathbf{p}) 
 = \sum_{k=0}^\infty
(A_k(\mathbf{q}, \mathbf{p}) D_k(\mathbf{q}, \mathbf{p}) -
C_k(\mathbf{q}, \mathbf{p}) B _k(\mathbf{q}, \mathbf{p}))
\end{align*}
which coincides with \eqref{symplectic form}. It is  in this sense
that the writing in \eqref{symplectic form in coordinates} represents
the weak symplectic form \eqref{symplectic form}.

In this case we can determine explicitly the space $C ^{\infty}_{
\omega}( \ell ^{\infty} \times \ell ^1)$. To do this, we observe
that for any $h \in C ^{\infty}_{ \omega} (\ell ^{\infty} \times
\ell ^1)$ its partial derivatives $\partial h/ \partial \mathbf{q}
\in (\ell ^{\infty})^\ast$ and $\partial h/ \partial \mathbf{p} \in
(\ell^1)^\ast = \ell ^{\infty}$, respectively. Thus the
Hamiltonian vector field $X_h$ defined by the weak symplectic form 
\eqref{symplectic form in coordinates} and the function $h $ exists
if and only if $\partial h / \partial\mathbf{q} \in \ell ^1 \subset
(\ell ^1)^{\ast \ast} = (\ell ^{\infty}) ^\ast$.
Therefore,
\begin{equation}
\label{Hamiltonian functions ell}
C^{\infty}_{\omega}(\ell ^{\infty} \times \ell ^1) =
\{f\in C ^{\infty}(\ell ^{\infty} \times \ell^1) \mid \{
\partial h /\partial q _k\}_{k=0} ^{\infty} \in \ell ^1\},
\end{equation}
and the Hamiltonian vector field defined by $h \in
C^{\infty}_{\omega}(\ell ^{\infty} \times \ell^1)$ has the expression
\begin{equation}
\label{v f ell}
X _h(\mathbf{q} , \mathbf{p}) = \frac{\partial h}{\partial
p_k}\frac{ \partial}{\partial q _k} - \frac{\partial h}{\partial
q_k}\frac{ \partial}{\partial p_k}.
\end{equation}
The canonical Poisson bracket of $f, h \in
C^{\infty}_{\omega}(\ell ^{\infty}  \times  \ell ^1)$ makes sense
and is given by 
\begin{equation}
\label{canonical Poisson bracket on ell}
\{f, g\}_\omega( \mathbf{q}, \mathbf{p}) = \sum_{k=0} ^{\infty}
\left(\frac{\partial f}{\partial q_k}\frac{\partial g}{\partial p_k} -
\frac{\partial g}{\partial q_k} \frac{\partial f}{\partial p_k}\right).
\end{equation}

For the weak symplectic Banach vector space  $(\ell^{\infty}\times
\ell^1, \omega) $, a direct computation shows that the Poisson
bracket  of any two functions from the set
\begin{align*}
&\left\{ f \in C ^{\infty}_{\omega}(\ell^{\infty} \times \ell^1)
\,\Big{|}\,
\left\{\sum_{j=0} ^{\infty}\frac{\partial^2 f}{ \partial q_i \partial q_j}
q_j' \right\}_{i=0}^{\infty} \in \ell^1,   \right. \nonumber \\ 
& \qquad 
\left.  \left\{\sum_{j=0} ^{\infty}\frac{\partial^2 f}{ \partial q_i
\partial p_j}p_j' \right\}_{i=0}^{\infty} \in \ell^1
 \text{~for all~} \{q_j'\}_{j=0}^{\infty} \in \ell ^{\infty}, 
\{p_j'\}_{j=0}^{\infty} \in \ell ^1 \right\}.
\end{align*}
is again in $C^{\infty}_{\omega}(\ell ^{\infty} \times \ell ^1)$.

\paragraph{Momentum maps on weak symplectic manifolds.} Throughout
this section, $G$ denotes a Banach Lie group and $\mathfrak{g}$ its
Lie algebra. We shall assume that $\mathfrak{g}$ admits a predual
$\mathfrak{g}_\ast$ and that the coadjoint action of $G $ on the
dual space $\mathfrak{g}^\ast$ leaves $\mathfrak{g}_\ast
\subset \mathfrak{g}^\ast$ invariant, that is
$\operatorname{Ad}^\ast_g \mathfrak{g}_\ast \subset
\mathfrak{g}_\ast $, for any $g \in G $. Recall from \cite{OR} that
$\mathfrak{g}_\ast $ is a Banach Lie-Poisson space (whose bracket is
hence given by \eqref{general LP}).

\begin{definition}
\label{def momentum map weak}
Let $(P, \omega)$ be a weak symplectic manifold and $G $ a Banach Lie
group satisfying the conditions above. A smooth map $\mathbf{J}: P
\rightarrow
\mathfrak{g}_\ast$ is a {\bfi momentum map\/} if whenever $\varphi,
\psi$ are locally defined smooth functions on $\mathfrak{g}_\ast $
such that $\varphi \circ \mathbf{J}, \psi \circ \mathbf{J}$ are
locally defined elements of $C ^{\infty}_{\omega}(P)$, we have $\{
\varphi \circ \mathbf{J}, \psi \circ \mathbf{J} \}_{ \omega} =
\{\varphi, \psi\} \circ \mathbf{J}$. Here
$\{\cdot ,\cdot \}_{\omega}$ denotes the Poisson bracket on
functions in $C ^{\infty}_{ \omega}(P)$ and
$\{ \cdot , \cdot \} $ is the Lie-Poisson bracket on
$\mathfrak{g}_\ast$.
\end{definition}

Momentum maps usually appear by the following construction. 

\begin{proposition} 
\label{prop: mom}
Let $\Phi:G \times P \rightarrow P $ be a smooth
symplectic action of the Banach Lie group $G $ on the weak
symplectic  Banach manifold $(P,\omega)$. Assume that the smooth map
$\mathbf{J}: P \rightarrow \mathfrak{g}_\ast $ is $G $-equivariant
and is such that for all $z \in \mathfrak{g}$ we have
$z \circ \mathbf{J} \in C ^{\infty}_{\omega}(P)$ and $z_P = X_{z
\circ \mathbf{J}}$, where $z_P(p): = \left.\frac{d}{dt}\right|_{t=0} \Phi\left(\exp(tz), 
p \right)$ denote the infinitesimal generator of the action. Then
$\mathbf{J}$ is a momentum map.
\end{proposition}

\noindent \textbf{Proof.} We proceed as in finite dimensions (see,
e.g., \cite{MR}). First note that if $\varphi$ is a smooth locally
defined function on $\mathfrak{g}_\ast$ and $p
\in P $, denoting $y: = \mathbf{d}\varphi(\mathbf{J}(p)) \in
\mathfrak{g}$, we have $\mathbf{d}(\varphi \circ \mathbf{J})(p) =
\mathbf{d}(y \circ \mathbf{J} )(p)$. But the Poisson bracket
evaluated at $p
$ depends only on the first derivatives of the functions  at $p $
which means that if
$\psi $ is another locally defined function on $\mathfrak{g}_\ast $
and $z: = \mathbf{d}\psi(\mathbf{J}(p)) \in \mathfrak{g}$ we have
\[
\{\varphi \circ \mathbf{J}, \psi\circ \mathbf{J}\}_{ \omega}(p) =
\{y \circ \mathbf{J}, z \circ \mathbf{J} \}_{ \omega}(p).
\]
On the other hand, the derivative at $g=e$ of the equivariance
identity $\mathbf{J}(g \cdot p) = \operatorname{Ad}^\ast _{g ^{-1}}
\mathbf{J}(p) $ for any $g \in G$ and $p \in P$ yields the relation
$T_p \mathbf{J} (z_P(p)) = -\operatorname{ad}^\ast_z \mathbf{J}(p)$
for any $z \in \mathfrak{g}$. Therefore, by \eqref{general LP} we get
\begin{align*}
(\{ \varphi, \psi\} \circ \mathbf{J})(p) 
&= \left\langle [ \mathbf{d}\varphi(\mathbf{J}(p)),
\mathbf{d}\psi(\mathbf{J}(p))], \mathbf{J}(p) \right\rangle
= \left\langle [y, z], \mathbf{J}(p) \right\rangle 
= \left\langle z, \operatorname{ad}^\ast_y \mathbf{J}(p)
\right\rangle \\ 
& = - \left\langle z, T_p\mathbf{J}(y_P(p)) \right\rangle
= - \mathbf{d}(z \circ \mathbf{J})(p) (y_P(p)) 
= - \mathbf{d}(z \circ \mathbf{J})(p)\left(X_{y \circ \mathbf{J}}(p)
\right) \\ 
& = \{y \circ \mathbf{J}, z \circ \mathbf{J}\}_{ \omega}(p)
\end{align*}
which shows that $\{\varphi \circ \mathbf{J}, \psi\circ
\mathbf{J}\}_{\omega} = \{
\varphi, \psi\} \circ \mathbf{J}$ and hence $\mathbf{J} :P
\rightarrow
\mathfrak{g}_\ast $ is a momentum map by Definition 
\ref{def momentum map weak}.
\quad $\blacksquare$
\medskip

Note that $C ^{\infty}_ \omega(P) $ is invariant by the $G $-action.
Indeed, the Hamiltonian vector field of the smooth function $f \circ
\Phi_g $ for $f \in C ^{\infty}_ \omega(P)$ exists and equals
$\Phi_g^\ast X_f$, where $\Phi_g: P \rightarrow P $ denotes the $G
$-action on $P $. Similarly, for any $z \in \mathfrak{g}$, the
Hamiltonian vector field of $\mathbf{d}f(z_P)$ exists and equals $[
z_P,X_f]$.

Propositions 7.3 and
7.4 in \cite{OR} show that if the coadjoint isotropy subgroup of
$\rho \in \mathfrak{g}_\ast$ is a closed Lie subgroup of $G $, the
coadjoint orbit is a weak symplectic manifold and the inclusion is a
momentum map in the sense of Definition \ref{def momentum map weak}.
We shall give study other momentum maps in subsequent sections.

\paragraph{The symplectic induced space.} Symplectic 
induction is a technique that associates to a given
Hamiltonian $H$-space a Hamiltonian $G$-space whenever
$H$ is a Lie subgroup of the Lie group $G$; see
\cite{DuElTu1992, GuSt1982a, GuSt1983, KaKoSt1978,
Weinstein1978, Za1986} for various versions of this
construction and several applications. We shall review this
method below in the category of Banach manifolds and shall 
impose also  certain splitting assumptions that are
satisfied in the example studied later. 

Let $G$ be a Banach Lie group with Banach Lie algebra
$\mathfrak{g}$. Let $H$ be a closed
Banach Lie subgroup of $G$ with Banach Lie algebra
$\mathfrak{h}$. Assume that both $\mathfrak{g}$ and
$\mathfrak{h}$ admit preduals $\mathfrak{g}_\ast$ and
$\mathfrak{h}_\ast$, which are invariant under the coadjoint 
actions of $G$ and $H$, respectively (see \cite{OR} for
various consequences of this assumption).
Throughout this section
we shall make the following hypotheses:
\begin{itemize}
\item $\mathfrak{h}_\ast \subset \mathfrak{g}_\ast$,
\item there is an $\operatorname{Ad}^\ast_H$-invariant
splitting 
\begin{equation}
\label{predual splitting}
\mathfrak{g}_\ast = \mathfrak{h}_\ast \oplus
\mathfrak{h}_\ast ^\perp,
\end{equation}
 where $\mathfrak{h}_\ast ^\perp$ is a
Banach $\operatorname{Ad}^\ast_H$-invariant subspace of
$\mathfrak{g}_\ast$, which means that
$\operatorname{Ad}^\ast_h \mathfrak{h}_\ast ^\perp \subset
\mathfrak{h}_\ast ^\perp$ for any $h \in H $, where
$\operatorname{Ad}^\ast: G \rightarrow
\operatorname{Aut}( \mathfrak{g}_\ast) $ is the $G
$-coadjoint action,
\item $\left(\mathfrak{h}_\ast ^\perp \right) ^\circ  =
\mathfrak{h}$, where $\left(\mathfrak{h}_\ast ^\perp \right) ^\circ$
is the annihilator of $\mathfrak{h}_\ast ^\perp$,
\item the Banach Lie group $H $ acts symplectically on the weak
symplectic Banach manifold $(P, \omega)$ and there  is a $H
$-equivariant  map $\mathbf{J}^H_P: P \rightarrow \mathfrak{h}_\ast $
satisfying the hypothesis of  Proposition \ref{prop: mom} for the
Lie group $H $ and hence $\mathbf{J}^H_P$ is a momentum map.
\end{itemize}
Dualizing the splitting \eqref{predual splitting}, we get an
$\operatorname{Ad}_H$-invariant splitting 
\begin{equation}
\label{LA splitting}
\mathfrak{g} =
\mathfrak{h} \oplus \mathfrak{h}^\perp,
\end{equation} 
where $\mathfrak{h}^\perp := ( \mathfrak{h}_\ast ) ^ \circ $ is
the annihilator of the Banach Lie-Poisson space
$\mathfrak{h}_\ast $. 
\medskip

The induction method produces a Hamiltonian
$G$-space by constructing a reduced manifold in the
following way. Form the product $P \times G \times
\mathfrak{g}_\ast$ of weak symplectic manifolds,
where $G \times \mathfrak{g}_\ast$ has the weak symplectic 
form
\begin{align}
\label{left trivialized canonical two form}
\omega_L(g, \tilde{\rho})\left((u_g, \tilde{\mu}), (v_g,
\tilde{\nu}) \right) 
  &= \langle \tilde{\nu}, T_gL_{g^{-1}}u_g \rangle  -
\langle \tilde{\mu}, T_gL_{g^{-1}}v_g \rangle 
\nonumber\\
& \qquad
+ \langle \tilde{\rho},
[T_gL_{g^{-1}}u_g , T_gL_{g^{-1}}v_g ] \rangle,
\end{align}
for $g \in G$, $u_g, v_g \in T_gG$, and $\tilde{\rho},
\tilde{\mu}, \tilde{\nu} \in \mathfrak{g}_\ast$. This
formula was introduced in \cite{OR} and it looks formally
the same as the left trivialized canonical symplectic form
on the cotangent bundle of a finite dimensional Lie group 
(see \cite{A-M}, \S 4.4, Proposition 4.4.1). From 
\eqref{left trivialized canonical two form} it follows that 
\[
C ^{\infty}_{\omega_L}\left(G \times \mathfrak{g}_\ast \right)
= \{k \in C^{\infty}(G \times \mathfrak{g}_\ast ) \mid T_e^\ast L_g
d_1k(g, \tilde{\rho}) \in \mathfrak{g}_\ast \},
\]
where $d_1k(g, \rho) \in T^\ast_g G$ and $d_2k(g, \rho) \in
(\mathfrak{g}_\ast)^\ast = \mathfrak{g}$ are the first and second
partial derivatives of $k$. If $k \in C ^{\infty}_{ \omega_L}(G
\times \mathfrak{g}_\ast) $, the Hamiltonian vector field $X_k \in
\mathfrak{X}\left(G \times \mathfrak{g}_\ast \right)$ has the
expression
\begin{equation}
\label{Hamiltonian vector field left trivialized}
X_k(g, \tilde{\rho}) = \left(T_eL_g d_2k(g, \tilde{\rho}), 
\operatorname{ad}^\ast_{d_2k(g,\tilde{\rho})} \tilde{\rho} -
 T_e^\ast L_g d_1k(g, \tilde{\rho})\right).
\end{equation}
Therefore the canonical Poisson bracket of $f, k \in C ^{\infty}_{
\omega_L} \left(G \times \mathfrak{g}_\ast \right)$ equals
\begin{multline}
\label{Poisson bracket in left trivialization}
\{f, k\}(g, \rho) = \left\langle d_1f(g, \rho), T_eL_g d_2k(g, \rho)
\right\rangle - \left\langle d_1k(g, \rho), T_eL_g d_2f(g, \rho)
\right\rangle \\ 
- \left\langle \rho, \left[d_2 f(g, \rho), d_2 k(g,
\rho)\right] \right\rangle.
\end{multline}
The left $G$-action on $G \times \mathfrak{g}_\ast$ given by $g'
\cdot (g, \rho) : = (g'g, \rho)$ induces the momentum map $(g, \rho)
\mapsto \operatorname{Ad}^\ast_{g^{-1}} \rho$ which is $G
$-equivariant.

The weak symplectic  form $\omega\oplus \omega_L \in
\Omega^2(P \times G \times \mathfrak{g}_\ast)$ is defined by
\begin{align}
\label{sum symplectic}
&(\omega\oplus \omega_L)(p, g, \tilde{\rho})\left((a_p, T_e
L_g \tilde{x}, \tilde{\mu}), (b_p, T_e L_g \tilde{y},
\tilde{\nu})\right) =    \nonumber \\
& \qquad \qquad \qquad 
\omega(p)(a_p, b _p) +  \langle \tilde{\nu},
\tilde{x} \rangle  - \langle \tilde{\mu}, \tilde{y} \rangle 
+ \langle \tilde{\rho}, [\tilde{x}, \tilde{y}] \rangle,
\end{align}
where $p  \in P$, $g \in G $, $\tilde{\rho}, \tilde{\mu},
\tilde{\nu} \in
\mathfrak{g}_\ast$, $\tilde{x}, \tilde{y} \in
\mathfrak{g}$, and
$a_p , b _p\in T _p P $.

The Banach Lie group $H$ acts on $P \times G \times
\mathfrak{g}_\ast$  by 
\begin{equation}
\label{h action}
h \cdot (p, g, \tilde{\rho} ) : = ( h \cdot p, gh ^{-1},
\operatorname{Ad}^\ast_{h ^{-1}} \tilde{\rho}).
\end{equation}
The infinitesimal generator of this action defined by $z \in
\mathfrak{h}$ equals
\[
z_{P \times G \times \mathfrak{g}_\ast}(p,g,\tilde{ \rho}) 
= \left(z_P(p), - T _e L_g z , - \operatorname{ad}^\ast_z
\tilde{\rho}\right)
\]
which, by \eqref{Hamiltonian vector field left trivialized} and the
assumption of the existence of a momentum map induced by the action
of $H $ on $P $, is a Hamiltonian vector field relative to the
function $z \circ \left(\mathbf{J}^H_P (p) - \Pi\tilde{ \rho}
\right)$,  where $\Pi: \mathfrak{g}_\ast \rightarrow
\mathfrak{h}_\ast$ is the projection defined by the
splitting $\mathfrak{g}_\ast = \mathfrak{h}_\ast \oplus
\mathfrak{h}_\ast ^\perp$.
Therefore, the action \eqref{h action} admits the equivariant
momentum map $J^H_{P \times G \times \mathfrak{g}_\ast}: P \times G
\times \mathfrak{g}_\ast \rightarrow \mathfrak{h}_\ast$ given by
\begin{equation}
\label{big mom map}
\mathbf{J}^H_{P \times G \times \mathfrak{g}_\ast}(p, g,
\tilde{
\rho}) = \mathbf{J}^H_P(p) - \Pi \tilde{\rho}.
\end{equation}

The $H$-action on $P \times G \times \mathfrak{g}_\ast$ is free and
proper because $H$ is a closed Banach Lie subgroup of $G$.
Therefore its restriction to the closed invariant subset
$(\mathbf{J}^H_{P \times G \times \mathfrak{g}_\ast}) ^{-1}(0)$ is
also free and proper. Let us assume at this point that $0 $ is a
regular value and hence that  $(\mathbf{J}^H_{P \times G \times
\mathfrak{g}_\ast}) ^{-1}(0)$ is a submanifold. In concrete
applications, such as gravity or Yang-Mills theory, the
proof of the regularity  of $0 $ is usually achieved by appealing to
elliptic operator theory. With the assumption that $0 $ is a regular
value and that for each $(p, g, \tilde{ \rho}) \in (\mathbf{J}^H_{P
\times G \times \mathfrak{g}_\ast}) ^{-1}(0)$ the map $h \in H
\mapsto h
\cdot  (p, g, \tilde{\rho} ) : = ( h \cdot p, gh ^{-1},
\operatorname{Ad}^\ast_{h ^{-1}} \tilde{\rho}) \in 
(\mathbf{J}^H_{P \times G \times \mathfrak{g}_\ast}) ^{-1}(0)$ is an
immersion, it follows that the quotient topological space $M:= 
(\mathbf{J}^H_{P \times G \times \mathfrak{g}_\ast}) ^{-1}(0)/H$
carries a unique smooth manifold structure relative to which the
quotient projection is a submersion. This underlying manifold
topology is that of the quotient topological space and it is
Hausdorff (see \cite{Bou1972}, Chapter III, \S1, Proposition 10
for a proof of these statements). Once these topological conditions
are satisfied, a technical lemma (stating that the double
symplectic  orthogonal of a closed subspace in a weak
symplectic Banach space is equal to the original subspace)  allows
one to extend the original proof of the reduction theorem in finite
dimensions (see \cite{mwr}) to the case of weak
symplectic Banach manifolds. We shall not dwell here on these
technicalities because in the example of interest to us, treated
later, the reduction process will be carried out by hand without any
appeal to general theorems. Summarizing, we can form the 
{\bfi induced space} $(M,\Omega_M)$ which is a smooth Hausdorff weak
symplectic Banach manifold, where $\Omega_M$ is the reduced
symplectic form on $\mathbf{J}^H_{P \times G \times
\mathfrak{g}_\ast})^{-1}(0)/H$.

Now note that if we denote $\tilde{ \rho} = \rho + \rho^\perp
\in \mathfrak{h}_\ast \oplus \mathfrak{h}_\ast ^\perp$ we get
\begin{align*}
(\mathbf{J}^H_{P \times G \times \mathfrak{g}_\ast})
^{-1}(0)  & = \left\{(p, g, \tilde{ \rho}) \in P \times G
\times \mathfrak{g}_\ast \mid \mathbf{J}^H_P(p) = \Pi
\tilde{\rho} \right\} \\
& = G \times \left\{(p, \rho) \in P \times
\mathfrak{h}_\ast \mid \mathbf{J}^H_P(p) = \rho \right\}
\times \mathfrak{h}_\ast ^\perp \\
& \cong G \times P \times \mathfrak{h}_\ast ^\perp,
\end{align*}
where the $H$-equivariant diffeomorphism in the last line
is given by 
\[
(p, \rho) \in \left\{(p, \rho) \in P \times
\mathfrak{h}_\ast \mid \mathbf{J}^H_P(p) = \rho \right\}
\longmapsto  p \in P.
\]

Therefore the weak symplectic Banach manifold $M
= (\mathbf{J}^H_{P \times G \times \mathfrak{g}_\ast})
^{-1}(0)/H$ is diffeomorphic to the fiber
bundle $G \times _H (P \times \mathfrak{h}_\ast ^\perp)
\rightarrow G/H$ associated to $G \rightarrow G/H$.

\paragraph{The weak symplectic form on the induced space.}
Let us denote by $\pi_0: G \times P \times
\mathfrak{h}_\ast ^\perp \rightarrow G \times _H (P \times
\mathfrak{h}_\ast ^\perp)$ the projection onto the  $H$-orbit
space. The next statement gives the weak symplectic form on $M$.

\begin{proposition}
\label{induction theorem}
The associated fiber bundle $G \times _H
(P \times \mathfrak{h}_\ast ^\perp) \rightarrow G/H$ has a weak
symplectic form $\Omega$  given by
\begin{align}
&\Omega(\pi_0(g, p, \rho^\perp)) \left(T_{(g, p, \rho^\perp)} \pi_0
(T_e L_g (x + x^\perp), a_p, \mu^\perp), T_{(g, p, \rho^\perp)} \pi_0
(T_e L_g (y + y^\perp), b_p, \nu^\perp) \right) \nonumber\\
& \qquad = \omega(p)(a_p, b_p) + \left\langle
T_p\mathbf{J}^H_P(b_p), x \right\rangle + 
\left\langle \nu^\perp, x^\perp \right\rangle
- \left\langle T_p\mathbf{J}^H_P(a_p), y \right\rangle - 
\left\langle \mu^\perp, y^\perp \right\rangle \nonumber\\
& \qquad \qquad \quad  +
\left\langle \mathbf{J}^H_P (p), [x,y]
\right\rangle + \left\langle \rho^\perp, [x^\perp, y] + [x, y^\perp]
\right\rangle + \left\langle \mathbf{J}^H_P(p) + \rho^\perp,
[x^\perp, y^\perp] \right\rangle 
\label{quotient symplectic form}\\
&\qquad = \omega(p)(a_p - x_P(p), b_p - y_P(p)) +
\left\langle T_p
\mathbf{J}^H_P(b_p - y_P(p)), 2x \right\rangle +
\left\langle\nu^\perp + \operatorname{ad}^\ast_y \rho^\perp, x^\perp
\right\rangle \nonumber \\
& \qquad \quad \qquad 
- \left\langle T_p \mathbf{J}^H_P(a_p - x_P(p)),
2y \right\rangle - \left\langle \mu^\perp +
\operatorname{ad}^\ast_x \rho^\perp, y^\perp \right\rangle + 
\left\langle \mathbf{J}^H_P(p), [2x, 2y] \right\rangle
\nonumber \\ 
& \qquad \qquad \quad 
+ \left\langle\rho^\perp, [x^\perp, 2y] + [2x,
y^\perp] \right\rangle + \left\langle \mathbf{J}^H_P(p) 
+ \rho^\perp, [x^\perp, y^\perp] \right\rangle
\label{quotient symplectic form transversal}
\end{align}
for $g \in G $, $p \in P $, $\rho^\perp, \mu^\perp, \nu^\perp \in
\mathfrak{h}_\ast ^\perp$, $x, y \in
\mathfrak{h}$, $x^\perp, y^\perp \in \mathfrak{h}^\perp$, and $a_p,
b_p\in T_p P$. The second expression uses only tangent
vectors of the form 
\[
\left(a_p - x_P(p), T_e L_g (2x +  x^\perp), 
\mu^\perp + \operatorname{ad}^\ast_x \rho^\perp \right) 
\]
 which are transversal to the $H $-orbits in the
zero level set of the momentum map and hence represent
the tangent space $T_{\pi_0(g, p, \rho^\perp)}M$ to the reduced
manifold $M$. 
\end{proposition}

\noindent \textbf{Proof.} We begin with the proof 
\eqref{quotient symplectic form}. Let $i_0: G \times P \times
\mathfrak{h}_\ast ^\perp \hookrightarrow  P
\times G \times \mathfrak{g}_\ast$ be the inclusion
$i_0(g, p, \rho^\perp): = (p, g, \mathbf{J}^H_P(p) + \rho^\perp)$.
For $p \in P $, $\rho^\perp, \mu^\perp, \nu^\perp\in
\mathfrak{h}_\ast^\perp$,
$g \in G $, $\tilde{x} = x + x^\perp, \tilde{y} = y + y^\perp \in
\mathfrak{g}$, $x, y \in \mathfrak{h}$, $x^\perp, y^\perp \in
\mathfrak{h}^\perp$, and $a_p, b_p \in T_p P$, the reduction
theorem and \eqref{sum symplectic} 
give
\begin{align*}
&\Omega(\pi_0(g, p, \rho^\perp)) \left(T_{(g, p, \rho^\perp)} \pi_0
(T_e L_g
\tilde{x}, a_p, \mu^\perp), T_{(g, p, \rho^\perp)} \pi_0 (T_e L_g
\tilde{y}, b_p, \nu^\perp) \right) \nonumber\\
&\quad = i_0 ^\ast (\omega\oplus \omega_L)(p, g, \rho^\perp)
\left((a_p, T_e L_g \tilde{x}, \mu^\perp), (b_p, T_e L_g
\tilde{y}, \nu^\perp)\right)   \\
& \quad = (\omega\oplus \omega_L)(p, g,
\mathbf{J}^H_P(p) + \rho^\perp) \left((a_p, T_e L_g \tilde{x},
T_p \mathbf{J}^H_P(a_p) + \mu^\perp), (b_p, T_e L_g \tilde{y},
T_p \mathbf{J}^H_P(b_p) + \nu^\perp) \right)
\\ 
& \quad  =  \omega(p)(a_p, b_p) +  \left\langle T_p
\mathbf{J}^H_P(b_p) + \nu^\perp, x + x^\perp \right\rangle  -
\left\langle T_p \mathbf{J}^H_P(a_p) + \mu^\perp, y + y^\perp
\right\rangle
\nonumber \\
& \qquad \qquad   
+ \left\langle \mathbf{J}^H_P(p) + \rho^\perp, [x + x^\perp, y +
y^\perp]
\right\rangle.
\end{align*}
Since $[x+x^\perp, y+y^\perp] = [x,y] + [x^\perp, y] + [x, y^\perp] +
[x^\perp, y^\perp]$, $[x,y] \in
\mathfrak{h} = (\mathfrak{h}_\ast ^\perp)^ \circ$, $[x^\perp, y] +
[x, y^\perp] \in \mathfrak{h}^\perp = (\mathfrak{h}_\ast )^\circ$
(because the splitting $\mathfrak{g} =
\mathfrak{h} \oplus \mathfrak{h}^\perp$ is 
$\operatorname{Ad}^\ast_H $-invariant), 
$\rho^\perp \in \mathfrak{h}_\ast ^\perp$, and $\mathbf{J}^H_P(p)
\in
\mathfrak{h}_\ast$, the last term becomes
\begin{align*}
\left\langle \mathbf{J}^H_P(p) + \rho^\perp, [x + x^\perp, y + y^\perp]
\right\rangle = &\left\langle \mathbf{J}^H_P (p), [x,y]
\right\rangle + \left\langle \rho^\perp, [x^\perp, y] + [x, y^\perp]
\right\rangle \\
& + \left\langle \mathbf{J}^H_P(p) + \rho^\perp,
[x^\perp, y^\perp] \right\rangle.
\end{align*} 
Since $T_p \mathbf{J}^H_P(b_p) \in \mathfrak{h}_\ast$,
$\nu^\perp \in \mathfrak{h}_\ast ^\perp$, $x \in \mathfrak{h} =
(\mathfrak{h}_\ast ^\perp)^ \circ$, and $x^\perp \in
\mathfrak{h}^\perp = (\mathfrak{h}_\ast )^\circ$, the second
term becomes
\[
\left\langle T_p
\mathbf{J}^H_P(b_p) + \nu^\perp, x + x^\perp \right\rangle
= \left\langle T_p\mathbf{J}^H_P(b_p), x \right\rangle + 
\left\langle \nu^\perp, x^\perp \right\rangle.
\]
Similarly, the third term is
\[
\left\langle T_p \mathbf{J}^H_P(a_p) + \mu^\perp, y + y^\perp
\right\rangle =  \left\langle T_p\mathbf{J}^H_P(a_p), y
\right\rangle + 
\left\langle \mu^\perp, y^\perp \right\rangle.
\]
Thus we get
\begin{align*}
&\Omega(\pi_0(g, p, \rho^\perp)) \left(T_{(g, p, \rho^\perp)} \pi_0
(T_e L_g (x + x^\perp), a_p, \mu^\perp), T_{(g, p, \rho^\perp)} \pi_0
(T_e L_g (y+y^\perp), b_p, \nu^\perp) \right) \\
& \qquad = \omega(p)(a_p, b_p) + \left\langle
T_p\mathbf{J}^H_P(b_p), x \right\rangle + 
\left\langle \nu^\perp, x^\perp \right\rangle
- \left\langle T_p\mathbf{J}^H_P(a_p), y \right\rangle - 
\left\langle \mu^\perp, y^\perp \right\rangle \\
& \qquad \qquad +
\left\langle \mathbf{J}^H_P (p), [x,y]
\right\rangle + \left\langle \rho^\perp, [x^\perp, y] + [x, y^\perp]
\right\rangle + \left\langle \mathbf{J}^H_P(p) + \rho^\perp,
[x^\perp, y^\perp] \right\rangle
\end{align*}
which proves \eqref{quotient symplectic form}.

We want to simplify this expression by taking  advantage
of the $H$-action on the zero level set of the momentum
map. For $x \in \mathfrak{h}$ we have by $H$-equivariance
of $\mathbf{J}^H_P$ and the $\operatorname{Ad}^\ast
_H$-invariance of the splitting $\mathfrak{g}_\ast =
\mathfrak{h}_\ast \oplus \mathfrak{h}_\ast ^\perp$
\begin{align*}
&x_{P \times G \times \mathfrak{g}_\ast}(p, g,
\mathbf{J}^H_P(p) + \rho^\perp )  
= \left.\frac{d}{dt}\right|_{t=0} \left( \exp t x
\cdot p, g \exp(-t x ), \operatorname{Ad}^\ast_{\exp
(-tx)} (\mathbf{J}^H_P(p) + \rho^\perp) \right) \\
&\qquad = \left.\frac{d}{dt}\right|_{t=0} \left( \exp t x
\cdot p, g \exp(-t x ), \mathbf{J}^H_P(\exp t x \cdot p) +
\operatorname{Ad}^\ast_{\exp (-tx)} \rho^\perp \right) \\
& \qquad = 
\left(x_P(p), -T_e L_g x,  T_p \mathbf{J}^H_P(x_P(p)) -
\operatorname{ad}^\ast_x \rho^\perp
\right).
\end{align*}
Now decompose
\begin{align*}
&\left(a_p, T_e L_g(x + x^\perp), T_p \mathbf{J}^H_P(a_p) +
\mu^\perp \right) \\
& \qquad = 
\left(x_P(p), -T_e L_g x, T_p \mathbf{J}^H_P(x_P(p)) -
\operatorname{ad}^\ast_x \rho^\perp \right) \\
& \qquad \qquad + 
\left(a_p - x_P(p), T_e L_g (2x +  x^\perp), T_p
\mathbf{J}^H_P(a_p - x_P(p)) + \mu^\perp +
\operatorname{ad}^\ast_x \rho^\perp \right).
\end{align*}
Since the form $\Omega$ does not depend on the first 
summand, this means that we can replace everywhere in 
\eqref{quotient symplectic form} $a_p$
by $a_p - x_P(p)$, $x$ by $2x$, and $\mu^\perp$ by $\mu^\perp +
\operatorname{ad}^\ast_x \rho^\perp$. Similarly, we can
replace $b_p$ by $b_p - y_P(p)$, $y $ by $2y$, and $\nu^\perp$
by $\nu^\perp + \operatorname{ad}^\ast_y \rho^\perp$. Thus
\eqref{quotient symplectic form} becomes
\begin{align*}
&\omega(p)(a_p - x_P(p), b_p - y_P(p)) + \left\langle T_p
\mathbf{J}^H_P(b_p - y_P(p)), 2x \right\rangle +
\left\langle\nu^\perp + \operatorname{ad}^\ast_y \rho^\perp, x^\perp
\right\rangle \\
& \qquad 
- \left\langle T_p \mathbf{J}^H_P(a_p - x_P(p)),
2y \right\rangle - \left\langle \mu^\perp +
\operatorname{ad}^\ast_x \rho^\perp, y^\perp \right\rangle + 
\left\langle \mathbf{J}^H_P(p), [2x, 2y] \right\rangle \\
& \qquad \qquad 
+ \left\langle\rho^\perp, [x^\perp, 2y] + [2x,
y^\perp] \right\rangle + \left\langle \mathbf{J}^H_P(p) 
+ \rho^\perp, [x^\perp, y^\perp] \right\rangle 
\end{align*}
which proves \eqref{quotient symplectic form
transversal}.
\quad $\blacksquare$

\paragraph{Remark.} If $H = G $, then one can verify directly
that the map $\Psi:G \times _H (P \times \{0\}) \rightarrow P$
given by $\Psi(\pi_0(g, p,0)): = g\cdot p $ is a diffeomorphism
between the weak symplectic manifolds $(G\times _H (P \times \{0\}),
\Omega)$ (the induced space) and
$(P,\omega)$ (the original manifold). 

\paragraph{The momentum map on the induced space.} Now we shall
construct a $G $-action on the induced space $\left(G \times _H
(P \times \mathfrak{h}_\ast ^\perp) , \omega \right)$ and a $G
$-equivariant momentum map $\mathbf{J}_M^G: G \times _H
(P \times \mathfrak{h}_\ast ^\perp) \rightarrow \mathfrak{g}_\ast$.

The Banach Lie group $G$ acts on $G \times P\times
\mathfrak{h}_\ast^\perp$ by $g' \cdot (g, p , \rho^\perp) : = (g'g,
p, \rho^\perp) $. This $G$-action commutes with the
$H$-action and so $G $ acts  on the induced space
$G\times _H (P \times \mathfrak{h}_\ast ^\perp)$ by $g' \cdot [g,
p , \rho^\perp] : = [ g'g, p , \rho^\perp]$. It is routine to verify
that this action preserves the weak symplectic form $\Omega $ and
that the map  
\begin{equation}
\label{induced momentum}
\mathbf{J}^G_M([g, p, \rho^\perp]) =
\operatorname{Ad}^\ast_{g^{-1}} \left(\mathbf{J}_P^H(p) +
\rho^\perp \right)
\end{equation}
satisfies the hypotheses of Proposition \ref{prop: mom}. We conclude
hence the following result.

\begin{proposition}
\label{induction theorem second}
The map $\mathbf{J}_M^G: G \times _H
(P \times \mathfrak{h}_\ast ^\perp) \rightarrow \mathfrak{g}_\ast$
given by \eqref{induced momentum} is a $G $-equivariant momentum map.
\end{proposition}

The goal of the induction construction has now been achieved:
starting with the Hamiltonian $H $-space $(P, \omega)$, where $H $
is a closed Lie subgroup of a Lie group $G $, a new Hamiltonian $G
$-space has been constructed, namely $(G \times _H
(P \times \mathfrak{h}_\ast ^\perp), \Omega)$.

\section{Induction and coinduction from $L^1(\mathcal{H})$}
\label{section: induced and coinduced from ell one}

\paragraph{The Banach Lie-Poisson space
$L^1(\mathcal{H}) $.} The Banach space of trace class
operators
$(L^1(\mathcal{H}), \|\cdot \|_1)$ on a separable Hilbert
space $\mathcal{H}$ has a canonical Banach Lie-Poisson
bracket  defined by
\begin{equation}
\label{LP on ell one}
\{f, g\}(\rho) = \operatorname{Tr}(\rho[Df(\rho),
Dg(\rho)]),
\end{equation}
where $f, g  \in C ^{\infty}(L^1(\mathcal{H}))$ and the
Fr\'echet derivatives $Df(\rho), Dg(\rho) $ are regarded
as elements of the Banach Lie algebra
$(L^{\infty}(\mathcal{H}), \|\cdot \|_\infty)$ of bounded
operators on $\mathcal{H}$, identified with the dual of
$L^1(\mathcal{H}) $ by the strongly nondegenerate pairing
\begin{equation}
\label{pairing}
\langle \rho, x \rangle = \operatorname{Tr}(\rho x),
\quad \text{for} \quad  
\quad \rho \in L^1(\mathcal{H}), \; x \in L ^{\infty}
(\mathcal{H} ).
\end{equation}
Hamilton's equations defined by the Poisson bracket
\eqref{LP on ell one} are easily verified to be given in
Lax form (see \cite{OR} for details)
\begin{equation}
\label{l one Hamiltonian equation}
\frac{d \rho}{dt} = [ Dh(\rho),  \rho ].
\end{equation}

 The orthonormal basis
$\{|n\rangle\}_{n = 0} ^{\infty}$ of $\mathcal{H}$, that
is, $\langle n | m \rangle = \delta_{nm} $ for $n, m \in
\mathbb{N} \cup \{0\}$, induces the Schauder basis
$\{|n\rangle \langle m|\}_{n, m = 0} ^{\infty}$ of
$L^1(\mathcal{H})$ since it is orthonormal in the Hilbert
space $L ^2(\mathcal{H})$ of Hilbert-Schmidt operators 
and $L ^1(\mathcal{H}) \subset L ^2(\mathcal{H}) $. Thus,
every trace class operator $\rho\in L ^1(\mathcal{H}) $ can
be uniquely expressed as 
\begin{equation}
\label{rho coordinate}
\rho = \sum_{n,m = 0} ^{\infty}\rho_{nm} 
|n\rangle \langle m| ,
\end{equation}
where the series is convergent in the $\|\cdot
\|_1$ topology. The coordinates $\rho_{nm} \in
\mathbb{R}$ are given by
$\rho_{nm} = \operatorname{Tr}(\rho |m\rangle \langle
n|)$.  The rank one projectors $|l \rangle \langle k|$
thought of as elements of $L ^{\infty}(\mathcal{H})$, by
giving their values on the Schauder basis of $L
^1(\mathcal{H}) $ as
$\operatorname{Tr}(|l\rangle \langle k|\,|n\rangle
\langle m|) = \delta_{kn} \delta_{l m}$, form a
biorthogonal family of functionals (see \cite{LiTz1977})
in $L ^{\infty}(\mathcal{H}) $ associated to the given
Schauder basis $\{|n\rangle \langle m|\}_{n,m = 0}
^{\infty}$ of $L ^1(\mathcal{H}) $. Therefore, each bounded
operator $x \in L ^{\infty}(\mathcal{H}) $ can be uniquely
expressed as 
\begin{equation}
\label{x coordinate}
x = \sum_{l,k = 0} ^{\infty}x_{lk} |l\rangle \langle k|,
\end{equation}
where the series is convergent in the
$w^\ast$-topology. The coordinates
$x_{lk} \in \mathbb{R}$ are also given by $x_{lk} =
\operatorname{Tr}(x |k\rangle \langle l|)$.
Recall that $w ^\ast$-convergence of the series
\eqref{x coordinate} means that the numerical series
\[
\sum_{l,k=0}^{\infty} x_{lk} \operatorname{Tr}(\rho
|l\rangle \langle k|) = \sum_{l,k=0}^{\infty} x_{lk}
\rho_{kl} = \operatorname{Tr}(x \rho)
\] 
is convergent for any $\rho \in  L ^1(\mathcal{H})$.

Since the separable Hilbert space
$\mathcal{H}$ is fixed throughout this paper we shall
simplify the notation by writing $L^1: = L^1(\mathcal{H})$
and $L^{\infty}: = L ^{\infty}(\mathcal{H})$.
\medskip 

\noindent \textbf{Shift operator notation.} The {\bfi shift
operator\/}  
\begin{equation}
\label{shift operator}
S: = \sum_{n=0} ^{\infty} |n \rangle \langle n+1|,
\end{equation}
and its adjoint 
\begin{equation}
\label{transpose shift operator}
S^T: = \sum_{n=0} ^{\infty} |n + 1 \rangle \langle n|,
\end{equation} 
turn out to give a very convenient coordinate description
of various objects that we shall study in this paper. Note 
that the matrix of $S$  has all entries of the upper
diagonal equal to one and all other entries equal to
zero whereas the matrix of $S^T$ has all entries of the
lower diagonal equal to one and all other entries equal to
zero. To facilitate various subsequent computations, we
note that
\begin{equation}
\label{useful shift identities}
S^k (S^T)^k = \mathbb{I}, \qquad (S^T)^k S^k = \mathbb{I} -
\sum_{i=0}^{k-1}p_i, \quad \text{for} \quad k = 1,2,
\dots, 
\end{equation}
where $p_i  = |i\rangle \langle i|: \mathcal{H}\rightarrow
\mathcal{H}$ are the orthogonal projectors on
$\mathbb{R}|i\rangle \subset \mathcal{H}$ for any
$i\in \mathbb{N} \cup \{0\}$. Let $L_0 ^{\infty} \subset
L^{\infty}$ and
$L _0 ^1 \subset L ^1$ denote the closed subspaces of
diagonal operators and define the bounded linear
operators $s , \tilde{s}$ on both $L_0 ^{\infty}$ and
$L_0 ^1 $ by 
\begin{equation}
\label{def of s tilde s}
\left.
\begin{array}{lclc}
&Sx = s(x) S \quad &\text{or} \quad &x S^T = S^T s(x) \\
&S^T x = \tilde{s}(x) S^T
\quad &\text{or} \quad  &x S = S \tilde{s}(x) 
\end{array}
\right\}
\end{equation}
for $x \in L_0 ^{\infty}$ or $x \in L_0 ^1$. The effect
of the  map $s$ is that the $i$th coordinate of $s(x)$
equals the $(i+1)$st coordinate of $x$, that is, $s(x_0,
x_1,x_2 \dots, x_n, \dots) : = (x_1, x_2, \dots, x_n,
\dots)$ for any $(x_0,x_1,x_2 \dots, x_n, \dots) \in \ell
^{\infty} \cong L ^{\infty} _0$. Similarly, the
effect of the map
$\tilde{s}$ is that the
$i$th coordinate of $\tilde{s}(x)$ equals the $(i-1)$st
coordinate of $x$ and the zero coordinate of
$\tilde{s}(x)$ is zero, that is, $\tilde{s}(x_0,
x_1,x_2 \dots, x_n, \dots) : = (0, x_0, x_1, x_2, \dots, x_n,
\dots)$. Thus
\begin{equation}
\label{useful s identities}
s^k \circ \tilde{s}^k = \operatorname{id} \qquad \text{and}
\qquad \tilde{s}^k \circ s^k =  
M_{\mathbb{I} - \sum_{i=0}^{k-1}p_i}, \qquad k = 1,2, \dots,
\end{equation}
where $M_y: L ^{\infty}_0 \rightarrow L ^{\infty}_0 $ is
defined by $M_y(x) := yx $ for any $y \in L ^{\infty}_0$.
The following identities are useful in several
computations later on:
\begin{equation}
\label{trace s tilde s}
\operatorname{Tr}(\rho s(x)) =
\operatorname{Tr}(\tilde{s}( \rho) x) \quad \text{and}
\quad \operatorname{Tr}(s(\rho) x) =
\operatorname{Tr}(\rho\tilde{s}(x))
\end{equation}
for any $\rho\in L ^1_0 $ and $x \in L ^{\infty}_0 $, which
means that $s $ and $\tilde{s}$ are mutually adjoint
operators.

Any $x \in L^{\infty}$ and $\rho \in L ^1$ can be written as 
\begin{equation}
\label{x in coordinates}
x = \sum_{j=1}^{\infty}
(S^T)^j x_{-j} + x_0 + \sum_{i=1}^{\infty}x_{i} S^i ,
\end{equation}
\begin{equation}
\label{rho in coordinates}
\rho = \sum_{j=1}^{\infty}
(S^T)^j \rho_j + \rho_0 + \sum_{i=1}^{\infty}\rho_{-i} S^i ,
\end{equation}
where $x_i, x_0 , x_{-j} \in L^{\infty}_0$ and $\rho_j,
\rho_0, \rho_{-i} \in L ^1_0$. \textit{Note the different
conventions\/}: the indices of the lower diagonals for the
bounded operators are negative whereas for the trace class
operators they are positive. This convention simplifies many
formulas later on.  

The expressions \eqref{x in coordinates} and 
\eqref{rho in coordinates} suggest the introduction, for
every $k \in \mathbb{Z}$, of the Banach subspaces 
\begin{align}
\label{def of l infinity k}
L^{\infty}_k &: = \{\rho \in L^{\infty} \mid \rho_{nm}
= 0\; \text{for} \; m \neq n +k \} \subset L ^{\infty}\\
\label{def of l one  k}
L ^1_k &: = \{\rho \in L ^1 \mid \rho_{nm} = 0\;
\text{for} \; m \neq n +k \} \subset L ^1
\end{align}
consisting of operators whose only  non-zero elements lie on
the $k$th diagonal. We have the following
Schauder decompositions
\begin{equation}
\label{dir sum dec on k}
L^{\infty} =
\bigoplus_{k \in \mathbb{Z}} L^{\infty}_k \qquad \text{and}
\qquad  L^1 =
\bigoplus_{k \in \mathbb{Z}} L^1_k.
\end{equation}
See \cite{singer1981} Ch. III, \S15, namely Definition 15.1 (page 485),
Defintion 15.3 (page 487), and Theorem 15.1 (page 489) for a detailed discussion of this concept and generalizations.
The duality relations between the various spaces
$L^{\infty}_n$ and $L ^1_k$ is given by
\begin{equation}
\label{duality k and  n}
\operatorname{Tr}(\rho_k x_n) =
\delta_{kn} \operatorname{Tr}(\rho_k x _k) \quad \text{if}
\quad \rho_k \in  L ^1_k \quad \text{and} \quad 
x_n   \in L ^{\infty} _{-n}.
\end{equation}
Finally, note that if $k \geq 0$ then $S^k\in L^{\infty}_k$,
$\left(S^T\right)^k\in L ^{\infty}_{-k}$, and  
\begin{equation}
\label{powers of S}
S^l(S^T) ^j = \left\{
\begin{aligned}
&S^{l-j}, \qquad \;\text{if} \quad l \geq j\\
&(S^T)^{j-l}, \quad  \text{if} \quad l \leq j
\end{aligned}
\right.
\end{equation}
which implies 
\begin{equation}
\label{trace formula for y}
\left\langle\rho, x \right\rangle = \sum_{k \in \mathbb{Z}}
\operatorname{Tr} \rho_i x_i
\end{equation}
if $\rho$ and $x $ are expressed in the form 
\eqref{rho in coordinates} and \eqref{x in coordinates}.

\medskip

\noindent \textbf{Banach subspaces of $L^1( \mathcal{H})$
and $L^{\infty}(\mathcal{H})$.} Given the Schauder basis
$\{|n\rangle \langle m|\}_{n, m = 0} ^{\infty}$ of
$L ^1$ (or biorthogonal family of $L ^{\infty}$) inducing the
direct sum splitting  \eqref{dir sum dec on k}, define the
{\bfi transposition operator\/}  $T: L^1
\rightarrow L^1$ (or $T :L ^{\infty}  \rightarrow
L^{\infty}$) by  $(\rho^T)_{ij}: = \rho_{ji}$ for any $i, j
\in \mathbb{N} \cup \{0\}$. We construct the following Banach
subspaces of
$L^1$:
\begin{itemize}
\item $L^1_-: = \oplus_{k = - \infty}^{0}  L ^1_k$ and
$L^1_+: = \oplus_{k = 0}^{\infty}  L ^1_k$ 
\item $L^1_S: = \{ \rho \in L^1 \mid
\rho = \rho^T\}$ and 
$L^1_A: = \{ \rho \in L^1 \mid
\rho = -\rho^T\}$ 
\item $L^1_{-,k}: = \oplus_{i=-k+1}^0 L^1_i$
and $L^1_{+,k}: = \oplus_{i=0}^{k-1} L^1_i$, for $k \geq 1$ 
\item $I^1_{-,k}: =\oplus_{i=-\infty}^{-k} L^1_i$ and 
$I^1_{+,k}: = \oplus_{i=k}^{\infty} L^1_i$, for $k
\geq 1$
\item $L ^1_{S,k} : = L ^1_S \cap \left(L ^1_{+,k}
+ L ^1_{-,k}\right)$ and $L ^1_{A,k} : = L ^1_A \cap
\left(L ^1_{+,k} + L^1_{-,k}  \right)$, for $k
\geq 1$.
\end{itemize}
Relative to operator multiplication, $I^1_{-,k}$ is an
ideal in $L ^1_{-}$, $I^1_{+,k}$ is an ideal in $L ^1_{+}$,
but neither is an ideal in $L^1$. Therefore, relative to the
commutator bracket, the same is true in the associated Banach
Lie algebras.
\medskip

Similarly, using the biorthogonal family of functionals
$\{|l\rangle \langle k|\}_{l,k = 0}^{\infty}$ in $L
^{\infty}$ inducing the direct sum splitting \eqref{dir sum
dec on k}, we construct the following Banach
subspaces of $L^{\infty}$:
\begin{itemize}
\item $L^\infty _-: = \oplus_{k = - \infty}^{0} 
L^{\infty}_k$ and $L^\infty_+: = \oplus_{k = 0}^{\infty} 
L^{\infty}_k$ 
\item $L^\infty_S: = \{ x \in L^\infty \mid
x^T = x\}$ and  $L^\infty_A: = \{ x \in L^\infty \mid
x^T = - x\}$ 
\item $L^\infty_{-,k}: = \oplus_{i=-k+1}^{0}
L^{\infty}_i$ and  $L^\infty_{+,k}: = \oplus_{i=0}^{k-1}
L^{\infty}_i$, for $k \geq 1$ 
\item $I^\infty_{-,k}: = \oplus_{i=-\infty}^{-k}
L^{\infty}_i$ and $I^\infty_{+,k}: = \oplus_{i=k}^{\infty}
L^{\infty}_i$, for $k \geq 1$
\item $L ^\infty_{S,k} : = L ^\infty_S \cap
\left(L^\infty_{+,k} + L ^\infty_{-,k}\right)$ and
$L^\infty_{A,k} : = L ^\infty_A \cap
\left(L ^\infty_{+,k} + L^\infty_{-,k}  \right)$, for $k
\geq 1$.
\end{itemize}

The following splittings of Banach spaces of trace class
operators
\begin{equation}
\label{trace class splittings}
L^1 = L^1_- \oplus I^1_{+,1}, \qquad 
L^1 = L^1_S \oplus I^1_{+,1}, \qquad 
L^1_- = L ^1_{-,k} \oplus I^1_{-,k}
\end{equation}
and of bounded operators
\begin{equation}
\label{bounded splittings}
L^{\infty} =  L^{\infty}_+ \oplus I^{\infty}_{-,1}, \qquad 
L^{\infty} = L ^{\infty}_+ \oplus L^\infty_A, \qquad 
L^{\infty}_+ = L^{\infty}_{+,k} \oplus I^{\infty}_{+,k}
\end{equation}
will be used below. The  strongly nondegenerate pairing
\eqref{pairing} relates the splittings \eqref{trace class
splittings} and
\eqref{bounded splittings} by 
\begin{equation}
\label{duality between the splittings}
\left.
\begin{array}{ccc}
(L ^1_-)^\ast \cong (I ^1_{+, 1})^\circ =
L^{\infty}_+,& \quad 
(L^1_S) ^\ast \cong (I ^1_{+,1}) ^ \circ 
= L ^{\infty}_+ &
\quad 
(L^1_{-,k})^\ast \cong (I ^1_{-, k}) ^\circ =
L^{\infty}_{+,k}\\ 
(I ^1_{+, 1}) ^\ast \cong (L^1_-)^\circ =
I^{\infty}_{-,1} &
\quad 
(I^1_{+,1})^\ast \cong (L^1_S)^\circ = L^\infty_A,& \quad 
(I ^1_{-, k}) ^\ast \cong (L^1_{-,k})^\circ =
I^{\infty}_{+,k} 
\end{array}
\right.
\end{equation}
where, as usual, $^\circ $ denotes the annihilator of the
Banach subspace in the dual of the ambient space.

The splittings \eqref{trace class splittings} and
\eqref{bounded splittings} define six projectors of
$L^1$ and $L ^{\infty}$, respectively. Let
$P^1_-, P^1_0, P^1_+: L^1 \rightarrow L^1$ be the
projectors whose ranges are $I^1_{-,1}$, $L^1_0$, and
$I^1_{+,1}$ defined by the splitting
$L ^1 = I^1_{-,1} \oplus L^1_0 \oplus I^1_{+,1}$. In
particular
$P^1_- + P^1_0 + P^1_+ =
\mathbb{I}$.  Let
$P_{-,k}^1: L^1_- \rightarrow L ^1_{-} $ be the
projector whose range is $L ^1_{-,k}$ defined by the
splitting $L^1_- = L^1_{-,k} \oplus I^1_{-,k} $.  Define
the six projectors 
\begin{equation}
\label{six trace class projectors}
\left.
\begin{array}{ccc}
R_{-}:= P^1_- + P ^1_0 ,&
\qquad 
R_S: = P ^1 _- + P^1_0 + T \circ P ^1_-,&
\qquad 
R_{-,k}:= P ^1_{-,k}\\  
R_+: = P ^1_+,&
\qquad 
R_{S, +}: = P ^1_+ - T \circ P^1_-, & \qquad 
R_{ik}: = R_{-}|_{L ^1_-} - R_{-,k}
\end{array}
\right.
\end{equation}
associated to the
splittings \eqref{trace class splittings}. The order of
presentation of these projectors corresponds to the order
of the splittings in
\eqref{trace class splittings}.

Similarly, the six projectors associated to the dual
splittings
\eqref{bounded splittings} are given by
\begin{equation}
\label{six bounded projectors}
\begin{array}{ccc}
R_{-}^\ast:= P^{\infty}_+ + P^{\infty}_0 ,&
\quad 
R_S ^\ast: = P^{\infty} _+ + P ^{\infty}_0 + T \circ
P^{\infty}_-,&
\quad 
R_{-,k} ^\ast := P^{\infty}_{+,k}\\  
R_+ ^\ast: = P^{\infty}_-,&
\quad 
R_{S, +} ^\ast: = P^{\infty}_- - T \circ P^{\infty}_-,&
\quad  R_{ik}^\ast: = R _- ^\ast|_{L^{\infty}_+}  -
P^{\infty}_{+,k}
\end{array}
\end{equation}
where $P ^\infty _-, P^\infty_0, P^\infty_+:
L^\infty \rightarrow L^\infty$ are the projectors
whose ranges are $I ^{\infty}_{-,1}, L^{\infty}_0,
I^{\infty}_{+,1}$ defined by the splitting
$L^\infty = I ^\infty_{-,1} \oplus L ^\infty_0
\oplus I ^\infty_{+,1} $ and $P_{+,k} ^{\infty} :
L^{\infty}_+ \rightarrow L^{\infty}_{+}$ is the projector
with range
$L^{\infty}_{+,k}$ defined by the splitting $L^{\infty}_+
= L^{\infty}_{+,k} \oplus I^{\infty}_{+,k}$.
 
All Banach spaces appearing in \eqref{bounded
splittings}, with the exception of
$L^{\infty}_{k}$ and $L^{\infty}_{+,k}$, are Banach
subalgebras of $L^{\infty}$ or $L^{\infty}_+$ whereas
$I^{\infty}_{+,k}$, for $k \in \mathbb{N}$, are ideals of
the Banach algebra $L^{\infty}_+$ (but not of
$L^{\infty}$). Therefore, $I^{\infty}_{+,k}$ define a 
filtration of $L^{\infty}_+$ and hence $L ^{\infty}_{+,k}
\cong L ^{\infty}_+/I ^{\infty}_{+,k}$ inherits the
structure of an associative Banach algebra. Thus all these
associative Banach algebras are naturally Banach Lie
algebras. The same considerations apply to the
Banach ideals $I^{\infty}_{-,k} \subset L^{\infty}_-$.

It will be useful in our subsequent development to
distinguish between the projectors defined in 
\eqref{six trace class projectors} and
\eqref{six bounded projectors} and the corresponding maps
onto their ranges. We shall denote by $\pi_-, \pi_+,
\pi_S$, and $\pi_{S,+}$ the maps on $L ^1$ equal to $R_-,
R_+, R_S$, and $R_{S,+}$ but viewed as
taking values in $\operatorname{im}R_- = L^1_-$,
$\operatorname{im}R_+ = I ^1_{+,1}$, $\operatorname{im}R_S =
L_S^1$, and $\operatorname{im}R_{S,+} = I ^1_{+,1}$,
respectively. Similarly, denote by $\pi_{-,k}$ and
$\pi_{ik}$ the maps on $L ^1_-$ equal to $R_{-,k}$ and
$R_{ik}$, but viewed as having values in
$\operatorname{im}R_{-,k} = L ^1_{-,k}$ and
$\operatorname{im}R_{ik} = I ^1_{-,k}$, respectively. For
the projectors on $L^{\infty}$ we shall denote by
$\pi^{\infty}_+, \pi^{\infty}_-, \pi^{\infty}_S$, and
$\pi^{\infty}_A$ the maps equal to
$R^\ast_-, R^\ast_+, R ^\ast_S$, and $R ^\ast_{S,+}$ viewed
as having values in 
$\operatorname{im}R^\ast_- = L^{\infty}_+$,
$\operatorname{im}R^\ast_+ = I^{\infty}_{-,1}$,
$\operatorname{im}R^\ast_S = L^{\infty}_S$, and
$\operatorname{im}R ^\ast_{S,+} = L ^{\infty}_A $, 
respectively. Finally, let $\pi^{\infty}_{+,k}$ and
$\pi^{\infty}_{ik}$ denote the maps on $L^{\infty}_+$
equal to $R ^\ast_{-,k}$ and $R ^\ast_{ik}$ viewed as
having values in $\operatorname{im}R
^\ast_{-,k} = L^{\infty}_{+,k}$ and
$\operatorname{im}R ^\ast_{ik} = I^{\infty}_{+,k}$,
respectively.

\paragraph{Associated Banach Lie groups.}
Note that the Banach Lie group 
\begin{equation}
\label{bounded invertible group}
GL ^{\infty}: = \{x
\in L ^{\infty} \mid x \; \text{is invertible}\}
\end{equation}
has Banach Lie algebra $L^{\infty}$ and is open in $L
^{\infty}$. Define the closed Banach Lie subgroup of upper
triangular operators in $GL^{\infty}$ by
\begin{equation}
\label{gl infinity plus}
GL ^{\infty}_+: = GL ^{\infty} \cap L ^{\infty}_+ .
\end{equation}
Since $GL ^{\infty}_{+} $ is open in $L^{\infty}_{+}$, we
can conclude that its  Banach Lie algebra is
$L^{\infty}_{+}$.  Define the closed Banach Lie subgroup
of orthogonal operators in $GL^{\infty}$ by
\begin{equation}
\label{orthogonal group}
O^{\infty}: = \{x \in L^{\infty} \mid x x^T = x^T x =
\mathbb{I}\}.
\end{equation}
The Banach Lie algebra $L^\infty_A$ of $O ^{\infty}$
consists of all bounded skew-symmetric operators.

Denote by  
\begin{align}
\label{normal subgroup k}
GI_{+,k}^{\infty} :=& (\mathbb{I} +
I_{+,k}^{\infty} ) \cap GL_+^{\infty} \nonumber \\
=&  \{ \mathbb{I} + \varphi
\mid \varphi \in I_{+,k}^{\infty},\;
\mathbb{I} + \varphi \text{ is
invertible in } GL_+^{\infty} \}
\end{align} 
the open subset of
$\mathbb{I}+ I_{+,k}^{\infty}$
formed by the group of all bounded invertible upper 
triangular operators whose strictly upper $(k-1)$-diagonals
are identically zero and whose diagonal is the identity.
This is a closed normal Banach Lie subgroup of
$GL_+^{\infty}$ whose Lie algebra is the closed ideal
$I_{+,k}^{\infty}$. 
\smallskip

\noindent \textbf{Remark.} Unlike the situation encountered
in finite dimensions, the set $\mathbb {I}+
I_{+,k}^{\infty}$ does not consist
only of invertible bounded linear isomorphisms. An example
of  an operator in $\mathbb{I} + I_{+,2}^{\infty}$ that is
not onto is given by  $\mathbb {I} - S^2$, where $S$ is
the shift operator defined in \eqref{shift operator},
since $\sum_{n=0}^{\infty}\frac{1}{n+1}|n\rangle
\notin \operatorname{im}(\mathbb{I} - S^2)$. 
\medskip

Returning to the general case, define the product 
\begin{equation}
\label{k product}
x \circ _k y : = \sum_{l = 0}^{k - 1}\left(\sum_{i = 0}^l
x_i s^i(y_{l-i}) \right) S ^l
\end{equation}
of the elements $x = \sum_{i=0}^{k-1}x_i S^i$ and $y =
\sum_{i=0}^{k-1}y_i S^i \in L ^{\infty}_{+,k}$,  where
$x_i, y_i $ are diagonal operators.
Relative to $\circ _k $, the Banach space
$L^{\infty}_{+,k}$ is an associative Banach algebra
 with unity. It is easy to see that the projection map 
$\pi^{\infty}_{+,k}: L^{\infty}_+
\rightarrow (L^{\infty}_{+,k}, \circ_k)$ is an
associative Banach algebra homomorphism whose kernel is
$I_{+,k} ^{\infty}$. So, it defines a Banach algebra
isomorphism
$[\pi^{\infty}_{+,k}]: L^{\infty}_+/I^{\infty}_{+,k}
\rightarrow (L^{\infty}_{+,k}, \circ_k)$ of the factor
Banach algebra $L^{\infty}_+/I^{\infty}_{+,k} $ with
$(L^{\infty}_{+,k}, \circ_k)$.

The associative  algebra $L ^{\infty}_{+,k}$ with
the commutator bracket 
\begin{equation}
\label{k bracket general}
[x,y]_k := x \circ _k y - y \circ _k x
= \sum_{l=0}^{k-1} \sum_{i=0}^l \left(x_i
s^i(y_{l-i}) - y_i s^i(x_{l-i}) \right)S^l  
\end{equation} is
the Banach Lie algebra of the group
\begin{equation}
\label{gl plus k}
GL^{\infty}_{+,k} = \left\{g = \sum_{i=0}^{k-1}
g_i S^i
\; \Big| \; g_i \in L ^{\infty}_0, |g_0| \geq
\varepsilon(g_0)\mathbb{I} \quad \text{for
some}\quad \varepsilon(g_0) >0 \right\}
\end{equation}
of invertible elements in $(L^{\infty}_{+,k}, \circ_k)$. 

\medskip

\noindent\textbf{Remark.} It is important to note that
invertibility in the Banach algebra $(L^{\infty}_{+,k},
\circ_k)$ does not mean invertibility of the operator on
$\mathcal{H}$. For example, $\mathbb{I} - S^2
\in GL^{\infty}_{+,3}$, that is, $\mathbb{I} - S^2$ is an
invertible element in $(L^{\infty}_{+,3}, \circ_3)$, but 
$\mathbb{I} - S^2$ is not an invertible operator, as noted in
the previous remark.

\medskip

Note that $(L ^{\infty}_{+,k}, [\cdot,
\cdot]_k)$ is not a Banach Lie subalgebra of
$L^{\infty}_{+}$. Since $\pi^{\infty}_{+,k}: L^{\infty}_{+}
\rightarrow L^{\infty}_{+,k}$ is also a Banach Lie algebra
homomorphism one has 
\begin{equation}
\label{pm k Lie bracket}
[x,y]_{k}  =
\pi^{\infty}_{+,k}([x,y])
\quad
\text{for} \quad x, y \in L ^{\infty}_{+,k}.
\end{equation}

Note that $\pi^{\infty}_{+,k}(GL^{\infty}_+) \subset
GL^{\infty}_{+,k}$, since every invertible operator in
$L^{\infty}_+$ is mapped by the
homomorphism $\pi^{\infty}_{+,k}$ to an invertible element
of $L^{\infty}_{+,k}$. Moreover, if $x \in
\pi^{\infty}_{+,k}\left(GL ^{\infty}_+\right) \subset
GL^{\infty}_{+,k}$, then
\[
\left(\pi^{\infty}_{+,k}|_{GL ^{\infty}_+}\right)^{-1}(x) =
\left\{g(\mathbb{I} + \psi) \mid \mathbb{I} + \psi \in GI
^{\infty}_{+,k}
\right\} \quad \text{for some} \quad g \in
\left(\pi^{\infty}_{+,k}|_{GL^{\infty}_+}\right)^{-1}(x).
\]
Indeed, if $g' \in
\left(\pi^{\infty}_{+,k}|_{GL^{\infty}_+}\right)^{-1}(x)$,
then there exists some $g\psi \in I^{\infty}_{+,k}$, since
$g$ is invertible, such that $g^{-1} g' = \mathbb{I} +
\psi \in GI^{\infty}_{+,k}$. The next proposition
shows that the restriction of $\pi^{\infty}_{+,k}$ to
$GL^{\infty}_+ $ has range equal to $GL^{\infty}_{+,k}$.

\begin{proposition}
\label{lie group k isom}
The Banach Lie group homomorphism $\pi^{\infty}_{+,k}|_{GL^{\infty}_+}:
GL^{\infty}_+ \rightarrow GL^{\infty}_{+,k}$ is
surjective and induces a Banach Lie group
isomorphism  $\widetilde{\pi_{+,k} ^{\infty}}:
GL^{\infty}_+/GI^{\infty}_{+,k}
\rightarrow GL^{\infty}_{+,k}$ for any $k = 1,2, \dots $.
\end{proposition}

\noindent \textbf{Proof.} To show that $\pi^{\infty}_{+,k}:
GL^{\infty}_+ \rightarrow GL^{\infty}_{+,k}$ is
surjective is equivalent to proving that for any 
$g_0 + g_1 S + \dots + g_{k-1}
S^{k-1} \in GL^{\infty}_{+,k}$ there exists $\varphi_k \in
I^{\infty}_{+,k}$ such that 
\begin{equation}
\label{first gamma equation}
g_0 + g_1 S + \dots + g_{k-1} S^{k-1} +
\varphi_k \in GL ^{\infty}_+.
\end{equation}

Assume for the moment that \eqref{first gamma equation}
holds. We shall draw a consequence from it. By
\eqref{gl plus k}, $g_0 + g_1 S + \dots +
g_{k-1} S^{k-1}$ is in $GL ^{\infty}_{+,k} $
if and only if $g_0$ is invertible.
Decompose $\varphi_k = \alpha_k S ^k g_0 +
\alpha_{k+1}$, where
$\alpha_{k+1} \in I ^{\infty}_{+,k+1}$. Choosing $N \in
\mathbb{N}$ large enough so that $\mathbb{I} - \frac{1}{N}
\alpha_k S^k\in GL^{\infty}_+$, we obtain
\begin{align}
\label{alpha equation}
GL ^{\infty}_+ \ni &\left(\mathbb{I} - \frac{1}{N}
\alpha_k S^k \right)^N \left(g_0 + g_1 S
+ \dots  + g_{k-1} S^{k-1} + \alpha_k S^k g_0 +
\alpha_{k+1}\right)
\nonumber \\
& = g_0 + g_1 S + \dots  + g_{k-1} S^{k-1}
+ \varphi_{k+1},
\end{align}
where
\begin{align}
\label{phi k plus one}
\varphi_{k+1} =& \left(\sum_{j = 2}^N \left(
\begin{array}{c}
N\\ j
\end{array} 
\right) (-1)^j \frac{1}{N^j}\left(\alpha_k S^k \right)^j
\right) \left(g_0 + g_1 S
+ \dots  + g_{k-1} S^{k-1} + \alpha_k S^k g_0 +
\alpha_{k+1} \right) \nonumber \\
& \quad  + \alpha_{k+1}  - \alpha_k S^k \left(g_1 S
+ \dots  + g_{k-1} S^{k-1} + \alpha_k S^k g_0 +
\alpha_{k+1} \right) \in I ^{\infty}_{+,k+1}.
\end{align}
Therefore, if $g_0 + g_1 S + \dots  +
g_{k-1} S^{k-1} + \varphi_k \in GL ^{\infty}_+ $ for
some $\varphi_k \in I ^{\infty}_{+,k} $, then there exists
some $\varphi_{k+1} \in I ^{\infty}_{+,k+1} $ such that
$g_0 + g_1 S + \dots  + g_{k-1} S^{k-1} +
\varphi_{k+1} \in GL^{\infty}_+ $.

Now we prove the proposition by induction on $k$.

If $k =1$, then $g_0 \in GL^{\infty}_+ $ by
definition. Next, let us assume that \eqref{first gamma
equation} holds. As we just saw, it follows that
\eqref{alpha equation} holds. Consider then $g_0 +
g_1 S + \dots + g_{k-1} S ^{k-1} + g_k S ^k \in
GL^{\infty}_{+,k}$ and decompose it in the group
$GL^{\infty}_{+,k}$  as $g_0 + g_1 S +
\dots + g_{k-1} S ^{k-1} + g_k S ^k = (\mathbb{I}
+ g_k S ^k g_0 ^{-1})\circ _k (g_0 +
g_1 S + \dots + g_{k-1} S ^{k-1})$. Let us
assume, that
$\|g_k\| <\operatorname{min}(1, \|g_0 \|)$ which
implies that $\|g_k S g_0 ^{-1}\| <1$ and hence
that $\mathbb{I} + g_k S^k g_0^{-1} \in
GL^{\infty}_+$. By \eqref{alpha equation} there exists
$\varphi_{k+1} \in I ^{\infty}_{+, k+1}$ such that
$g_0 + g_1 S + \dots + g_{k-1} S ^{k-1} +
\varphi_{k+1} \in GL^{\infty}_+$. Thus we get
\[
(\mathbb{I}
+ g_k S ^k g_0 ^{-1}) (g_0 +
g_1 S + \dots + g_{k-1} S ^{k-1} + \varphi_{k+1})
= 
g_0 + g_1 S + \dots + g_k S ^k + \psi_{k+1}
\in GL ^{\infty}_+
\]
for 
\[
\psi_{k+1} = (\mathbb{I} + g_k S^kg_0^{-1})
\varphi_{k+1} + g_k S^kg_0^{-1}(g_1 S +
\dots + g_{k-1} S^{k-1}) \in I ^{\infty}_{+,k+1}
\]
which proves the assertion \eqref{first gamma equation}
for any element in the connected component of
$GL^{\infty}_{+,k}$. Since
$\{\mathbb{I} + g_1 S + \dots + g_k S^k \mid g_1,
\dots , g_k \; \text{diagonal operators in} \;
L^{\infty}\}$ is a connected Banach Lie subgroup of the
connected component of $GL^{\infty}_{+,k}$ and any element
of $GL ^{\infty}_{+,k} $ can be written as a product of an
element of this group and the Banach Lie subgroup
$GL^{\infty}_{+,1}$ of diagonal operators, it follows
that \eqref{first gamma equation} holds for any element
in $GL^{\infty}_{+,k}$. 
\quad $\blacksquare$

\medskip

In the Banach Lie group $(GL^{\infty}_{+,k}, \circ_k)$, the
inverse $g^{-1} = g_0 ^{-1} + h_1 S + \dots + h_{k-1}
S^{k-1}$ of $g = g_0 + g_1 S + \dots + g_{k-1}
S^{k-1} \in GL^{\infty}_{+,k}$ is given by
\begin{equation}
\label{inverse k explicit}
h_p = - g_0 ^{-1} \left[\sum_{r=1}^{p-1} \sum (-1)^{r-1}
g_{i_1} s^{j_1}(g_0 ^{-1}g_{i_2}) \dots s^{j_q}(g_0
^{-1}g_{i_q}) \dots s^{j_r}(g_0 ^{-1}g_{i_r})
\right] s^p(g_0 ^{-1}),
\end{equation}
$1 \leq p \leq k-1$, where the second sum is taken over all
indices $\{i_1,
\dots , i_r, j_1, \dots, j_r\}$ such that $i_1 +
\dots + i_r = p$ (equality between the $i_q$ is permitted),
$0\leq i_1,
\dots , i_r \leq p$, $1 \leq i_1 = j_1 < j_2 < \dots <
j_r = p - i_r \leq p-1$. For example, here are the first
elements:
\begin{align*}
h_1 &= -g_0^{-1} g_1 s(g_0 ^{-1})\\
h_2 & = -g_0^{-1} \left[g_2 - g_1 s(g_0 ^{-1}g_1)\right]
s^2(g_0 ^{-1}) \\
h_3 & =  -g_0^{-1} \left[g_3 - g_2 s^2(g_0 ^{-1}g_1) - g_1
s(g_0 ^{-1}g_2) + g_1s(g_0 ^{-1}g_1) s^2(g_0 ^{-1}g_1)
\right] s^3(g_0 ^{-1}).
\end{align*}

\paragraph{Coinduced Banach Lie-Poisson structures.}
After these preliminary remarks and notations let us apply
the results of the previous section to the Banach
Lie-Poisson space $L^1$. We shall drop the upper indices
``ind" and ``coind" on the Poisson brackets because it will
be clear from the context which brackets are induced and
coinduced on various subspaces.

We start with points (i)
of Proposition
\ref{decomposition proposition} and Proposition
\ref{double diagram}.  So let us consider the diagram

\unitlength=5mm
\begin{center}
\begin{picture}(18,5.5)
\put(4.7,5.1){\makebox(0,0){$L ^1$}}
\put(1,0.5){\makebox(0,0){$L^1_S$}}
\put(9,0.5){\makebox(0,0){$I ^1_{+, 1}$}}
\put(5,4.5){\vector(1,-1){3}}
\put(4,4.5){\vector(-1,-1){3}}
\put(1.4,1.5){\vector(1,1){3}}
\put(7.6,1.5){\vector(-1,1){3}}
\put(2.3,3.5){\makebox(0,0){$\pi_S$}}
\put(3.5,2.5){\makebox(0,0){$\iota_S$}}
\put(5.5,2.5){\makebox(0,0){$\iota_{S,+}$}}
\put(7.3,3.5){\makebox(0,0){$\pi_{S,+}$}}
\put(13,5.1){\makebox(0,0){$L ^1$}}
\put(17,0.5){\makebox(0,0){$L ^1_-$}}
\put(12.3,4.5){\vector(-1,-1){3}}
\put(9.7,1.4){\vector(1,1){3}}
\put(13.4,4.5){\vector(1,-1){3}}
\put(16,1.5){\vector(-1,1){3}}
\put(10.4,3.5){\makebox(0,0){$\pi_{+}$}}
\put(12,2.5){\makebox(0,0){$\iota_+$}}
\put(15.3,3.5){\makebox(0,0){$\pi_-$}}
\put(14,2.5){\makebox(0,0){$\iota_-$}}
\end{picture}
\end{center}
where we recall that $\pi_S, \pi_{S,+}, \pi_+ $ and $\pi_-$
are the projections onto the ranges of $R_S, R_{S+}, R_+
$, and
$R_-$ respectively and $\iota_S, \iota_{S,+}, \iota_+ $,
and
$\iota_-$ are inclusions. We see from the above that
the assumptions in part (i) of Proposition
\ref{double diagram} are satisfied because
$(I^1_{+,1})^\circ = L^{\infty}_+$ is a Banach Lie
subalgebra of $(L ^1) ^\ast = L ^{\infty}$.  Thus we can
conclude the following facts.

\begin{itemize}
\item[{\bf (i)}]
By Proposition \ref{double diagram} (i) it follows that
$L^1_S$ and $L ^1_-$ are isomorphic Banach Lie-Poisson
spaces with the Poisson brackets defined by formula
\eqref{coinduced bracket with projector evaluated}. They
are given, respectively, by
\begin{align}
\label{r s l p bracket}
\{f,g\}_{S}(\sigma) &=
\operatorname{Tr}\left(\iota_S(\sigma)\left[
D(f\circ \pi_S)(\iota_S(\sigma)),  D(g\circ
\pi_S)(\iota_S(\sigma)) 
\right] \right)
\end{align}
for $\sigma \in L^1_S$ and $f, g  \in C ^{\infty}(L^1_S)$
and 
\begin{align}
\label{r minus l p bracket}
\{f,g\}_{-}(\rho) &=
\operatorname{Tr}\left(\iota_-(\rho)\left[
D(f \circ \pi_-)(\iota_-(\rho)), D(g \circ
\pi_-)(\iota_-(\rho)) 
\right] \right)
\end{align}
for $\rho\in L ^1_-$ and $f, g \in C ^{\infty}(L ^1_-)$.

The linear continuous maps $\Phi_{-,S}: = \pi_- \circ
\iota_S: L^1_S \rightarrow L_- ^1$ and $\Phi_{S,-}: = 
\pi_S \circ \iota_-: L ^1_- \rightarrow L^1_S$ are mutually
inverse isomorphisms of the Banach Lie-Poisson spaces
$(L^1_S, \{\,,\}_{S}) $ and $(L^1_-, \{\,,\}_{-})$. 
The coadjoint actions of the Banach Lie group
$GL_+^{\infty}$ on $L^1_-$ and $L ^1_S$ are given by
\begin{equation}
\label{concrete minus coadjoint action}
(\operatorname{Ad}^{+})_{g^{-1}}^{\ast} \rho =
\pi_-(g \iota_-(\rho) g^{-1}) \quad  \text{for} \quad
\rho\in L ^1_- 
\end{equation}
\begin{align}
\label{concrete s coadjoint action}
(\operatorname{Ad}^{S})_{g ^{-1}}^{\ast} \sigma &=
\pi_S(g\iota_S(\sigma) g^{-1} ) 
\quad \text{for}
\quad \sigma \in L^1_S
\end{align}
and $g \in GL^{\infty}_+$. Differentiating these formulas
relative to $g $ at the identity, we get
\begin{equation}
\label{concrete minus coadjoint algebra action}
(\operatorname{ad}^{+})_{x}^{\ast} \rho =
- \pi_-([x,\iota_-(\rho)]) \quad  \text{for} \quad
\rho\in L ^1_- 
\end{equation}
\begin{align}
\label{concrete s coadjoint algebra action}
(\operatorname{ad}^{S})_{x}^{\ast} \sigma &=
- \pi_S([x,\iota_S(\sigma)]) 
\quad \text{for}
\quad \sigma \in L^1_S
\end{align}
for $x \in L ^{\infty}_+ $.
The isomorphisms
$\Phi_{-,S}: L^1_S \rightarrow L_- ^1$ and
$\Phi_{S,-} : L ^1_- \rightarrow L^1_S$ are
equivariant relative to these coadjoint actions, that is,
\begin{equation}
\label{concrete coadjoint s minus equivariance}
(\operatorname{Ad}^{S})_{g ^{-1}}^{\ast} \circ \Phi_{S,-}
= \Phi_{S,-} \circ (\operatorname{Ad}^{+})_{g
^{-1}}^{\ast} 
\end{equation}
\begin{equation}
\label{concrete coadjoint minus s equivariance}
(\operatorname{Ad}^{+})_{g ^{-1}}^{\ast}
\circ \Phi_{-,S}
= \Phi_{-,S} \circ
(\operatorname{Ad}^{S})_{g ^{-1}}^{\ast}
\end{equation} 
for any $g \in GL ^{\infty}_+ $.

\item[{\bf (ii)}] By \eqref{duality between the
splittings}, $I_{+,1}^1$ is the predual of the two Banach
Lie algebras 
$I_{-, 1} ^{\infty}$ and $L^\infty_A$. Thus 
\eqref{trace class splittings} -
\eqref{six bounded projectors} and point (ii) of
Proposition \ref{double diagram}  imply that
$I_{+,1}$ carries two different Lie-Poisson
brackets, namely by \eqref{coinduced bracket with
projector evaluated} we have
\begin{align}
\label{first l p on i plus one}
\{f, g\}_{+}(\rho) &= 
\operatorname{Tr}\left(\iota_+(\rho)\left[
D(f \circ \pi_+)(\iota_+(\rho)), D(g \circ
\pi_+)(\iota_+(\rho)) 
\right] \right) 
\end{align} 
and
\begin{align}
\label{second l p on i plus one}
\{f, g\}_{{S,+}}(\rho) &= 
\operatorname{Tr}\left(\iota_{S+}(\rho)\left[
D(f \circ \pi_{S,+})(\iota_{S,+}(\rho)), 
D(g \circ \pi_{S,+})(\iota_{S,+}(\rho))  \right] \right),
\end{align} 
where $\rho \in I ^1_{+,1}$, $f,g \in
C^{\infty}(I^1_{+,1})$. 

The coadjoint actions
$(\operatorname{Ad}^{-})^{\ast}$ and
$(\operatorname{Ad}^{A})^{\ast} $ of the groups
$GI_{-,1} ^{\infty}$ and $O ^{\infty}$ respectively on
$I_{+,1}^1$ are given by
\begin{equation}
\label{concrete coadjoint action of g i minus one}
(\operatorname{Ad}^{-})_{h^{-1}}^{\ast} \rho =
\pi_+(h\iota_+(\rho) h^{-1} )  \quad \text{for}
\quad h \in GI ^{\infty}_{-,1}
\end{equation}
and
\begin{align}
\label{concrete coadjoint action of o} 
&(\operatorname{Ad}^{A})_{g^{-1}}^{\ast} \rho =
\pi_{S+}(g \iota_{S,+}(\rho) g^{-1}) 
\quad \text{for} \quad g \in O^{\infty}
\end{align}
where $\rho \in I_{+,1}^1$. We shall not pursue the
investigation of this interesting case in this paper.
\end{itemize}

\paragraph{Induced Banach Lie-Poisson structures.} We
begin with the study of  the lower triangular case. Denote by
$\iota_{-,k} : L ^1_{-,k}
\hookrightarrow L ^1_-$ the inclusion and let $ \iota _{-,k}
^{-1}: \iota _{-,k}
\left( L ^1_{-,k} \right) \rightarrow  L ^1_{-,k}$ be its
inverse (defined, of course, only on the range of
$\iota_{-,k}$). Then $\iota_{-,k}^\ast: L ^{\infty}_+
\rightarrow L ^{\infty}_{+,k}$. Since
$\ker \iota_{-,k} ^\ast = I_{+,k} ^{\infty}$ is an ideal in
$L^{\infty}_+$, by Proposition
\ref{induction proposition} we have 
$(\operatorname{Ad}^+)^\ast_{g ^{-1}} \iota_{-,k}( L
^1_{-,k})
\subset \iota_{-,k}( L ^1_{-,k})$ for any $g \in
GL^{\infty}_+$. Therefore there are $GL^{\infty}_+$ and
 $L^{\infty}_+$ coadjoint actions on $L^1_{-,k}$ defined by
\begin{align}
\label{group coadjoint minus k action}
(\operatorname{Ad}^{+,k})^\ast_{g ^{-1}} \rho
&:= \iota_{-,k}^{-1} \left(\pi_{-}\left(g (\iota_- \circ
\iota_{-,k})(\rho) g^{-1}\right) \right)
\quad \text{for} \quad \rho \in L^1_{-,k} \quad \text{and}
\quad g \in GL ^{\infty}_{+}\\
\label{algebra coadjoint minus k action}
(\operatorname{ad}^{+,k})^\ast_x \rho &:= 
\iota_{-,k}^{-1}\left(\pi_{-}[x, (\iota_- \circ
\iota_{-,k})(\rho)]\right)
\quad \quad  \;\; \text{for} \quad \rho \in L^1_{-,k} \quad \,
\text{and}
\quad x \in L^{\infty}_{+}.
\end{align}

Since the action \eqref{group coadjoint minus k action} is
trivial for all elements of the closed normal Lie
subgroup $GI^{\infty}_{+,k}$, it induces the coadjoint
action of the group $GL ^{\infty}_{+,k}
\cong GL ^{\infty}_+/GI^{\infty}_{+,k}$ given by
\eqref{group coadjoint minus k action} that will be also
denoted by $(\operatorname{Ad}^{+,k})^\ast$. Similarly, the
Lie algebra action \eqref{algebra coadjoint minus k action}
is trivial for all elements in the closed ideal
$I^{\infty}_{+,k}$ so it induces the coadjoint action of
the Lie algebra $L^{\infty}_{+,k} \cong
L^{\infty}/I^{\infty}_{+,k}$ on $L ^1_{-,k}$ denoted also
by $(\operatorname{ad}^{+,k})^\ast$.

 One can express
\eqref{group coadjoint minus k action} and \eqref{algebra
coadjoint minus k action}
in terms of the expansions $\rho = \rho_0 + S^T \rho_1  +
\dots + (S^T)^{k-1}\rho_{k-1} \in L ^1_{-,k}$, $x = x_0 +
x_1 S + \dots + x_{k-1} S^{k-1} \in L^{\infty}_{+,k}$,
and $g = g_0 + g_1 S + \dots + g_{k-1}S^{k-1} \in
GL^{\infty}_{+,k}$  in the following way
\begin{equation}
\label{group coordinate coadjoint minus k action}
(\operatorname{Ad}^{+,k})^\ast_{g ^{-1}} \rho = 
\sum_{i,j,l = 0,\, j \geq i+l}^{k-1} (S^T)^{j-i-l}
\tilde{s}^l[s^j(\tilde{s}^i(g_i))\rho_j h_l],
\end{equation}
where the diagonal operators $h_l$ are expressed in terms
of the $g_i$ in \eqref{inverse k explicit}, and (using
\eqref{powers of S})
\begin{equation}
\label{algebra coordinate coadjoint minus k action}
\left(\operatorname{ad}^{+,k}\right)^\ast_x \rho=
\sum_{j=0}^{k-1} (S^T)^j
\sum_{i=j}^{k-1}\left(\tilde{s}^{i-j}( \rho_{i} x_{i-j})
 - \rho_i s^j(x_{i-j}) \right).
\end{equation}

By \eqref{k bracket general} and \eqref{trace formula for
y}, the Lie-Poisson bracket on $L^1_{-,k}$ is
given by
\begin{align}
\label{k poisson bracket}
&\{f, g\}_{k}(\rho) = \operatorname{Tr} \left( \rho
\left[Df(\rho), Dg(\rho) \right]_k \right) \nonumber \\
&\quad = \sum_{l=0}^{k-1} \sum_{i=0}^l
\operatorname{Tr} \left[\rho_{l} \left(\frac{
\delta f}{\delta \rho_{i}}(\rho)s^i\left(\frac{\delta
g}{\delta
\rho_{l-i}}(\rho)\right) -
\frac{ \delta g}{\delta
\rho_{i}}(\rho)s^i\left(\frac{\delta f}{\delta
\rho_{l-i}}(\rho)\right) \right) \right]
\end{align} 
for $f, g \in C^{\infty}(L^1_{-,k})$, where 
$\frac{\delta f}{\delta \rho_{i}}(\rho)$ denotes the
partial functional derivative of $f $ relative to
$\rho_i$ defined by 
$Df(\rho) = \frac{\delta f}{\delta
\rho_0} (\rho) + \frac{\delta f}{\delta \rho_1} (\rho)S +
\dots + \frac{\delta f}{\delta \rho_{k-1}} (\rho)S^{k-1}$.

If in the previous formulas we let $k = \infty$
one obtains the Lie-Poisson bracket on $L^1_-$. Indeed,
the Lie-Poisson bracket $\{f, g \}_-$ on $L^1_-$ given by 
\eqref{r minus l p bracket} expressed in the coordinates 
$\{\rho_i\}_{i =0}^{\infty}$ equals  
\eqref{k poisson bracket} for $k = \infty$.

\begin{proposition}
\label{k-induced iso}
The Lie-Poisson bracket \eqref{k poisson bracket} on
$L ^1 _{-,k}$ coincides with the induced bracket 
\eqref{induced bracket third form} determined by the
inclusion $\iota_{-,k}: L ^1_{-,k} \hookrightarrow L^1_-$
and the Lie-Poisson bracket \eqref{r minus l p bracket} on
$L^1_- $.
\end{proposition}
 
\noindent \textbf{Proof.} We need to prove that the induced
bracket \eqref{induced bracket third form} evaluated on two
linear functionals $x, y \in L^{\infty}_{+,k} \cong
(L^1_{-,k})^\ast \subset C ^{\infty}(L ^1_{-,k})$ coincides
with $[x, y]_k$. To see this we note that $D(x \circ
\iota_{-,k} ^{-1} \circ R_{-,k})(\iota_{-,k}( \rho)) =
\iota_{+,k} x \in L ^{\infty}_+$, where $\iota_{+,k} : 
L^{\infty}_{+,k} \hookrightarrow L^{\infty}_+$ is the
inclusion. Then, a direct verification shows that for any
$\rho\in L ^1_{-,k}$ we have
\begin{align*}
\{x, y\}^{ \operatorname{ind}}( \rho) = \left\langle
[\iota_{+,k} x, \iota_{+,k} y], \iota_{-,k}
\rho \right\rangle 
= \operatorname{Tr} \left([x,y] \rho \right)
= \operatorname{Tr} \left([x,y]_k \rho   \right)
\end{align*}
by \eqref{k bracket general}. 
\quad $\blacksquare$

\medskip

Let us study now the symmetric representation of
$\left(L^1_{-,k}, \{ \cdot , \cdot \}_{-,k} \right)$ for $k
\in \mathbb{N} \cup \{ \infty\}$. This will be done by using
the Banach Lie-Poisson space isomorphism $\Phi_{S,-}: = 
\pi_S \circ \iota_-: L ^1_- \rightarrow L^1_S$. Let
$\pi_{-,k}: L^1_- \rightarrow L^1_{-,k}$ and $\pi_{S,k}: L
^1_S \rightarrow L^1_{S,k}$ be the projections with the
indicated ranges and $\iota_{S, k}: L^1_{S,k} \rightarrow
L^1_S$ the inclusion. Define
$\Phi_{S,-,k}: = \pi_{S,k} \circ \Phi_{S,-}
\circ \iota_{-,k}: L ^1_{-,k} \rightarrow L^1_{S,k}$. The
following commutative diagram illustrates these maps:
\unitlength=5mm
\begin{center}
\begin{picture}(9,8.6)
    \put(1,7){\makebox(0,0){$L^1_-$}} 
    \put(9,7){\makebox(0,0){$L^1_S$}}
    \put(1,2){\makebox(0,0){$L ^1_{-,k}$}} 
    \put(9,2){\makebox(0,0){$L ^1_{S,k}$}}
    \put(1.2,6){\vector(0,-1){3}}
    \put(0.7,3){\vector(0,1){3}}
    \put(9.2,6){\vector(0,-1){3}} 
    \put(8.7,3){\vector(0,1){3}}
    \put(2.5,7){\vector(1,0){5}}
    \put(2,2){\vector(1,0){5.7}}
    \put(-0.4,4.3){\makebox(0,0){$\iota_{-,k}$}}
    \put(2.4,4.3){\makebox(0,0){$\pi_{-,k}$}}
    \put(10.6,4.4){\makebox(0,0){$\pi_{S,k}$}}
    \put(7.6,4.4){\makebox(0,0){$\iota_{S,k}$}}
    \put(5,7.8){\makebox(0,0){$\Phi_{S,-}$}}
    \put(5,2.5){\makebox(0,0){$\Phi_{S,-,k}$}}
    \end{picture}
\end{center}
Pushing forward the Poisson bracket $\{ \cdot , \cdot \}_k$ on
$L ^1_{-,k}$ by the Banach space isomorphism $\Phi_{S,-,k}$
endows $L ^1_{S,k}$ with an isomorphic Poisson structure
denoted by $\{ \cdot , \cdot \}_{S,k}$. From Propositions
\ref{double diagram} and
\ref{k-induced iso}, all the maps in the diagram above
are linear Poisson maps, with the exception of $\pi_{-,k}$
and $\pi_{S,k}$ which are not Poisson. Recall that
$GL^{\infty}_+$ acts on
$L ^1_-
$ and $L ^1_S$ by \eqref{concrete minus coadjoint action} and
\eqref{concrete s coadjoint action} respectively, and that 
$GL^{\infty}_{+}$ (and hence $GL^{\infty}_{+,k}$) acts on 
$L^1_{-,k}$ by  \eqref{group coadjoint minus k action}. Using
the isomorphisms $\Phi_{S,-}$ and $\Phi_{S,-,k}$ to push
forward these actions to $L ^1_S $ and $L ^1_{S,k}$,
respectively, all the maps in the diagram above are also
$GL^{\infty}_{+}$-equivariant. Consequently, one has the $GL
^{\infty}_+
$-invariant filtrations
\begin{equation}
\label{blp minus filtration}
\iota_{-,1}(L^1_{-,1}) \hookrightarrow
\iota_{-,2}(L^1_{-,2}) \hookrightarrow
\ldots \hookrightarrow
\iota_{-,k}(L^1_{-,k}) \hookrightarrow
\iota_{-,k+1}(L^1_{-,k+1})
\hookrightarrow
\dots \hookrightarrow L^1_-
\end{equation}
\begin{equation}
\label{blp s filtration}
\iota_{S,1}(L^1_{S,1}) \hookrightarrow
\iota_{S,2}(L^1_{S,2}) \hookrightarrow
\ldots \hookrightarrow
\iota_{S,k}(L^1_{S,k}) \hookrightarrow
\iota_{S,k+1}(L^1_{S,k+1})
\hookrightarrow
\dots \hookrightarrow L^1_S
\end{equation}
of Banach Lie-Poisson spaces predual to the sequence
\begin{equation}
\label{bla filtration}
L^{\infty}_+ \longrightarrow  \dots
\longrightarrow
L^{\infty}_{+,k} \longrightarrow
L^{\infty}_{+,k-1} \longrightarrow
\dots \longrightarrow
L^{\infty}_{+,2} \longrightarrow
L^{\infty}_{+,1} 
\end{equation}
of Banach Lie algebras in which each arrow is the surjective
projector $\pi^{\infty}_{+,k,k-1}: L ^{\infty}_{+,k} 
\rightarrow L ^{\infty}_{+,k-1}$ that maps $k$-diagonal upper
triangular operators to $(k-1)$-diagonal upper triangular 
operators by eliminating the $k$th diagonal. We have
$\pi^{\infty}_{+,k,k-1} \circ \pi^{\infty}_{+, k} =
\pi^{\infty}_{+,k-1} $.

\section{Dynamics generated by Casimirs of
$L^1(\mathcal{H})$}
\label{section: dynamics generated by Casimirs of ell one}

We begin by presenting
Hamilton's equations on $L^1_-$ and $L^1_S$
given by arbitrary smooth functions
$h$ and $f$ defined on the relevant Banach Lie-Poisson
spaces. Using formula \eqref{coinduced vector field} of
Proposition \ref{decomposition proposition}, one
obtains Hamilton's equations
\begin{align}
\frac{d}{dt}{\rho} &= \pi_-\left(\left[
 D(h \circ \pi_-)(\iota_-(\rho)),
\iota_-(\rho)\right]\right) \quad \text{for} \quad \rho
\in L ^1_- \quad \text{and} \quad h \in C ^{\infty}(L ^1_-),
\label{equ on l one minus}\\
\frac{d}{dt}{\sigma} & = \pi_S
\left(\left[ D(f \circ \pi_S)(\iota_S(\sigma)),
\iota_S(\sigma)\right]\right) \quad \, \text{for} \quad 
\sigma \in L^1_S \quad \, \text{and} \quad f \in 
C^{\infty}(L^1_S),
\label{equ on s one}
\end{align}
on the isomorphic Banach Lie-Poisson spaces $(L ^1_-, 
\{ \cdot , \cdot \}_- )$ and $(L ^1_S, \{ \cdot ,
\cdot \}_S)$;  from \S 
\ref{section: induced and coinduced from ell one} we know
that this isomorphism is $\Phi_{S,-}: = 
\pi_S \circ \iota_-: (L ^1_-, 
\{ \cdot , \cdot \}_- ) \stackrel{\sim}\longrightarrow
(L^1_S, \{\cdot ,\cdot \}_S)$. Therefore, if $f \circ 
\Phi_{S,-} = h$ then equations \eqref{equ on l one minus} and
\eqref{equ on s one} give the same dynamics. Recall that
$\pi_-: L ^1 \rightarrow L^1_-$ and $\pi_S: L^1 \rightarrow
L ^1_S$ are, by definition, the projectors $P ^1_- + P _0:
L^1 \rightarrow L^1$ and $\pi_S:= P^1_- + P _0 + T \circ
P^1_-: L^1 \rightarrow L^1$ considered as maps on their
ranges (see \eqref{six trace class projectors} and the
subsequent comments) and $\iota_-: L ^1_- \hookrightarrow
L^1$, $\iota_S: L ^1_S \hookrightarrow L ^1$ are the
inclusions. 

Now let us observe that the family of functions $I _l \in
C^{\infty}( L^1)$ defined by 
\begin{equation}
\label{casimir on l one}
I_l(\rho): = \frac{1}{l} \operatorname{Tr} \rho ^l \qquad
\text{for} \qquad  l \in \mathbb{N}
\end{equation}
are Casimir functions on the Banach Lie-Poisson space $(L^1,
\{\cdot, \cdot \})$. This follows from \eqref{LP on ell one}
since one has 
\begin{equation}
\label{derivative of i l}
DI_l (\rho) = \rho^{l-1} \in L ^1 \subset L ^{\infty} \cong
(L^1) ^\ast.
\end{equation}
Restricting $I _l$ to $\iota_-: L_-^1 \hookrightarrow L ^1$
and $\iota_S: L ^1_S \hookrightarrow L ^1$ we obtain for all
$l \in \mathbb{N}$
\begin{align}
\label{i l minus}
I_l^-(\rho) &:= I_l(\iota_-(\rho)) \qquad
\text{for} \qquad \rho \in L^1_-\\
\label{i l s}
I_l^S (\sigma) &:= I_l(\iota_S (\sigma))
\qquad \text{for} \qquad \sigma \in L^1_S.
\end{align} 
According to Corollary
\ref{involution corollary}(i), \eqref{i l minus} and 
\eqref{i l s}  form two infinite families of integrals in
involution 
\begin{equation}
\label{involution integrals}
\{I^-_l, I^-_m\}_- = 0 \quad \text{and} \quad 
\{ I^S_l, I^S_m\}_S =0 \quad \text{for} \quad l,m \in
\mathbb{N}.
\end{equation}
Since $I _l^S \circ \Phi_{S, -} \neq I _l^-$, the
Hamiltonians $I _l^-$ and $I _l^S$ define on $(L ^1_-, 
\{ \cdot , \cdot \}_- )$ (or $(L ^1_S, \{ \cdot ,
\cdot \}_S)$) different families of dynamical systems.
\medskip

Firstly, we shall investigate the systems associated to 
the Hamiltonians $I _l^-$ given by 
\eqref{i l minus}. As we shall see, the framework of the
Banach Lie-Poisson space $(L^1_-, \{ \cdot , \cdot \}_-)$
is more natural in this case. Hence, taking into account 
Corollary \ref{involution corollary}(ii), substituting
$I_l^{-}$ into \eqref{equ on l one minus}, then applying
$\iota_-$ to \eqref{equ on l one minus}, and
using \eqref{derivative of i l}, yields the family of
Hamilton equations on $L^1_-$ 
\begin{equation}
\label{l minus one i l}
\frac{\partial \iota_-(\rho)}{\partial t_l} 
= \left(P^1_- + P^1_0 \right)\left[\left(P^{\infty}_+ +
P^{\infty}_0 \right)\left([\iota_-(\rho)]^{l-1}\right),
\iota_-(\rho)
\right]
\end{equation}
or, equivalently, in Lax form
\begin{equation}
\label{l minus one i l 1}
\frac{\partial \iota_-(\rho)}{\partial t_l}
= - \left[P^{\infty}_- 
\left([\iota_-(\rho)]^{l-1}\right), \iota_-(\rho) \right]
= \left[P_0 ^{\infty}\left([\iota_-(\rho)]^{l-1}\right),
\iota_-(\rho) \right],
\end{equation}
where $t_l$ denotes the time parameter for the
$l$th flow.  

Equation \eqref{l minus one i l} implies that its solution
is given by the coadjoint action of the group $GL
^{\infty}_+$ on the dual $L ^1_-$ of its Lie algebra. Hence,
there is some smooth curve $\mathbb{R} \ni t_l \mapsto 
h_+(t_l) \in GL^{\infty}_+$ satisfying $\left(\operatorname{Ad}^+
\right)^\ast_{h_+(t _l) ^{-1}} \circ \left(\operatorname{Ad}^+
\right)^\ast_{h_+(s _l) ^{-1}} = \left(\operatorname{Ad}^+
\right)^\ast_{h_+(t _l + s _l) ^{-1}}$ such that 
\begin{equation}
\label{solution l minus one i l minus}
\iota_-(\rho(t_l)) = 
(\operatorname{Ad}^+) ^\ast_{h_+(t_l)^{-1}} \rho(0)
=(P^1_- + P^1_0)\left(h_+(t_l)
\iota_-(\rho(0)) h_+(t_l)^{-1}\right)
\end{equation}
is the solution of \eqref{l minus one i l} with initial
condition $\rho(0)$ for $t_l = 0$. 

On the other hand,  the  solution of                 
\eqref{l minus one i l 1} is given by
\begin{equation}
\label{solution l minus one i l s}
\iota_-(\rho(t_l)) = h_-(t_l)^{-1} \iota_-(\rho(0))
h_-(t_l),
\end{equation}
for a  smooth one-parameter subgroup $\mathbb{R} \ni t_l \mapsto
h_-(t_l)
\in GL^{\infty}_-$ that can be explicitly determined. We
shall do this by using  the decomposition 
$\iota_-( \rho)  = \rho_0 + \rho_- $, where $\rho_- =
\sum_{i=1}^{\infty} (S^T)^i \rho_i$ and $\rho_i \in L ^1_0 $
if $i \in \mathbb{N}\cup \{ 0 \}$. Since
$P_0^{\infty}\left([\iota_-(\rho)]^{l-1}\right) =
\rho_0^{l-1}$, equation
\eqref{l minus one i l 1} becomes
\[
\frac{\partial}{\partial t_l} \iota_- ( \rho) 
 = [ \rho_0 ^{l-1}, \rho_0 + \rho_- ] = [ \rho_0 ^{l-1},
\rho_- ]
\]
which is equivalent to 
\begin{equation}
\label{easy equation}
\frac{\partial}{\partial t_l} \rho_- = [\rho_0 ^{l-1},
\rho_-] \quad \text{and} \quad 
\frac{\partial}{\partial t_l} \rho_0 =0.
\end{equation}
It immediately follows that its solution is given by 
\eqref{solution l minus one i l s} with 
\begin{equation}
\label{g solution}
h_-(t_l) = e ^{-t_l \rho_0(0)^{l-1}},
\end{equation} 
where $\rho(0) = \rho_0(0) + \rho_-(0)$ is the initial
value of $\rho$ at time $t_l = 0 $.

Note that $h_-(t _l) \in GL ^{\infty}_-$ is in fact a
diagonal operator which can also be obtained from the
decomposition
\begin{equation}
\label{first decomposition diagonal}
e^{t_l [\iota_-(\rho(0))]^{l-1}} = k_-(t_l)h_-(t_l)^{-1},
\end{equation}
where $k_-(t_l) \in GI ^{\infty}_{-,1}$. It
follows that we can write the solution also in the form
\begin{equation}
\label{two form solution}
\iota_-( \rho(t_l))  
 = k_-(t_l) ^{-1}\left[\iota_-( \rho(0))
\right]k_-(t_l).
\end{equation}

Finally, note that in \eqref{solution l minus one i l minus}
we can choose $h_+(t _l) = h_-(t _l)$ since also $h_-(t _l) \in
GL^{\infty}_+$.
\medskip

Let us analyze the system \eqref{l minus one i l 1} in more detail.
We begin by noting that there is an isometry
between $\ell^\infty$ and the diagonal bounded linear
operators $L^{\infty}_0 \subset L^{\infty}$ and between
$\ell^1$ and the diagonal trace class operators $L^1_0
\subset L^1$. Fix a strictly lower triangular element 
\begin{equation}
\label{def of lambda}\nu_-
=\sum_{i=1}^{k-1} (S^T)^i \nu_i \in L^1_{-,k} \quad
\text{where} \quad k \in \mathbb{N} \cup \{ \infty\}
\end{equation}
and define the map $\mathcal{J}_{\nu_-}:
\ell ^{\infty} \times \ell ^1 \rightarrow L ^1_{-,k}$ by
\begin{equation}
\label{def of cal j}
\mathcal{J}_{\nu_-}(\mathbf{q}, \mathbf{p}): =
\mathbf{p} + e^{\mathbf{q}} \nu_- e^{-\mathbf{q}},
\end{equation}
where, on the right hand side, we identify $\mathbf{p}$ and
$\mathbf{q}$ with diagonal operators and $e^{\mathbf{q}}$
is the exponential of $\mathbf{q}$. It is easy to see that
this map is smooth and that
$\mathcal{J}_{\nu_-}(\mathbf{q}, \mathbf{p}) =
\mathcal{J}_{\nu_-}(\mathbf{q}+ \alpha\mathbb{I},
\mathbf{p})$, for any $\alpha \in \mathbb{R}$. We shall prove in
Proposition \ref{Flaschka map theorem} that if $\nu_- = (S^T)^{k-1}
\nu_{k-1} \in L ^1_{-k+1} \subset L ^1_{-,k}$, the map
$\mathcal{J}_{n_-}: \ell ^{\infty} \times \ell ^1 \rightarrow
I^1_{-,0,k-1}$, the space of bidiagonal trace class operators having
non-zero entires only on the main and the lower $(k-1)$st diagonal,
is a momentum map in the sense of Definition 
\ref{def momentum map weak}

We shall argue below, in analogy with the finite dimensional case,
that $(\mathbf{q}, \mathbf{p})$ can be considered as 
angle-action coordinates for the Hamiltonian system 
\eqref{l minus one i l 1}. We begin by recalling that the  solution
of \eqref{l minus one i l 1} is given by $\iota_-(\rho(t_l)) =
h_-(t_l)^{-1} \iota_-(\rho(0)) h_-(t_l)$, where $h_-(t_l) = 
e^{-t_l\rho_0(0)^{l-1}}$, $\rho(0) =
\rho_0(0) + \rho_-(0) \in L ^1_-$ is the initial value of the
variable $\rho$ at $t _l = 0$, $\rho_0 \in L^1_0$ a diagonal
operator, and $\rho_-$ a strictly lower triangular operator.
Therefore,
$h_-(t_l) h_-(t_m) = h_-(t_m)h_-(t_l)$ for any $l, m \in
\mathbb{N}$ and hence the product
\begin{equation}
\label{inf product}
h_-(t) := h_-(t_1, t_2, \dots) :
= \prod_{l=1}^{\infty}h_-(t_l) 
\end{equation}
is independent on the order of the factors and it exists as
an invertible bounded operator if we assume that $t: = (t_1,
t_2, \dots) \in \ell ^{\infty}_0$ which means that $t $ has
only finitely many non-zero elements. 

One also has
\begin{equation}
\label{aa condition}
h_-(t)^{-1}\mathcal{J}_{\nu_-}\big(\mathbf{q},
\mathbf{p}\big) h_-(t) = 
\mathcal{J}_{\nu_-}\left(\mathbf{q} +
\sum_{l=1}^{\infty} t_l\rho_0(0)^{l-1}, \mathbf{p}\right)
\qquad \text{for}
\qquad t \in \ell_0^{\infty},
\end{equation}
which shows that the flow in the coordinates $(\mathbf{q},
\mathbf{p})$ is described by a straight line motion in
$\mathbf{q}$ with $\mathbf{p}$ conserved. If this would be a
finite dimensional system, since $(\mathbf{q}, \mathbf{p})$ are
also Darboux coordinates (see \eqref{symplectic form} or
\eqref{symplectic form in coordinates}), we would say that they are
action-angle coordinates on $\mathcal{J}_{\nu_-} (\ell ^{\infty}
\times \ell ^1)$.

In infinite dimensions, even the definition of action-angle
coordinates presents problems. First, if the symplectic form is
strong, the Darboux theorem (that is, the symplectic  form is
locally constant) is valid; see the proof of Theorem 3.2.2 in
\cite{A-M}. Second, if the symplectic form is weak, which is our
case, the Darboux theorem fails in general, even if the manifold
is a Hilbert space; Marsden's classical counterexample can be
found and discussed in Exercise 3.2H of \cite{A-M}. Third, even if
one could show  in a particular case that the Darboux theorem
holds, there still is the problem of coordinates. In the case
presented above, the action-angle coordinates were constructed
explicitly. In general, on Banach weak symplectic manifolds this
may well be impossible.

\medskip

We return now to the systems described by the family of
integrals in involution $I _l^S$ given by \eqref{i l s}.
By Corollary \ref{involution corollary}(ii), substituting
$I_l^S$ into
\eqref{equ on s one}, applying $\iota_S $ to 
\eqref{equ on s one}, and
using \eqref{derivative of i l}, yields the family of
Hamilton equations on
$L^1_S$ 
\begin{equation}
\label{l s one i l}
\frac{\partial \iota_S(\sigma)}{\partial t_l} 
= \left(P ^1_- + P ^1_0 + T \circ P ^1_-\right)
\left[\left(P^{\infty}_+ + P^{\infty}_0 + T \circ
P^{\infty}_-\right) \left([\iota_S(\sigma)]^{l-1}\right),
\iota_S(\sigma)\right]
\end{equation}
or, equivalently, in Lax form
\begin{equation}
\label{l s one i l 1} 
\frac{\partial \iota_S(\sigma)}{\partial t_l}
= - \left[(P^{\infty}_-
- T \circ P ^{\infty}_-)
\left([\iota_S(\sigma)]^{l-1}\right),
\iota_S(\sigma) \right],
\end{equation}
where $t_l$ denotes the time parameter for the
$l$th flow. 

From \eqref{l s one i l} it follows that the solution of
this equation can be written in terms of the coadjoint
action of the Banach Lie group $GL ^{\infty}_+$ on the dual
$L ^1_S$ of its Lie algebra. More precisely, the solution is
necessarily of the form 
\begin{equation}
\label{solution l minus one i l plus}
\iota_S(\sigma(t_l)) = 
\left(\operatorname{Ad}^S \right)^\ast_{\bar{g}_+(t_l) ^{-1}}
\sigma(0)
 = \left(P ^1_- + P ^1_0 + T
\circ P ^1_-\right) \left(\bar{g}_+(t_l)
\iota_S(\sigma(0))
\bar{g}_+(t_l)^{-1}\right)
\end{equation}
for some smooth curve $\mathbb{R} \ni t_l \mapsto
\bar{g}_+(t_l) \in GL^{\infty}_+$ and $\sigma(0)$ the
initial condition for $t_l = 0$. 

On the other hand, the solution of \eqref{l s one i l 1} is
\begin{equation}
\label{solution l minus one i l o}
\iota_S(\sigma(t_l)) = g_S(t_l)^T
\iota_S(\sigma(0)) g_S(t_l),
\end{equation} 
where $\mathbb{R} \ni t_l \mapsto g_S(t_l) \in O^{\infty}$
is a smooth curve that will be determined in the next
proposition by the same method as in the finite dimensional  case
(see, e.g., \cite{F,K,N,S}).

\begin{proposition}
\label{solution method}
Assume that we have the decomposition (we set here $t = t_l$)
\begin{equation}
\label{toda start}
e^{t [\iota_S(\sigma(0))]^{l-1}} = g_S(t)g_+(t)
\end{equation}
for $g_S(t) \in O ^{\infty}$ and $g_+(t) \in
GL^{\infty}_+$. Then
\begin{align}
\label{toda dec sol}
\iota_S(\sigma(t)) : &= g_S(t)^T
[\iota_S(\sigma(0))] g_S(t) =
g_+(t)[\iota_S(\sigma(0))] g_+(t)^{-1}
\end{align} 
is the solution of  
\eqref{l s one i l 1} with initial condition
$\iota_S(\sigma(0))$. 
\end{proposition}

\noindent \textbf{Proof.} 
To prove the first equality in 
\eqref{toda dec sol}, use
\eqref{toda start} to get
\[
g_S(t) = e^{t [\iota_S(\sigma(0))]^{l-1}} g_+(t)^{-1} 
\]
and hence
\[
g_+(t) e^{-t[\iota_S(\sigma(0))]^{l-1}}
 [\iota_S(\sigma(0))]
 e^{t [\iota_S(\sigma(0))]^{l-1}}g_+(t)^{-1}
= g_+(t)[\iota_S(\sigma(0))]g_+(t)^{-1}
\] 
since $\iota_S(\sigma(0))$ commutes with
$e^{t [\iota_S(\sigma(0))]^{l-1}}$.

 Let $\iota_S(\sigma(t)) : = 
g_S(t)^{-1}[\iota_S(\sigma(0))] g_S(t)$. Taking the
time derivative of
\eqref{toda start} and multiplying on the right by
$g_S(t)^{-1}$ and on the left by $g_+(t)^{-1}$ we get
\[
[\iota_S(\sigma(t))]^{l-1} = g_S(t)^{-1} \dot{g}_S(t) +
\dot{g}_+(t)g_+(t)^{-1} 
\] 
which is equivalent to the equations
\begin{align}
\label{skew g s}
g_S(t)^{-1} \dot{g}_S(t) & = \left(P ^{\infty}_- - T \circ
P^{\infty}_-\right)
\left([\iota_S(\sigma(t))]^{l-1}\right) \\
\label{plus g}
\dot{g}_+(t) g_+(t)^{-1} & = \left(P^{\infty}_+ +
P^{\infty} _0 + T \circ P^{\infty}_-\right)
\left([\iota_S(\sigma(t))]^{l-1}\right).
\end{align}
Therefore
\begin{align*}
&\frac{d}{dt}\iota_S(\sigma(t)) 
= - g_S(t)^{-1} \dot{g}_S(t) g_S(t)^{-1}
[\iota_S(\sigma(0))] g_S(t) +
g_S(t)^{-1}[\iota_S(\sigma(0))]
\dot{g}_S(t)\\
&\quad = - \left(P ^{\infty}_- - T \circ
P^{\infty}_-\right)
\left([\iota_S(\sigma(t))]^{l-1}\right)
\iota_S(\sigma(t)) + \iota_S(\sigma(t))
\left(P ^{\infty}_- - T \circ P^{\infty}_-\right)
\left([\iota_S(\sigma(t))]^{l-1}\right) \\
& \quad = - \left[(P^{\infty}_-
- T \circ P ^{\infty}_-)
\left([\iota_S(\sigma)]^{l-1}\right),
\iota_S(\sigma) \right]
\end{align*}
which is \eqref{l s one i l 1}.
\quad $\blacksquare$
\medskip

This proposition shows that the solution 
\eqref{toda dec sol} of the system \eqref{l s one i l 1}
could be expressed using the analogue of the
Iwasawa decomposition $GL ^{\infty} = O ^{\infty} \cdot
GL_0 ^{\infty}\cdot GI_{+,1} ^{\infty}$ for the Banach Lie 
group $GL^{\infty}$.  To our knowledge, there is no proof of this
decomposition and there could be technical difficulties that may 
even render it impossible. However, see the appendix in
\cite{neeb} for the polar decomposition theorem. 

Note also that \eqref{toda dec sol} produces a smooth curve
$g_+(t) \in GL ^{\infty}_+$ satisfying  
\eqref{solution l minus one i l plus} even without the
projection operator in that formula. This follows also
directly from \eqref{toda dec sol} and 
\eqref{concrete s coadjoint action}.

\medskip

The previous general considerations involving
Proposition \ref{induction proposition}, imply that the families
of flows given by \eqref{equ on l one minus} or 
\eqref{equ on s one} and, in particular by
\eqref{l minus one i l 1} or \eqref{l s one i l 1}, not only
preserve the symplectic leaves of $L^1_-$ and $L^1_S$, but
also the filtrations \eqref{blp minus filtration} and 
\eqref{blp s filtration}, respectively. This remark has some
 important consequences which we  discussed below.
\medskip

We turn now  to the study of Hamiltonian systems induced on the
filtrations \eqref{blp minus filtration} and 
\eqref{blp s filtration}.
A {\bfi $k $-diagonal Hamiltonian system\/} is, by definition, a
Hamiltonian system on $\left(L^1_{-,k}, \{ \cdot , \cdot \}_k
\right)$. Since the map $\Phi_{S,-,k}: \left(L ^1_{-,k}, \{
\cdot , \cdot \}_k \right) \rightarrow \left(L ^1_{S,k}, \{
\cdot , \cdot \}_{S,k} \right)$ introduced at the end of \S
\ref{section: induced and coinduced from ell one} is a Banach
Lie-Poisson space isomorphism, we can regard $k $-diagonal
Hamiltonian systems as being defined also on $\left(L ^1_{S,k},
\{\cdot , \cdot \}_{S,k} \right)$. 
From \eqref{algebra coordinate coadjoint minus k action},
Hamilton's equations on $\left(L^1_{-,k}, \{ \cdot , \cdot \}_k
\right)$ defined by an arbitrary function $h_k \in C
^{\infty}(L^1_{-,k})$ are given by
\begin{equation}
\label{hamiltonian v f minus k action}
\frac{d}{dt} \rho_j= -
\sum_{l=j}^{k-1}\left(\tilde{s}^{l-j} \left(\rho_l 
 \frac{\delta h_k}{\delta\rho_{l-j}}\right)
- \rho_l s^{j}\left(\frac{\delta h_k}{\delta\rho_{l-j}}\right)
\right) \quad \text{for} \quad j = 0,1,2,\dots , k-1.
\end{equation}

Note that for all $n>k$ (including $n = \infty $), any $h_k \in
C^{\infty}(L^1_{-,k})$ can be smoothly extended to $h_n : =
h_k\circ \pi_{kn}\in C ^{\infty}(L^1_{-,n})$, where $\pi_{kn}:
L^1_{-,n} \rightarrow L ^1_{-,k}$ is the projection that
eliminates the last lower $n-k $ diagonals of an operator in
$L^1_{-,n} : = \oplus_{i=-n+1}^0 L^1_i$. Conversely, any $h_n \in
C^{\infty} (L ^1_{-,n})$ gives rise to a smooth function $h_k : =
h _n\circ \iota_{nk} \in C ^{\infty}(L ^1_{-,k})$, where
$\iota_{nk}: L^1_{-,k} \hookrightarrow L ^1_{-,n}$ is the natural
inclusion. Since the flow defined  by $h \in C^\infty(L^1_-)$
preserves the filtration \eqref{blp minus filtration} (see 
Proposition \ref{induction proposition}) it follows that if the
initial condition $\rho(0) \in L ^1_{-,k}$ its trajectory is
necessarily contained in $L ^1_{-,k}$. This means that in order to
solve the system
\eqref{hamiltonian v f minus k action} for a given $k \in
\mathbb{N}$, it suffices to solve the Hamiltonian system given by
the extension of $h_k $ to
$\left(L ^1_-, \{\cdot , \cdot \}_- \right)$ for initial
conditions in $L^1_{-,k}$.  

Let us now specialize the functions $h_k \in
C^{\infty}(L^1_{-,k})$ and $f_k \in C^{\infty}(L ^1_{S,k})$
to 
\begin{align}
\label{i l minus k}
I^{-,k}_l (\rho) &: = I_l^-\left(\iota_{-,k} (\rho)\right)
= I_l \left((\iota_- \circ \iota_{-,k})( \rho)  \right) \quad
\text{for} \quad \rho \in L^1_{-,k}\\ 
\label{i l s k}
I^{S,k}_l(\sigma) &: = I_l^S
\left(\iota_{S,k}(\sigma)\right) = I_l\left((\iota_S
\circ \iota_{S,k})(\sigma)\right)
\quad \; \text{for} \quad\sigma \in L^1_{S,k},
\end{align}
respectively, where $\iota_{-,k}: L^1_{-,k} \hookrightarrow
L^1_-$ and $\iota_{S,k}: L^1_{S,k} \hookrightarrow L ^1_S$ are
the inclusions. Note that since $I^{S,k}_l \circ \Phi_{S,-,k}
\neq I_l^{-,k}$, the dynamics induced by the functions
$I^{-,k}_l$ and $I^{S,k}_l$ are different in spite of the
fact that the Poisson structures on $L^1_{-,k}$ and
$L ^1_{S,k}$ are isomorphic. Therefore, we see that  on has the 
family of Hamiltonian systems indexed by $k \in \mathbb{N}$ which
have an infinite  number of integrals in involution indexed
by $l \in \mathbb{N}$. For
$k=2$ the system is the semi-infinite Toda lattice.
Therefore, the  {\bfi $k$-diagonal semi-infinite Toda systems\/} are
defined to be the Hamiltonian systems on $L^1_{S,k}$
associated to the functions
$I^{S,k}_l$, $l \in \mathbb{N}$.

An important consequence of the fact that the Poisson
brackets on $L^1_{-,k}$ and $L^1_{S,k}$ are induced is that
the method of solution of the corresponding Hamilton equations
for $I_l^{-,k}$ and $I_l^{S,k}$, respectively, can be
obtained by solving  these equations on $L ^1_-$ and $L
^1_S$ respectively. Namely, it suffices to work with the
equations of motion 
\eqref{l minus one i l 1}  and \eqref{l s one i l 1} with
initial conditions $\rho(0) \in L ^1_{-,k}$ and
$\sigma(0) \in L^1_{S,k}$, respectively, and use
Proposition \ref{solution method}.
 We shall do this
in the rest of the paper for a special case related to the
semi-infinite Toda system.

\section{The bidiagonal case}
\label{sect: bidiagonal case}

In this section we shall study in great detail the
bidiagonal case consisting of operators that have only two
non-zero diagonals: the main one and the lower $k-1$
diagonal. The results obtained  in this section will be used
later to give a rigorous functional analytic formulation of
the integrability of the semi-infinite Toda lattice.

\paragraph{The coordinate description of the bidiagonal
subcase.} Due to their usefulness in the study of
the Toda lattice, we shall express in coordinates several
formulas from \S \ref{section: induced and coinduced from
ell one} adapted to the subalgebra $I^{\infty}_{+, 0,k-1}
\subset L ^{\infty}_{+,k}$, $k \geq 2$, consisting of
bidiagonal elements
\begin{equation}
\label{x bidiagonal element}
x : = x_0 + x_{k-1} S^{k-1} = \sum_{i =0}^{\infty}\left(
x_{0,ii} |i \rangle \langle i| + x_{k-1,ii} |i \rangle
\langle i+k-1 | \right),
\end{equation}
where $x_0, x_{k-1}$ are diagonal operators whose
entries are  given by the sequences
$\{x_{0,ii}\}_{i=0}^{\infty}$,
$\{x_{k-1,ii}\}_{i=0}^{\infty} \in \ell^{\infty}$,
respectively.  The subalgebra $I^{\infty}_{+, 0,k-1}$ of
$L^{\infty}_{+,k}$ is hence formed by upper triangular
bounded operators that have only two non-zero diagonals,
namely the main diagonal and the strictly upper $k-1$
diagonal. 

The predual of $I^{\infty}_{+, 0,k-1}$ is $I^1_{-,0,k-1}$
which consists of lower triangular trace class operators
having only two non-vanishing diagonals, namely the main 
one and the strictly lower $k-1$ diagonal ($k \geq 2$), that
is, they are of the type
\begin{equation}
\label{rho bidiagonal element}
\rho = \rho_0 + (S^{k-1})^T \rho_{k-1} = \sum_{i=0}
^{\infty} \left(\rho_{0,ii} |i \rangle \langle i| +
\rho_{k-1,ii}|i+k-1 \rangle \langle i| \right),
\end{equation}
where $\rho_0$ and $\rho_{k-1}$ are diagonal operators whose
entries are  given by the sequences
$\{\rho_{0,ii}\}_{i=0}^{\infty}$,
$\{\rho_{k-1,ii}\}_{i=0}^{\infty}  \in \ell^1$,
respectively. 
The Banach Lie subgroup $GI^{\infty}_{+,0,k-1}$ of
$GL ^{\infty}_{+,k}$ whose Banach Lie algebra is
$I^{\infty}_{+,0,k-1}$ has elements given by
\begin{equation}
\label{bidiagonal k group element}
g = g_0 + g_{k-1}S^{k-1} = \sum_{i =0}^{\infty}\left(
g_{0,ii} |i \rangle \langle i| + g_{k-1,ii} |i \rangle
\langle i+k-1 | \right),
\end{equation}
where $g_0$ and $g_{k-1}$ are diagonal operators whose
entries are  given by the sequences
$\{g_{0,ii}\}_{i=0}^{\infty}$,
$\{g_{k-1,ii}\}_{i=0}^{\infty} \in \ell^{\infty}$,
respectively,  and the sequence
$\{g_{0,ii}\}_{i=0}^{\infty}$ is bounded below by a
strictly positive number (that depends on $g_0$). 

The product of $g,h \in GI^{\infty}_{+,0,k-1}$ in
$GL ^{\infty}_{+,k}$ is given by
\begin{align}
\label{composition k bidiagonal}
g \circ_k h &= g_0h_0 + (g_0h_{k-1} + g_{k-1}
s^{k-1}(h_0))S^{k-1} \nonumber\\
&=  \sum_{i=0}^\infty g_{0,ii}h_{0,ii}|i\rangle
\langle i| + \sum_{i=0}^\infty \left(g_{ii}h_{k-1,ii}
+ g_{k-1,ii} h_{0,i+k-1, i+k-1}\right)|i\rangle \langle
i+k-1|.
\end{align}
and the  inverse of $g$ in $GL ^{\infty}_{+,k}$ is given by
\begin{equation}
\label{inverse in k bidiagonal}
g^{-1} = g_0^{-1} - g_0 ^{-1}g_{k-1} s^{k-1}(g_0 ^{-1})
S^{k-1} = 
\sum_{i=0}^\infty \frac{1}{g_{0,ii}}|i\rangle
\langle i| - \sum_{i=0}^\infty \frac{g_{k-1,ii}}{g_{0,ii}
g_{0,i+k-1,i+k-1}}|i\rangle \langle i+k-1|.
\end{equation}
The Lie bracket of $x, y \in I^\infty_{+,0,k-1}$ has the
expression 
\begin{align}
\label{Lie bracket on k bidiagonals}
&[x, y]_{k}  = \big(x_{k-1}(s^{k-1}(y_0) - y_0) -
y_{k-1}(s^{k-1}(x_0) -x_0) \big)S^{k-1}
\nonumber\\
 & = \sum_{i=0}^\infty \big(x_{k-1,ii}(y_{0,i+k-1,
i+k-1} - y_{0,ii}) - y_{k-1,ii}(x_{0,i+k-1, i+k-1}
 - x_{0,ii})\big)|i \rangle \langle i+k-1|.
\end{align}

The group coadjoint action
$\left(\operatorname{Ad}^{+,k} \right)^\ast_{g ^{-1}}:
I^1_{-,0, k-1} \rightarrow I^1_{-,0, k-1}$ for
$g : = g_0 + g_{k-1}S^{k-1}\in GI^\infty_{+,0, k-1} \subset
GL ^{\infty}_{+,k}$ and Lie algebra coadjoint action
$(\operatorname{ad}^{+,k})_x ^\ast: I^1_{-,0, k-1}
\rightarrow I^1_{-,0,k-1}$, for $x: =x_0 + x_{k-1}S^{k-1} \in
I^\infty_{+,0,k-1}
\subset L^{\infty}_{+,k}$ are given by
\begin{align}
\label{group k coadjoint action}
\left(\operatorname{Ad}^{+,k}\right)_{g^{-1}}^\ast \rho 
& = \rho_0 + g_0^{-1} g_{k-1} \rho_{k-1} 
-\tilde{s}^{k-1}\left(g_0^{-1} g_{k-1} \rho_{k-1}
\right)\left(\mathbb{I} - \sum_{j=0}^{k-2}p_j \right) 
\nonumber \\ 
& \qquad  + 
\left(S^T\right)^{k-1} s^{k-1}(g_0) g_0 ^{-1}
\rho_{k-1}     \nonumber\\  
&= \sum_{i=0}^\infty
\left(\rho_{0,ii} + \rho_{k-1,ii}
\frac{g_{k-1,ii}}{g_{0,ii}} - \rho_{k-1,ii}
\frac{g_{k-1,ii}}{g_{0,i-k+1,i-k+1}}
\right)|i
\rangle \langle i | \nonumber \\
& \qquad + \sum_{i=0}^\infty \rho_{k-1,ii}
\frac{g_{0, i+k-1, i+k-1}}{g_{0,ii}} |i + k- 1 \rangle
\langle i |
\end{align}
and
\begin{align}
\label{algebra k coadjoint action}
\left(\operatorname{ad}^{+,k}\right)_x ^\ast \rho 
& =  \tilde{s}^{k-1}(\rho_{k-1} x_{k-1}) - \rho_{k-1}
x_{k-1} +  \left(S^T\right)^{k-1} \rho_{k-1}(x_0 -
s^{k-1}(x_0))   \nonumber\\ 
&= \sum_{i=0}^\infty(\rho_{k-1,ii} x_{k-1,ii} -
\rho_{k-1,ii} x_{k-1,ii} )|i \rangle \langle i|
\nonumber   \\  
& \qquad + \sum_{i=0}^\infty
\rho_{k-1,ii} (x_{0,ii} - x_{0,i+k-1, i+k-1}) |i + k - 1
\rangle \langle i |,
\end{align}
where $\rho: = \rho_0 + (S^T)^{k-1}\rho_{k-1} \in
I^1_{-,0,k-1}$.
\medskip

Since $\left(I ^1_{-,0,k-1} \right)^\ast = I
^{\infty}_{+,0,k-1}$ and the duality pairing  is given by
the trace of the product, it follows that the Lie-Poisson
bracket and its associated Hamiltonian vector field on $I
^1_{-,0,k-1}$ are given by
\begin{align}
\label{PB for k bidiagonal}
&\{f, h\}_{0, k-1}(\rho)  \nonumber  \\
& \qquad = 
\operatorname{Tr} \left[\rho_{k-1}
\left(\frac{\partial f}{\partial \rho_{k-1}} \left(s^{k-1}
\left(\frac{\partial h}{\partial \rho_0} \right) -
\frac{\partial h}{\partial \rho_0} \right) 
- \frac{\partial h}{\partial \rho_{k-1}} \left(s^{k-1}
\left(\frac{\partial f}{\partial
\rho_0} \right) - \frac{\partial f}{\partial \rho_0}
\right) \right) \right]
\nonumber\\
& \qquad = \sum_{i=0}^\infty \rho_{k-1,ii} 
\left[\frac{\partial f }{\partial \rho_{k-1,ii}}
\left(\frac{\partial h}{\partial \rho_{0,i+k-1,i+k-1}} -
\frac{\partial h}{\partial \rho_{0,ii}} \right)
\right. \nonumber\\
& \qquad \qquad \qquad \qquad \qquad   \left.
- \frac{\partial h}{\partial \rho_{k-1,ii}}
\left(\frac{\partial f}{\partial \rho_{0, i+k-1,i+k-1}} -
\frac{\partial f}{\partial \rho_{0,ii}} \right) \right]
\end{align}
and 
\begin{align}
\label{Hamiltonian vector field for k bidigonal}
X_h^{0,k-1}(\rho) 
& = \operatorname{Tr} \left[\rho_{k-1}
\left(s^{k-1}
\left(\frac{\partial h}{\partial \rho_0} \right) - 
\frac{\partial h}{\partial \rho_0} 
\right)\frac{\partial}{\partial \rho_{k-1}} - 
\frac{\partial h}{\partial \rho_{k-1}} \left(s^{k-1} \left(
\frac{\partial}{\partial \rho_0} \right) - 
\frac{\partial}{\partial \rho_0} \right)
 \right]
\nonumber\\
& = \sum_{i=0}^\infty \rho_{k-1,ii} \left[
\left(\frac{\partial h}{\partial \rho_{0,i+k-1,i+k-1}} -
\frac{\partial h}{\partial \rho_{0,ii}}
\right)\frac{\partial}{\partial \rho_{i+k-1,i}} 
\right. \nonumber\\
& \qquad \qquad \qquad \qquad    \left.
- \frac{\partial h}{\partial \rho_{k-1,ii}}
\left(\frac{\partial}{\partial \rho_{0,i+k-1,i+k-1}} -
\frac{\partial}{\partial \rho_{0,ii}}\right) \right] 
\end{align}
for $f,h \in C^\infty(I^1_{-,0,k-1})$. Like in 
\S\ref{section: symplectic induction}, in 
\eqref{Hamiltonian vector field for k bidigonal}  we have used the
standard coordinate conventions from finite dimensions to write a
vector field. The precise meaning of the symbols $\partial /
\partial\rho_{k-1} =
\{\partial/
\partial \rho_{i+k-1,i}\}_{i=0} ^{\infty}$ and $\partial/
\partial \rho_0 = \{\partial/
\partial_{0,ii}\} _{i=0} ^{\infty}$ is that they form the Schauder
basis of the tangent space $T _\rho I ^1_{-,0,k-1}$ corresponding to
the Schauder basis $\{|i+k-1 \rangle \langle i|, |i
\rangle \langle i|\}_{i=0} ^{\infty}$ of $I ^1_{-,0,k-1}$.
Thus Hamilton's equations in terms of diagonal operators are 
\begin{align}
\label{first Hamilton equations for k bidiagonal in diag}
\frac{d}{dt}\rho_0 &= \rho_{k-1} 
\frac{\partial h}{\partial \rho_{k-1}} - \tilde{s}^{k-1}
\left(\rho_{k-1} \frac{\partial h}{\partial
\rho_{k-1}} \right)  \\
\label{second Hamilton equations for k bidiagonal in diag}
\frac{d}{dt}\rho_{k-1} &= \rho_{k-1} \left(s^{k-1}
\left(\frac{\partial h}{\partial \rho_0} \right) - 
\frac{\partial h}{\partial \rho_0} 
\right)
\end{align}
or, in coordinates, for $i\in \mathbb{N} \cup \{0\}$, $k \geq
2$,
\begin{align}
\label{first Hamilton equations for k bidiagonal}
\frac{d}{dt}{\rho}_{0,ii} &=  
\rho_{k-1,ii}\frac{\partial h}{\partial
\rho_{k-1,ii}} -\rho_{k-1,ii} \frac{\partial
h}{\partial \rho_{k-1,ii}}
\\
\label{second Hamilton equations for k bidiagonal}
\frac{d}{dt}{\rho}_{k-1,ii} &= 
\rho_{k-1,ii}\left(\frac{\partial h}{\partial
\rho_{0,i+k-1,i+k-1}} - \frac{\partial h}{\partial
\rho_{0,ii}}  \right).
\end{align}

\paragraph{Structure of the generic coadjoint  orbit.}
By a generic coadjoint orbit we will understand the orbit
\[
\mathcal{O}_ \nu : = \left\{
\left(\operatorname{Ad}^{+,k} \right) ^\ast_{g^{-1}}
\nu \, \big{|}\, g \in GI^{\infty}_{+,0,k-1}
\right\},
\]
through the element $\nu = \nu_0 +
\left(S^T\right)^{k-1} \nu_{k-1} \in I^1 _{-,0, k-1}$ such
that $\nu_{k-1,ii}\neq 0 $ for $i = 0, 1, 2, \dots$. 

Let us denote by $GL ^{\infty,k-1}_0$ the Banach Lie
subgroup of $(k-1)$-periodic elements of $GL ^{\infty}_0$,
that is, $g_0 \in GL ^{\infty,k-1}_0$ if and only if
$s^{k-1}(g _0) = g _0$.  Denote by $L ^{\infty,k-1}_0$ the
Banach Lie algebra of $GL^{\infty,k-1}_0$. 

\begin{proposition}
\begin{itemize}
\item[{\rm \textbf{(i)}}] One has the following equalities
\begin{equation}
\label{centralizer number}
Z(GI^{\infty}_{+,0,k-1}) =
\left(GI^{\infty}_{+,0,k-1} \right)_ \nu = GL_0^{\infty,k-1},
\end{equation}
where $Z(GI^{\infty}_{+,0,k-1}) $ is the  center of
$GI^{\infty}_{+,0,k-1}$ and $\left(GI^{\infty}_{+,0,k-1}
\right)_ \nu $ is the stabilizer of the generic element $\nu
\in I^1 _{-,0, k-1}$.
\item[{\rm \textbf{(ii)}}] The generic orbit
\begin{equation}
\label{orbit group}
\mathcal{O}_\nu \cong
GI^{\infty}_{+,0,k-1}/GL_0^{\infty,k-1} 
\end{equation} 
is a Banach Lie group.
\item[{\rm \textbf{(iii)}}] One has the relation
\begin{equation}
\label{orbit number}
\mathcal{O}_ \nu = \nu_0 + \mathcal{O}_{(S^T)^{k-1}
\nu_{k-1}}
\end{equation}
between the coadjoint orbits through $\nu = \nu_0 +
(S^T)^{k-1} \nu_{k-1}$ and through
$\left(S^T\right)^{k-1}\nu_{k-1}$. 
\end{itemize}
\end{proposition}

\noindent \textbf{Proof.} Part \textbf{(i)} follows from a
direct verification. Since $GL_0^{\infty,k-1} $ is a normal
Banach Lie  group of $GI^{\infty}_{+,0,k-1}$ the quotient
$GI^{\infty}_{+,0,k-1}/GL_0^{\infty,k-1}$ is also a
Banach Lie group (see \cite{Bou1972}). This proves
\textbf{(ii)}. Part \textbf{(iii)} follows from \eqref{group
k coadjoint action}. \quad $\blacksquare$
\smallskip

We conclude from \eqref{orbit number} that to describe any
$\mathcal{O}_\nu$ it suffices to study  coadjoint orbits
through the $(k-1)$-lower diagonal elements, $k \geq 2$.

Since the Banach Lie group $GI_{+,0, k-1}^{\infty}$ and
the generic element $\nu 
\in I^1 _{-,0, k-1}$ satisfy all the hypotheses of
Theorems 7.3 and 7.4 in \cite{OR} we conclude:
\begin{itemize}
\item The map $\iota_{\nu}:
GI^{\infty}_{+, 0, k-1}/ GL_0^{\infty,k-1}
\rightarrow I^1 _{-,0, k-1}$ given by
$\iota_{\nu}([g]): =
\left(\operatorname{Ad}^{+,k} \right)^\ast_{g^{-1}}
\nu$  is a weak injective  immersion. This means
that its derivative is injective but no conditions on the
closedness of its range or the fact that it splits are
imposed. The map $\iota_ \nu$ is not an immersion
as we now show by using Theorem 7.5 in \cite{OR}.

Since the coadjoint stabilizer Lie algebra
$\left(I^{\infty} _{+,0, k-1}\right)_{\nu}$ is equal to the
center
\[
Z(I^{\infty} _{+,0, k-1}) = \left\{x = x_0 + x_{k-1}
S^{k-1}
\in I^{\infty}_{+,0,k-1}  \mid s^{k-1}(x_0) = x_0,\, 
x_{k-1} = 0 \right\}
\]
it follows that its annihilator is
\begin{align*}
\left(\left(I^{\infty} _{+,0, k-1}\right)_{\nu}\right) ^
\circ  &= \left\{ \rho = \rho_0 + (S^T)^{k-1}\rho_{k-1} \in
I^1_{-,0, k-1} \mid \operatorname{Tr}(x_0\rho_0) = 0,
\right. \\ 
&\left. \qquad \qquad \qquad \qquad 
 \text{ for all } x_0 \in L ^{\infty} _0  
\text{ such that } s^{k-1}(x_0) = x_0 \right\}.
\end{align*} 
Because
\begin{align*}
\operatorname{Tr}\left(x_0\left(\left(\operatorname{ad}^{+,k}
\right)^\ast_x \nu \right)_0\right) 
= \operatorname{Tr} 
\left(x_0\left(\operatorname{ad}^{+,k}
\right)^\ast_x \nu \right) 
= \operatorname{Tr} \left([x_0, x]_k \nu \right)= 0 
\end{align*}
for any $x_0 \in Z(I^{\infty}_{+,0,k-1})$ and any $x \in
I^{\infty}_{+,0,k-1}$, we have $S_{\nu}  \subset
\left(\left(I^{\infty} _{+,0, k-1}\right)_{\nu}\right)
^\circ$, where  $S_{\nu} : =
\left\{\left(\operatorname{ad}^{+,k} 
\right)^\ast_x \nu \mid x
\in I_{+,0,k-1}^{\infty} \right\}$ is the characteristic
subspace of the Banach Lie-Poisson structure of
$I_{-,0,k-1}^1$ at $\nu$. Moreover, the bounded operator
$K_{\nu}: x \in  I_{+,0,k-1} ^{\infty}
\mapsto \left(\operatorname{ad}^{+,k} \right)^\ast_x
\nu \in I_{-,0,k-1}^1$ has non-closed range
$\operatorname{im}K_{\nu} = S_{\nu}$ and thus the
inclusion $S_{\nu}  \subset \left(\left(I^{\infty} _{+,0,
k-1}\right)_{\nu}\right) ^\circ$ is strict. To see that
the range of $K_{\nu}$ is not closed, one uses the Banach
space isomorphisms $I_{-,0,k-1}^1 \cong \ell^1\times
\ell^1$ and $I_{+,0,k-1}^{\infty} \cong \ell ^{\infty}
\times
\ell^{\infty}$ and shows that the two components of
$K_{\nu}$ are both bounded linear operators with
non-closed range. Therefore, since Theorem 7.5 in \cite{OR}
states that $\iota_ \nu $ is an immersion if and
only if $S_{\nu}  =
\left(\left(I^{\infty} _{+,0, k-1}\right)_{\nu}\right)
^\circ$,
this argument shows that
$\iota_{\nu}$ is only a weak immersion.
\item The quotient space
$GI^{\infty}_{+, 0, k-1}/ GL_0^{\infty,k-1}$ is a weak
symplectic Banach manifold relative to the closed two-form
\begin{align}
\label{weak symplectic form on k quotient}
&\omega_{\nu}([g])(T_g\pi(g \circ_k x),T_g\pi(g
\circ_k y)) =
\operatorname{Tr}(\nu [x,y]_k) \nonumber \\
& =\sum_{i=0}^{\infty}\nu_{k-1,ii}\big( 
x_{k-1,ii}(y_{0,i+k-1,i+k-1} - y_{ii}) - 
y_{k-1,ii}(x_{0,i+k-1,i+k-1} - x_{0,ii}) \big), 
\end{align}
where $x,y \in I_{+,0,k-1}^{\infty}$, $g \in
GI_{+,0,k-1}^{\infty}$, $[g] :=
\pi(g)$, $\pi: GI^{\infty}_{+,0,k-1}
\longrightarrow GI^{\infty}_{+, 0, k-1}/
GL_0^{\infty,k-1}$ is the
canonical projection, and 
$T_g \pi: T_g GI^{\infty} _{+,0,k-1} \longrightarrow
T_{[g]}\left(GI^{\infty}_{+, 0, k-1}/ GL_0^{\infty,k-1}
\right)$ is its derivative  at $g$. In  this formula we have
used the fact that the value at $g$ of the left invariant
vector field $\xi_x$ on
$GI^{\infty}_{+,0,k-1}$ generated by $x$ is $g\circ _k x $.
\item Relative to the Banach manifold structure on
$\mathcal{O}_\nu$ making
$\iota_\nu: GI^{\infty}_{+, 0, k-1}/ GL_0^{\infty,k-1}
\longrightarrow \mathcal{O}_ \nu$ into a diffeomorphism, the
push forward of the weak symplectic  form \eqref{weak
symplectic form on k quotient} has the expression
\begin{align}
\label{weak symplectic form on k orbit}
&\omega_{ \mathcal{O}_ \nu}(\rho) \left(
\left(\operatorname{ad}^{+,k} \right)^\ast_x \rho, 
\left(\operatorname{ad}^{+,k} \right)^\ast_y \rho \right)
= \operatorname{Tr}(\rho[x,y]_k)
\nonumber \\
&=\sum_{i=0}^{\infty}\rho_{k-1,ii}\big( 
x_{k-1,ii}(y_{0,i+k-1,i+k-1} - y_{0,ii}) - 
y_{k-1,ii}(x_{0, i+k-1,i+k-1} - x_{0,ii}) \big), 
\end{align}
where $x,y \in I_{+,0,k-1}^{\infty}$ and $\rho \in
\mathcal{O}_ \nu$.
\end{itemize}

We shall express the pull back $\pi ^\ast
\omega_ \nu $ of the weak symplectic form $\omega_
\nu$ in terms of the diagonal operators represented by
$\{g_{0,ii}\}_{i=0}^{\infty} \in \ell^{\infty}$ and 
$\{g_{k-1,ii}\}_{i=0} ^{\infty} \in \ell^{\infty}$
defining the element $g \in GI^{\infty}_{+,0,k-1}$. If $x =
x_0 + x_{k-1} S^{k-1}$, $y = y_0 + y_{k-1} S^{k-1}
\in I^{\infty}_{+,0,k-1}$, and $\nu = \nu_0 + (S^T)^{k-1}
\nu_{k-1} \in I^1_{-,0,k-1}$, 
\eqref{weak symplectic form on k quotient} yields
\begin{align}
\label{k first intermediate weak symplectic form on orbit}
&(\pi^\ast \omega_ \nu)(g)\left(g \circ_k x , 
g \circ_k y \right) = 
\omega_{\nu}([g]) \left(T_g \pi(g
\circ _k x),  T_g \pi(g \circ _k y) \right)
= \operatorname{Tr}(\nu[x,y]_k)
\nonumber \\
&\qquad =\sum_{i=0}^{\infty}\nu_{k-1,ii}
\big(x_{k-1,ii}(y_{0,i+k-1,i+k-1} - y_{0,ii}) - 
y_{k-1,ii}(x_{0,i+k-1,i+k-1} - x_{0,ii}) \big), 
\end{align}
where $\nu_{k-1}$ has the diagonal entries $\{
\nu_{k-1,ii}\}_{i=0} ^{\infty}$. The left invariant
vector field $\xi_x$ on $GI_{+,0,k-1}^{\infty}$  generated
by $x$ has the expression
\begin{align*}
\xi_x = 
\sum_{i=0}^\infty g_{0,ii}x_{0,ii}\frac{\partial}{\partial
g_{0,ii}} + \sum_{i=0}^\infty
\left(g_{0,ii}x_{k-1,ii} + g_{k-1,ii}
x_{0,i+k-1,i+k-1}\right)\frac{\partial}{\partial
g_{k-1, ii}}. 
\end{align*}
The symbols $\{\partial/ \partial g_{0,ii},
\partial/ \partial g_{k-1,ii} \}_{i=0} ^{\infty}$ denote the
biorthogonal family in the tangent space $T_g I ^{\infty}_{+,0,k-1}$
corresponding to the standard biorthogonal family 
$\{|i\rangle \langle i|, |i\rangle \langle i+k-1|\}_{i=0} ^{\infty}$
in $I^{\infty}_{+,0,k-1}$. We shall use, as in finite dimensions, the
exterior derivative on real valued smooth functions, in
particular coordinates, to represent elements in the dual space. With
this convention, we have
\begin{equation}
\label{k pull back in coordinates}
\pi^\ast \omega_ \nu = \sum_{i=0} ^{\infty} \mathbf{d}
\log g_{0,ii} \wedge \mathbf{d}\left(\nu_{k-1,ii} 
\frac{g_{k-1,ii}}{g_{0,ii}} - \nu_{k-1,ii}
\frac{g_{k-1,ii}}{g_{0,i-k+1,i-k+1}}
\right),
\end{equation}
where, as usual, any element that has negative index is
set equal to zero.
To show this, we evaluate the right hand side of
\eqref{k pull back in coordinates} on $\xi_x $ and $\xi_y $
and observe that it equals the right hand side of
\eqref{k first intermediate weak symplectic form on orbit}. Note that
the computations make sense since $\nu_{k-1} \in \ell ^1$.

The action of the coadjoint isotropy subgroup 
$\left(GI ^{\infty}_{+, 0, k-1} \right)_ \nu =
GL_0^{\infty,k-1}$ on $GI^{\infty}_{+,0,k-1}$
is given by $g_{0,ii} \mapsto h_{0, ii}
g_{0,ii}$, $g_{k-1,ii} \mapsto
h_{0,ii} g_{k-1,ii}$, where $h_{0,ii} = h_{0, i+k-1, i+k-1}$. As
expected, the right hand side of \eqref{k pull back in coordinates}
is invariant under this transformation and its interior product with
any tangent vector to the orbit of the normal subgroup
$GL_0^{\infty,k-1}$ is zero. This shows, once again, that
\eqref{k pull back in coordinates} naturally descends to
the quotient group $GI ^{\infty}_{+, 0,k-1}/
GL_0^{\infty,k-1}$.

\medskip

In order to understand the structure of $\mathcal{O}_
\nu$, define the action $\alpha^k: GI_{+,0,k} ^{\infty}
\times L^1_{-k+1} \rightarrow  L^1_{-k+1}$ by 
\begin{equation}
\label{alpha action}
\alpha^k_g\left((S^T)^{k-1} \nu_{k-1} \right) : =
(S^T)^{k-1} s^{k-1}(g_0) g_0 ^{-1}\nu_{k-1}.
\end{equation} 
The projector $\delta^k: I^1_{-,0, k-1} \rightarrow 
L^1_{-k+1}$ defined by the splitting  $I^1_{-,0,k-1} = 
L^1_{-k+1}\oplus L^1_0 $ is a $GI^{\infty}_{+,0,k-1} $-
equivariant map relative to the coadjoint and the
$\alpha^k$-actions of $GI^{\infty}_{+,0,k-1}$, that is, the
diagram
\unitlength=5mm
\begin{center}
\begin{picture}(9,8.6)
    \put(1,7){\makebox(0,0){$I^1_{-,0,k-1}$}} 
    \put(9,7){\makebox(0,0){$I^1_{-,0,k-1}$}}
    \put(1,2){\makebox(0,0){$ L ^1_{-k+1}$}} 
    \put(9,2){\makebox(0,0){$ L ^1_{-k+1}$}}
    \put(1,6){\vector(0,-1){3}}
    \put(9,6){\vector(0,-1){3}} 
    \put(2.5,7){\vector(1,0){5}}
    \put(2,2){\vector(1,0){5.7}}
    \put(0.3,4.3){\makebox(0,0){$\delta^k$}}
    \put(9.6,4.4){\makebox(0,0){$\delta^k$}}
    \put(5,7.8){\makebox(0,0)
    {$\left(\operatorname{Ad}^{-,k}\right)^\ast_{g^{-1}}$}}
    \put(5,2.5){\makebox(0,0){$\alpha^k_g$}}
    \end{picture}
\end{center}
commutes for any $g \in GI^{\infty}_{+,0,k-1}$. We observe
that the stabilizer $GL_0^{\infty, k-1}$
of the $\alpha^k$-action does not depend on the choice of
the generic element $\left(S^T\right)^{k-1} \nu_{k-1} \in 
L^1_{k-1}$. The orbits of the coadjoint action of the
subgroup $GL_0^{\infty, k-1}$ on $(\delta^k)^{-1}((S^T)^{k-1}
\nu_{k-1})$ are of the form 
\[
\Delta_{ \nu_0, \nu_{k-1}} +
(S^T)^{k-1} \nu_{k-1} \subset (\delta^k)^{-1}((S^T)^{k-1}
\nu_{k-1}) \subset I^1_{-,0, k-1}, 
\]
where     
\begin{equation}
\label{orbits affine spaces}
\Delta_{ \nu_0, \nu_{k-1}}: =  \nu_0 + \operatorname{im}
\mathcal{N}_{\nu_{k-1}} \subset L ^1_0
\end{equation}
are affine spaces for each $\nu_0 \in L ^1_0 $ and
the linear operator
$\mathcal{N}_{\nu_{k-1}} : L _0 ^{\infty} \rightarrow L^1_0$
is defined by
\[
\mathcal{N}_{\nu_{k-1}}(g_{k-1}): = \nu_{k-1}
g_{k-1} + \tilde{s}(\nu_{k-1} g_{k-1})\left(\mathbb{I} -
\sum_{j=0}^{k-2}p_j\right).
\]

The orbits of the $\alpha^k$-action of
$GI^{\infty}_{+,0,k-1}$ on $L^1_{-k+1}$ are
\begin{equation}
\label{alpha orbits}
GI^{\infty}_{+,0,k-1} \cdot \left((S^T)^{k-1}
\nu_{k-1}\right) =
\{(S ^T)^{k-1} s^{k-1}(g_0) g_0 ^{-1} \nu_{k-1} \mid g_0 \in
GL _0 ^{\infty}\} = : \Delta_{\nu_{k-1}}.
\end{equation} 
Note that if $\Delta_{ \nu_{k-1}} = \Delta_{ \nu_{k-1} '}$
then $\operatorname{im} \mathcal{N}_{\nu_{k-1}}
=\operatorname{im} \mathcal{N}_{\nu_{k-1}'}$ and so
$\Delta_{ \nu_0,
\nu_{k-1} } =
\Delta_{
\nu_0, \nu'_{k-1}}$. These remarks show that the
coadjoint  orbit $\mathcal{O}_ \nu$ is diffeomorphic
to the product $\left(\nu_0+ \operatorname{im}
\mathcal{N}_{\nu_{k-1}}  \right)
\times
\Delta_{\nu_{k-1}}$ of the affine  space $\Delta_{
\nu_0, \nu_{k-1}} $ with the $\alpha^k$-orbit
$\Delta_{ \nu_{k-1}} $. This diffeomorphism does
not depend on the choice of $(S ^T)^{k-1} \nu_{k-1}' \in
\Delta_{\nu_{k-1}}$. Additionally, one identifies the set of
generic coadjoint orbits with the total space $\mathbb{L}_k$
of the vector bundle $\mathbb{L}_k \rightarrow
L^{\infty}_0 /\alpha^k(GL_{0}^{\infty})$, whose fiber
at $[\nu_{k-1}]$ is  $L ^1_0/\operatorname{im}
\mathcal{N}_{\nu_{k-1}}$. The vector space
$L^1_0/\operatorname{im} \mathcal{N}_{\nu_{k-1}}$ is not
Banach since
$\operatorname{im} \mathcal{N}_{\nu_{k-1}}$ is not closed in
$L ^1_0 $ because the operator $\mathcal{N}_ {\nu_{k-1}} : L
^{\infty}_0 \rightarrow L^1_0$ is compact. Consequently, the
bundle $\mathbb{L}_k \rightarrow L^{\infty}_0
/\alpha^k(GL_{0}^{\infty})$ does not have the structure of a
Banach vector bundle and does not have fixed typical fiber.

\paragraph{The momentum map.} 
Let us now study an important particular case of the  map
$\mathcal{J}_{ \nu_-}$ by taking in  \eqref{def of cal j} the
element $\nu_- = (S^T)^{k-1} \nu_{k-1} \in L ^1_{-k+1} \subset
L^1_{-,k}$. The map \eqref{def of cal j}, denoted in this
case $\mathcal{J}_{\nu_{k-1}}: \ell^{\infty} \times \ell^1
\rightarrow  I^1_{-,0,k-1}$, becomes
\begin{equation}
\label{flaschka map}
\mathcal{J}_{\nu_{k-1}}( \mathbf{q}, \mathbf{p}) =
\mathbf{p} + (S^T)^{k-1} \nu_{k-1} e ^{s^{k-1}( \mathbf{q}) -
\mathbf{q}}.
\end{equation}
Recall that  we identify
$\ell ^1$ with $L ^1_0 $ and  $\ell ^{\infty}$ with
$L^{\infty}_0 $. Having fixed
$(S^T)^{k-1} \nu_{k-1} \in L ^1_{-k+1}
$, define the action of  $GI^{\infty}_{+,0,k-1}$ on $\ell
^{\infty}\times \ell ^1$ by
\begin{align}
\label{sigma action}
&\boldsymbol{\sigma}_g ^{ \nu_{k-1}}( \mathbf{q}, \mathbf{p}): = 
\nonumber \\
& \left(\mathbf{q} + \log g_0, \mathbf{p}+
g_{k-1} g_0 ^{-1}
\nu_{k-1} e^{s^{k-1}(\mathbf{q}) - \mathbf{q}} -
\tilde{s}^{k-1}
\left(g_{k-1} g_0 ^{-1} \nu_{k-1} e^{s^{k-1}( \mathbf{q}) -
\mathbf{q}}\right) \right), 
\end{align}
where $g: = g_0 + g_{k-1} S^{k-1} \in GI^{\infty}_{+,0,k-1} $
and
$(\mathbf{q}, \mathbf{p}) \in \ell ^{\infty} \times \ell
^1$. The coordinate form of the action \eqref{sigma action}
is 
\begin{align}
\label{first sigma action}
q'_i &= q_i + \log g_{0,ii}\\
\label{second sigma action}
p'_i &= p_i + \frac{g_{k-1,ii}}{g_{0,ii}} \nu_{k-1,ii}
e^{q_{k+1} - q_k} - \frac{g_{k-1,ii}}{g_{0,k-1, k-1}}
\nu_{k-1,ii} e^{q_k - q_{k-1}}
\end{align}
for $i \in \mathbb{N}\cup \{0\}$. Using \eqref{first sigma
action} and \eqref{second sigma action} one shows that 
\[
\sum_{i=0} ^{\infty} p'_i \mathbf{d}q'_i = 
\sum_{i=0} ^{\infty} p_i \mathbf{d}q_i - \mathbf{d}Q,
\]
where the function $Q: \ell ^{\infty} \rightarrow
\mathbb{R}$ is given by 
\begin{equation}
\label{function Q}
Q(\mathbf{q}) : = \operatorname{Tr} \left(g_0 ^{-1}g_{k-1}
\nu_{k-1} e^{s^{k-1}( \mathbf{q}) - \mathbf{q}} \right)
= \sum_{i=0} ^{\infty} \frac{g_{k-1,ii}}{g_{0,ii}}
\nu_{k-1,ii} e^{q_{k+1} - q_k}. 
\end{equation}
Thus we see that $\omega$ is invariant relative to the
$\boldsymbol{\sigma}^{ \nu_{k-1}} $-action, that is , 
$\left(\boldsymbol{\sigma}_g
^{\nu_{k-1}} \right)^\ast \omega = \omega$ for any $g \in
GI^{\infty}_{+,0,k-1} $. 

\begin{proposition}
\label{Flaschka map theorem}
The smooth map
$\mathcal{J}_{\nu_{k-1}}: \ell^{\infty}\times
\ell ^1 \rightarrow I^1_{-,0,k-1}$ 
given by \eqref{flaschka map} is  constant on the
$\boldsymbol{\sigma}^{\nu_{k-1}}$-orbits of the subgroup
$GL_0^{\infty,k-1}$. In addition:
\begin{enumerate}
\item[{\rm \textbf{(i)}}] $\mathcal{J}_{\nu_{k-1}}$ is a
momentum map. More precisely, 
$\{f \circ \mathcal{J}_{\nu_{k-1}}, g \circ
\mathcal{J}_{\nu_{k-1}} \}_\omega =
 \{f, g \}_{0,k-1} \circ
\mathcal{J}_{\nu_{k-1}}$, for all $f,g \in C^{\infty}(I^1_{-,0,k-1})
$, where $\{ \cdot , \cdot \} _\omega$  is the canonical Poisson
bracket of the weak symplectic Banach  space $\left(\ell ^{\infty}
\times \ell ^1, \omega \right)$ given by 
\eqref{canonical Poisson bracket on ell} and $\{\,,\}_{0,k-1}$ is the
Lie-Poisson bracket on $I ^1_{-,0,k-1}$ given by 
\eqref{PB for k bidiagonal}.
\item[{\rm \textbf{(ii)}}] $\mathcal{J}_{\nu_{k-1}}$ is
$GI^{\infty}_{+,0,k-1}$-equivariant, that is,
$\mathcal{J}_{\nu_{k-1}} \circ
\boldsymbol{\sigma}^{ \nu_{k-1}}_{g} =
\left(\operatorname{Ad}^{-,k} \right)^\ast_{g^{-1}}
\circ \mathcal{J}_{\nu_{k-1}}$ for any
$g \in GI^{\infty}_{+,0,k-1}$. 
\item[{\rm \textbf{(iii)}}] One has $\mathcal{J}_{
\nu_{k-1}}(
\ell ^{\infty} \times \ell ^1) =
(\delta^k)^{-1}\left(\Delta_{
\nu_{k-1}} \right)$ and $ \mathcal{J}_{ \nu_{k-1}} ( \ell
^{\infty}\times \{0\}) = \Delta_{ \nu_{k-1}} $ and hence
$(\ell ^{\infty} \times \ell ^1)/ \boldsymbol{\sigma}^{
\nu_{k-1}} (GL_0^{\infty,k-1})
\cong \mathcal{J}_{\nu_{k-1}}( \ell ^{\infty} \times
\ell^1)$ consists of those coadjoint orbits which are
projected by  $\delta^k$ to the $\alpha^k$-orbit $\Delta_{
\nu_{k-1}} $.
\end{enumerate}
\end{proposition}

\noindent \textbf{Proof.} 
To prove \textbf{(i)}, let $f, g \in C^{\infty}(I^1_{-,0,k-1})$
and notice that 
\[
\frac{\partial(f \circ \mathcal{J}_{ \nu_{k-1}})}
{\partial \mathbf{q}} \in (L ^{\infty})^\ast \qquad \text{and}
\qquad
\frac{\partial(f \circ \mathcal{J}_{ \nu_{k-1}})}
{\partial \mathbf{p}} \in (L ^1)^\ast = L ^{\infty}
\]
because $\mathbf{q}\in L ^{\infty}$ and $\mathbf{p} \in L ^1$. 
However, by \eqref{flaschka map}, 
\begin{equation}
\label{partial f by partial  q}
\frac{\partial(f \circ \mathcal{J}_{\nu_{k-1}})}
{\partial \mathbf{q}} 
= \left(\frac{ \partial f}{ \partial \rho_{k-1}} \circ 
\mathcal{J}_{ \nu_{k-1}}\right)( \mathbf{q}, \mathbf{p})
\left(\rho_{k-1} \circ \mathcal{J}_{\nu_{k-1}}\right)( \mathbf{q},
\mathbf{p})(S^{k-1} - \mathbb{I}) \in L ^1
\end{equation}
since $(\rho_{k-1} \circ \mathcal{J}_{\nu_{k-1}})( \mathbf{q},
\mathbf{p}) \in L ^1$ and
\begin{equation}
\label{partial f by partial p}
\frac{\partial(f \circ \mathcal{J}_{\nu_{k-1}})}
{\partial \mathbf{p}}
= \left(\frac{ \partial f}{ \partial \rho_0} \circ \mathcal{J}_{
\nu_{k-1}} \right)(\mathbf{q}, \mathbf{p}) \in L ^{\infty}.
\end{equation}
Note that \eqref{partial f by partial q} implies that $f \circ
\mathcal{J}_{\nu_{k-1}} \in C ^{\infty}_{ \omega}(\ell ^{\infty}
\times \ell ^1)$ for any $f \in C^{\infty}(I^1_{-,0,k-1})$.

 Thus,  using the formula for the canonical
bracket on the weak symplectic Banach space  $(\ell^{\infty} \times
\ell^1, \omega)$ and the fact that the duality pairing 
$(L^{\infty})^\ast \times L ^{\infty} \rightarrow \mathbb{R}$
restricted to $L ^1 \times L ^{\infty}$ equals the trace of the
product, we get
\begin{align*}
&\{f \circ \mathcal{J}_{\nu_{k-1}}, g \circ \mathcal{J}_{
\nu_{k-1}} \}_ \omega (\mathbf{q}, \mathbf{p}) \\
& \quad = 
\left\langle \frac{\partial(f \circ
\mathcal{J}_{\nu_{k-1}})}{ \partial\mathbf{q}},
\frac{ \partial(g \circ \mathcal{J}_{\nu_{k-1}})}{
\partial \mathbf{p}} \right\rangle 
- \left\langle \frac{\partial(g \circ
\mathcal{J}_{\nu_{k-1}})}{ \partial\mathbf{q}},
\frac{ \partial(f \circ \mathcal{J}_{\nu_{k-1}})}{
\partial \mathbf{p}}  \right\rangle \\
& \quad = \operatorname{Tr} \left[\left(\rho_{k-1} \circ
\mathcal{J}_{\nu_{k-1}}\right)(\mathbf{q}, \mathbf{p}) 
\left(  (S^{k-1} -\mathbb{I})
\left(\frac{\partial g}{\partial\rho_0}
\circ \mathcal{J}_{\nu_{k-1}} \right)(\mathbf{q}, \mathbf{p})
\left(\frac{\partial f}{\partial \rho_{k-1}} \circ 
\mathcal{J}_{ \nu_{k-1}}\right)(\mathbf{q}, \mathbf{p}) 
\right. \right.\\ 
& \qquad \qquad \qquad \left.\left. 
-  (S^{k-1} - \mathbb{I})
\left(\frac{\partial f}{ \partial\rho_0} \circ
\mathcal{J}_{\nu_{k-1}}\right)(\mathbf{q}, \mathbf{p})
\left(\frac{\partial g}{\partial \rho_{k-1}} \circ 
\mathcal{J}_{\nu_{k-1}}\right)(\mathbf{q}, \mathbf{p})
\right) \right] \\
& \quad = \left(\{f, g\}_{0,k-1} \circ \mathcal{J}_{\nu_{k-1}} \right)
(\mathbf{q}, \mathbf{p})
\end{align*}
by \eqref{PB for k bidiagonal}.

Parts \textbf{(ii)} and
\textbf{(iii)} are proved by direct verifications.
\quad $\blacksquare$

\medskip

Let us define the map $\Phi^{\nu_{k-1}}(g):
GI_{+,0,k-1}^{\infty} \rightarrow \ell^{\infty}
\times \ell^1$ by
\begin{equation}
\label{definition of phi sub rho zero}
\Phi^ {\nu_{k-1}}(g): = \boldsymbol{\sigma}_g^{\nu_{k-1}}(\mathbf{0},
\mathbf{0}),
\end{equation}
or, in coordinates,
\begin{equation}
\label{phi sub rho zero in coordinates}
\Phi^ {\nu_{k-1}}(g_0, g_{k-1}) = \left(\log g_0, g_{k-1}g_0
^{-1}\nu_{k-1} - \tilde{s}^{k-1}(g_{k-1} g_0 ^{-1}
\nu_{k-1})
\right),
\end{equation}
which shows that  $\Phi^{\nu_{k-1}}$ is smooth
and injective.
\begin{proposition} 
\label{Flaschka map diagram}
The following diagram 

\unitlength=5mm
\begin{center}
\begin{picture}(17.5,6.2)
\put(-3.5,5){\makebox(0,0){$1$}}
\put(0.4,5){\makebox(0,0){$GL_0^{\infty,k-1}$}} 
\put(6.8,5){\makebox(0,0){$GI_{+,0,k-1}^{\infty}$}}
\put(15.5,5){\makebox(0,0){$GI_{+,0,k-1}^{\infty}/
                           GL_0^{\infty,k-1}$}}
\put(21.7,5){\makebox(0,0){$1$}}

\put(-4.5,0){\makebox(0,0){$0 \times \{\mathbf{p}\}$}}
\put(0.8,0){\makebox(0,0){$L_0^{\infty, k-1} \times 
\{\mathbf{p}\}$}} 
\put(7,0){\makebox(0,0){$\ell^{\infty} \times \ell^1$}} 
\put(16,0){\makebox(0,0){$(\delta^k)^{-1}(\Delta_{
\nu_{k-1}})$}}
\put(21.7,0){\makebox(0,0){$0$}}

\put(7,4){\vector(0,-1){3}}
\put(16,4){\vector(0,-1){3}}

\put(-3,5){\vector(2,0){1,7}}
\put(2,5){\vector(1,0){3}}
\put(8.8,5){\vector(1,0){3}}
\put(19,5){\vector(1,0){2.3}}

\put(-3,0){\vector(2,0){1,3}}
\put(3.4,0){\vector(1,0){2}}
\put(8.8,0){\vector(1,0){4.5}}
\put(18.5,0){\vector(1,0){2.6}}
\put(10.7, 5.5){\makebox(0,0){$\pi$}}
\put(11.5,.7){\makebox(0,0){$\mathcal{J}_{\nu_{k-1}}$}}
\put(17,2.4){\makebox(0,0){$\iota_{\nu_{k-1}}$}}
\put(6,2.4){\makebox(0,0){$\Phi^{\nu_{k-1}}$}}
\end{picture}
\end{center}
\bigskip
commutes. The first row is an exact sequence of Banach
Lie groups. The second row is also exact in the
following sense: the map
$\mathcal{J}_{\nu_{k-1}}$ is onto and its level sets are
all of the form $L_0^{\infty, k-1}\times \{ \mathbf{p}\} $,
where
$\mathbf{p}\in L ^1_0 $. In addition, 
\begin{equation}
\label{symplecticity}
(\Phi^{\nu_{k-1}}) ^\ast \omega = \pi^\ast
\omega_{\nu_{k-1}},
\end{equation}
where $\omega$ and $\omega_{\nu_{k-1}}$ are the weak
symplectic forms \eqref{symplectic form} and 
\eqref{k pull back in coordinates} on
$\ell^{\infty}\times \ell^1$ and
$GI_{+,0,k-1}^{\infty}/GL_0^{\infty, k-1}$
respectively. We also have 
\begin{equation}
\label{fiber preservation}
\Phi^{\nu_{k-1}}
\left( \pi^{-1}([g]) \right) =
\mathcal{J}_{\nu_{k-1}}^{-1}
\left(\iota_{\nu_{k-1}}([g])\right)
\end{equation}
for any $g \in GI_{+,0,k-1}^{\infty}$.
\end{proposition}

\noindent \textbf{Proof.} Commutativity is verified
using \eqref{group k coadjoint action}, \eqref{flaschka map},
and \eqref{definition of phi sub rho zero}.  The
identities \eqref{symplecticity} and \eqref{fiber
preservation} are obtained by direct verifications. \quad
$\blacksquare$
\medskip

\paragraph{Remarks.} (i) The analysis of the coadjoint orbit
$\mathcal{O}_\nu \cong GI ^{\infty}_{+,0,k-1}/
GL_0^{\infty,k-1}$ through the  generic element $ \nu\in
I_{-,0,k-1}^1$   carried out in this section shows that it
is diffeomorphic to $\Delta_{ \nu_0, \nu_{k-1}} \times
\Delta_{\nu_{k-1}}$. For an arbitrary
$(\nu'_0, \nu' _{k-1}) \in \Delta_{ \nu_0,
\nu_{k-1}} \times \Delta_{ \nu_{k-1}} $, the manifolds 
$\iota_{ \nu_{k-1}}^{-1}(
\{\nu'_0 \} \times  \Delta_{ \nu_{k-1}}))$ and
$\iota_{\nu_{k-1}}^{-1}(
\Delta_{\nu_0, \nu_{k-1}} \times \{\nu'_{k-1} \}))$
are Lagrangian submanifolds in the sense that their
tangent spaces are maximal isotropic. 

(ii) If $k=2$ we have $I ^1_{-,0,1} = L^1_{-,2}$ and
$GI^{\infty}_{+,0,1} = GL ^{\infty}_{+,2}$. If, in addition,
we consider the finite dimensional case, that is, instead of 
$L^1_{-,2}$ we work with the traceless $n \times n $ matrices
having  non-zero entries only on the main and the first lower
diagonals, then $\mathcal{J}_{ \nu_1}$ is a symplectic
diffeomorphism of $\mathbb{R}^{2(n-1)}$, endowed with the
canonical symplectic structure, with a single coadjoint
orbit of the upper bidiagonal group through a strictly lower
diagonal element all of whose entries are non-zero (see 
\cite{K} or, in tridiagonal symmetric formulation
\cite{A,S}).

(iii) If $k = 2$ and we consider the generic infinite
dimensional case, that is, $\nu_1$ has all entries
different from zero, then the map $\mathcal{J}_{\nu_1}$ does
not provide a morphism of weak symplectic manifolds between 
$\ell ^{\infty} \times \ell^1$ and a single 
coadjoint orbit of $ GL^{\infty}_{+,2} $. The relation
between these spaces is more complicated and is explained in
the diagram of Proposition \ref{Flaschka map diagram}. 
Each $GL^{\infty}_{+,2}$-coadjoint orbit through a generic
element $S^T\nu_1$ is only weakly symplectic and Poisson
injectively weakly immersed in $L^1_{-,2}$ but not equal to it.

(iv) If $k = 2$ and we consider the infinite
dimensional case with $\nu_1$ having also some vanishing
entries, the structure of the $GL^{\infty}_{+,2}$-coadjoint 
orbit through $S^T\nu_1$ reduces to the two previous cases
as we shall explain below. Let $i_0 $ be the first index for
which the entry $\nu_{1,i_0 i_0} = 0 $. Formula 
\eqref{group k coadjoint action} shows that the first $i_0
 \times i_0$ block of
$\mathcal{O}_{S^T \nu_1} $ is that of a
finite dimensional orbit of the upper bidiagonal group of
matrices of size $i_0\times i_0$  and that the coadjoint
action preserves this block. Let
$i_1$ be the next index for which $\nu_{1,i_1 i_1} = 0 $.
Again by
\eqref{group k coadjoint action} it follows that there is an
$i_1 \times i_1$ block of $\mathcal{O}_{S^T \nu_1}$ that is
preserved by the coadjoint action and that is equal to a
finite dimensional orbit of the upper bidiagonal group of
matrices of size $i_1\times i_1$. Continuing in this fashion
we arrive either at an infinite sequence of orbits of finite
dimensional upper bidiagonal groups (in the case that there
is an infinity of indices $i_s$ such that $\nu_{1,i_s i_s} =
0 $, $s \in \mathbb{N} \cup \{0\}$) or to a generic  infinite
dimensional orbit of $GL ^{\infty}_{+,2}$ (if there are only
finitely many indices $i_s$, $s = 0, 1, \ldots, r$, such
that $\nu_{1,i_s i_s} = 0 $). In the latter case, the last
infinite block is preserved by the coadjoint action and we
are in the generic case of an orbit of $GL ^{\infty}_{+,2}$
but on the space complementary to the $r +1$ finite dimensional
blocks of sizes $i_0 \times i_0 $, ..., $i_r \times i_r $.
Thus, decomposing the orbit as described, the problem of
classification of the general $GL ^{\infty}_{+,2}$-coadjoint
orbit is reduced to the finite dimensional case and to the
generic infinite dimensional case.

(v) One can restrict the Hamiltonians $I_{l}^{S,k}$
given by \eqref{i l s k} to $I ^1_{-,0,k-1}$
but these functions are not in involution because
the inclusion of $I ^1_{-,0,k-1}$ in $L ^1_{-, k}$ is not
Poisson. Indeed, as recalled in 
\S\ref{section: induced and coinduced}, the inclusion would
be Poisson if and only if the kernel of its dual map is an
ideal in $L ^{\infty}_{+,k}$ which is easily seen to be false
unless $k = 2 $, in which case we  have
\begin{equation}
\label{total linear momentum}
(I_{1}^{S,2} \circ  \mathcal{J}_{ \nu_{1}})( \mathbf{q},
\mathbf{p}) = \sum_{i=0} ^{\infty} p_i
\end{equation}
and 
\begin{equation}
\label{k toda hamiltonians}
H_2( \mathbf{q}, \mathbf{p}): = 
(I_{2}^{S,2} \circ  \mathcal{J}_{ \nu_{1}})( \mathbf{q},
\mathbf{p}) = \frac{1}{2}\sum_{i=0} ^{\infty} p_i ^2 + 
\sum_{i=0} ^{\infty} \nu_{1,ii} ^2 e^{2(q_{i+1} -
q_i)}.
\end{equation}
The function $H_2$ is, up to a renormalization of constants,
the Hamiltonian of the semi-infinite Toda
lattice. The first  integral $I_{1}^{S,2} \circ 
\mathcal{J}_{ \nu_{1}} $ is the total momentum of the system
which generates the translation action given by the subgroup
$\mathbb{R}_+ \mathbb{I}$. All integrals
$I_{l}^{S,2} \circ \mathcal{J}_{\nu_{1}}$, $l \in
\mathbb{N}$,  give the full Toda lattice hierarchy on 
$\ell^{\infty} \times \ell^1$; see \S\ref{toda example}.
\medskip

These considerations justify the name of
{\bfi Flaschka map\/} for the momentum map
$\mathcal{J}_{\nu_1}: \ell ^{\infty} \times \ell ^1
\rightarrow  I ^1_{-,0,1} = L ^1_{-,2}$. In the next section
we will present a momentum map from the weak symplectic
manifold $\left(\ell ^{\infty} \right)^{k-1} \times
\left(\ell ^1 \right)^{k-1}$, endowed with a weak magnetic
symplectic structure, to the Banach Lie-Poisson
space $L ^1_{-,k}$. This momentum map can be considered, as
we shall see, as a natural generalization of the Flaschka
map to the system of integrals in involution 
\eqref{i l s k} for $k \geq 2$.

\section{The Flaschka map for $\left(\ell ^{\infty} \right)^{k-1} 
\times \left(\ell ^1 \right)^{k-1}$ }
\label{sect: generalized Flaschka}

In this section we construct a $GL ^{\infty}_{+,k}$-equivariant
momentum map $\mathbf{J}_k:\left(\ell ^{\infty} \right)^{k-1} \times
\left(\ell ^1 \right)^{k-1} \rightarrow L ^1_{-,k}$ (see 
\eqref{momentum on concrete k quotient}) which can be interpreted as
a generalization of the Flaschka map \eqref{flaschka map} defined for
the bidiagonal case. We also construct a weak symplectic  form
$\Omega_k$  on
$\left(\ell^{\infty} \right)^{k-1} \times \left(\ell^1 \right)^{k-1}$
(see \eqref{concrete k symplectic form}) which has a
non-canonical term responsible for the
interaction of the Toda system with some kind of an  external
``field". We shall illustrate
the hierarchy of dynamical systems obtained in this way  by
studying the special case
$k=3$ in detail (see
\eqref{second degree 3 hamiltonian}). The simpler case $k=2$
does not add anything new since one recovers by the
symplectic  induction method the original semi-infinite Toda system
studied in the previous section.

 We shall apply the induction 
method discussed in \S\ref{section: symplectic induction}
to the weak symplectic manifold $(P, \omega) =  
(\ell^{\infty} \times \ell^1,\omega)$ with $\omega$  
given by \eqref{symplectic form}, the Banach Lie 
group $G : = (GL^{\infty}_{+,k}, \circ_k)$ defined 
in \eqref{gl plus k}, and the Banach Lie subgroup 
$H: = GI^{\infty}_{+,0,k-1}$ consisting of invertible 
bidiagonal elements of the form
\eqref{bidiagonal k group element}. As will be seen, 
the  abstract constructions presented in 
\S\ref{section: symplectic induction}
become completely explicit in this case.

We begin by listing the objects involved  in this
construction. The Banach Lie algebra is $\mathfrak{g}: = L
^{\infty}_{+,k} = \oplus_{i=0}^{k-1} L_i ^{\infty}$, the
subalgebra is $\mathfrak{h}: = I ^{\infty}_{+,0,k-1} = L
^{\infty}_0 \oplus L ^{\infty}_{k-1}$, and its closed split
complement is $\mathfrak{h}^\perp := \oplus_{i=1}^{k-2}L_i
^{\infty} =: (I ^{\infty}_{+,0,k-1})^\perp$.  At
the level of the preduals we have $\mathfrak{g}_\ast =
L^1_{-,k} = \oplus_{i=-k+1}^\perp L_i$, $\mathfrak{h}_\ast = I
^1_{-,0,k-1} = L ^1_0 \oplus L ^1_{-k+1}$, and its closed
split complement $\mathfrak{h}_\ast ^\perp = \oplus_{i=-k+2}^{-1}
L ^1_i = : (I ^1_{-,0,k-1})^\perp$. We have hence the
Banach space direct sums 
\begin{equation}
\label{first k concrete dec}
L ^{\infty}_{+,k} =
I^{\infty}_{+,0,k-1} \oplus (I ^{\infty}_{+,0,k-1})^\perp 
\end{equation} 
and 
\begin{equation}
\label{second k concrete dec}
L ^1_{-,k} = I^1_{-,0,k-1} \oplus (I^1_{-,0,k-1})^\perp.
\end{equation}
Thus any $\rho \in L ^1_{-,k}$ uniquely decomposes as $\rho=
\gamma + \gamma^\perp$, where $\gamma = \rho_0 + (S^T)^{k-1}
\rho_{k-1} \in I^1_{-,0,k-1}$ and $\gamma^\perp = S^T
\rho_1 +  \dots (S^T)^{k-2} \rho_{k-2} \in
(I^1_{-,0,k-1})^\perp$.
Let us show that the splitting \eqref{first k concrete dec}
is invariant relative  to the restriction of the adjoint
action $\operatorname{Ad}^{+,k}$ of the Banach Lie group
$GL^{\infty}_{+,k}$ to the Lie subgroup
$GI^{\infty}_{+,0,k-1}$. Clearly the factor
$I^{\infty}_{+,0,k-1}$ is preserved because it is the Lie
algebra of $GI^{\infty}_{+,0,k-1}$. To see that the second
factor $(I^{\infty}_{+,0,k-1})^\perp$ is also preserved, using
\eqref{inverse in k bidiagonal}, it suffices to show that for
any
$h = h_0 + h_{k-1} S^{k-1}
\in GI^{\infty}_{+,0,k-1}$ and any $x_1 S + \dots +
x_{k-2}S^{k-2}
\in (I ^{\infty}_{+,0,k-1})^\perp$ we have 
\begin{align}
\label{adjoint k action on complement}
&(\operatorname{Ad}^{+,k})_h(x_1 S + \dots
x_{k-2}S^{k-2})  \nonumber \\  
& \qquad =
\left(h_0 + h_{k-1}S^{k-1}\right) \circ _k \left(x_1 S +
\dots + x_{k-2}S^{k-2}\right) \circ_k \left(h_0 ^{-1} -
h_0^{-1} h_{k-1} s^{k-1}(h_0 ^{-1}) S^{k-1} \right)
\nonumber\\ 
& \qquad = h_0 s(h_0 ^{-1}) x_1 S + \dots + h_0
s^{k-2}(h_0 ^{-1}) x_{k-2} S^{k-2}
\end{align}
which is a straightforward verification.

Next we show that the splitting 
\eqref{second k concrete dec} is invariant relative  to the
restriction of the coadjoint action $(
\operatorname{Ad}^{+,k}) ^\ast$ of 
$GL^{\infty}_{+,k}$ to the Lie subgroup
$GI^{\infty}_{+,0,k-1}$. First, by \eqref{group k coadjoint
action} the $GI^{\infty}_{+,0,k-1}$ coadjoint action
preserves the predual $I^1_{-,0,k-1}$. Second, to show that
the second factor $(I^1_{-,0,k-1})^\perp$ is also preserved, one
verifies directly, using \eqref{inverse in k bidiagonal}, 
that for any $h = h_0 + h_{k-1} S^{k-1}
\in GI^{\infty}_{+,0,k-1}$ and $S^T \rho_1 + \dots +
(S^T)^{k-2} \rho_{k-2} \in (I^1_{-,0,k-1})^\perp$ we have 
\begin{align}
\label{coadjoint k action on complement}
&(\operatorname{Ad}^{+,k})^\ast_{h^{-1}}(S^T \rho_1 + \dots +
(S^T)^{k-2} \rho_{k-2})  \nonumber \\  
& \qquad = S^T s(h_0) h_0 ^{-1} \rho_1 + \dots + (S^T)^{k-2}
s^{k-2}(h_0) h_0 ^{-1}\rho_{k-2}.
\end{align}

According to the general theory we shall take the weak
symplectic manifolds
$GL ^{\infty}_{+,k} \times L ^1 _{-,k}$ and $\ell^{\infty}
\times \ell^1$, the canonical action
$\boldsymbol{\sigma}^{\nu_{k-1}}: GI^{\infty}_{+,0,k-1} \times
(\ell^{\infty}\times \ell^1)
\rightarrow \ell^{\infty}\times \ell^1$ defined in
\eqref{sigma action}, and its equivariant momentum map 
$\mathcal{J}_{\nu_{k-1}}: 
\ell^{\infty}\times \ell^1 \rightarrow I^1_{-,0,k-1}$ given
by \eqref{flaschka map} (see
Proposition \ref{Flaschka map theorem}). We fix in all
considerations below an element $\nu_{k-1} \in L ^1_0$. By
\eqref{h action}, the Banach Lie group
$GI^{\infty}_{+,0,k-1}$ acts on the product $(\ell^{\infty}
\times \ell^1) \times GL^{\infty}_{+,k} \times L^1 _{-,k}$
by
\[
h \cdot (( \mathbf{q}, \mathbf{p}), g, \rho): =
\left(\boldsymbol{\sigma}^{\nu_{k-1}} (\mathbf{q}, \mathbf{p}), g \circ_k
h^{-1}, (\operatorname{Ad}^{+,k}) ^\ast_{h ^{-1}} \rho
\right),
\]
where $h \in GI ^{\infty}_{+,0,k-1}$, $g \in
GL^{\infty}_{+,k}$, $( \mathbf{q}, \mathbf{p}) \in
\ell^{\infty} \times \ell ^1 $, and $\rho \in L ^1_{-,k}$.
This action admits the equivariant momentum map
\eqref{big mom map}, which in this case becomes
\begin{align*}
&((\mathbf{q}, \mathbf{p}), g, \gamma+ \gamma^\perp) \in
(\ell^{\infty} \times \ell ^1) \times GL^{\infty}_{+,k}
\times \left(I^1_{-,0,k-1} \oplus (I^1_{-,0,k-1})^\perp
\right) \\
& \qquad \qquad  \longmapsto
\mathcal{J}_{\nu_{k-1}}(\mathbf{q}, \mathbf{p}) - \gamma \in
I^1_{-,0,k-1}.
\end{align*}
The zero level set of this momentum map is a smooth manifold,
$GI^{\infty}_{+,0,k-1}$-equivariantly diffeomorphic to 
$GL^{\infty}_{+,k} \times (\ell^{\infty} \times \ell^1)
\times (I ^1_{-,0,k-1})^\perp$, the action on the target being
\[
h \cdot \left(g, \mathbf{q}, \mathbf{p}, \gamma^\perp \right) :
= \left(g \circ_k h ^{-1},
\boldsymbol{\sigma}^{\nu_{k-1}}_h(\mathbf{q},
\mathbf{p}),
\left(\operatorname{Ad}^{+,k}\right)^\ast_{h^{-1}} \gamma^\perp
\right).
\] 
The symplectically induced space is hence the fiber
bundle  
\[
GL^{\infty}_{+,k} \times_{GI^{\infty}_{+,0,k-1}} \left(
\ell^{\infty} \times \ell ^1 \times (I ^1_{-,0,k-1})^\perp
\right) \rightarrow 
GL^{\infty}_{+,k}/GI^{\infty}_{+,0,k-1}
\] 
associated to the principal bundle $GL^{\infty}_{+,k}
\rightarrow GL^{\infty}_{+,k}/GI^{\infty}_{+,0,k-1}$. 

We
begin by explicitly determining the base of this bundle. If
$g = g_0 + \dots + g_{k-1} S^{k-1} \in GL ^{\infty} _{+,k}$
and $h = h_0 + h_{k-1}S^{k-1} \in GI ^{\infty}_{+,0,k-1}$
then
\begin{align*}
g\circ_k h ^{-1} &= (g_0 + \dots + g_{k-1} S^{k-1}) \circ _k
(h_0 ^{-1} - h_0^{-1} h_{k-1} s^{k-1}(h_0 ^{-1}) S^{k-1} ) \\
& = g_0 h_0^{-1} + g_1 s( h_0 ^{-1}) S + \dots +
g_{k-2}s^{k-2}(h_0^{-1}) S^{k-2}  \\
& \qquad + \left(g_{k-1}s^{k-1}(h_0^{-1}) - g_0 h_0
^{-1}h_{k-1} s^{k-1}(h_0 ^{-1}) \right) S^{k-1}.
\end{align*}
Therefore, the smooth map $GL^{\infty} _{+,k} \rightarrow
\left(\ell^{\infty} \right)^{k-2}$ given by
\begin{align*}
GL^{\infty} _{+,k} \ni
g_0 + \dots + g_{k-1} S^{k-1} &\mapsto (g_0 + \dots + g_{k-1}
S^{k-1}) \circ _k (g_0 ^{-1} - g_0^{-1} g_{k-1}
s^{k-1}(h_0^{-1}) S^{k-1}) \\
& = \mathbb{I} + g_1 s(g_0 ^{-1}) S + \dots +
g_{k-2}s^{k-2}(g_0^{-1}) S^{k-2} \\
& \mapsto \left(g_1 s(g_0 ^{-1}), \dots ,
g_{k-2}s^{k-2}(g_0^{-1}) \right)
\in \left(\ell^{\infty} \right)^{k-2}
\end{align*}
factors through the $GI^{\infty}_{+,0,k-1}$-action thus
inducing a smooth map
$GL^{\infty}_{+,k}/GI^{\infty}_{+,0,k-1} \rightarrow
\left(\ell^{\infty} \right)^{k-2}$. Its inverse is the
smooth map 
\[
\left(\mathbf{q}_1, \dots ,
\mathbf{q}_{k-2} \right) \in
\left(\ell^{\infty} \right)^{k-2} \mapsto
[ \mathbb{I} +
\mathbf{q}_1 S + \dots + \mathbf{q}_{k-2}S^{k-2}] \in 
GL^{\infty}_{+,k}/GI^{\infty}_{+,0,k-1}
\]
which proves that $GL^{\infty}_{+,k}/GI^{\infty}_{+,0,k-1}$
is diffeomorphic to $\left(\ell^{\infty} \right)^{k-2}$. 

Next, we shall prove that the smooth map 
\[
\Phi: (\ell^{\infty} \times \ell ^1) \times 
\left(\ell^{\infty} \right)^{k-2} \times
\left(\ell^1\right)^{k-2} \rightarrow GL^{\infty}_{+,k}
\times_{GI^{\infty}_{+,0,k-1}} \left(
\ell^{\infty} \times \ell ^1 \times (I ^1_{-,0,k-1})^\perp
\right)
\]
given by
\begin{align*}
&\Phi\left((\mathbf{q},\mathbf{p}), \mathbf{q}_1,\dots,
\mathbf{q}_{k-2},  \mathbf{p}_1,\dots,
\mathbf{p}_{k-2} \right) \\  
&\qquad := \left[\left(\mathbb{I} + \mathbf{q}_1 S + \dots +
\mathbf{q}_{k-2}S^{k-2}, (\mathbf{q}, \mathbf{p}),
S^T \mathbf{p}_1 + \dots + (S^T)^{k-2} \mathbf{p}_{k-2}
\right)\right] 
\end{align*}
is a diffeomorphism thereby trivializing the associated
bundle,  which is the reduced space.
Indeed, this map has a smooth inverse given by
\begin{align*}
&\Phi^{-1}\left( \left[\left(g_0 + \dots + g_{k-1}S^{k-1},
(\mathbf{q}, \mathbf{p}), \gamma^\perp \right)\right]\right) \\
& \qquad  = 
\left(\boldsymbol{\sigma}^{\nu_{k-1}}_{g_0 +
g_{k-1}S^{k-1}} (\mathbf{q},\mathbf{p}), g_1 s(g_0 ^{-1}),
\dots , g_{k-2}s^{k-2}(g_0^{-1}), 
\left(\operatorname{Ad}^{+,k} \right)^\ast_{(g_0 +
g_{k-1}S^{k-1})^{-1}}\gamma^\perp \right),
\end{align*}
where, in the third component of the right hand side we have
identified $(I ^1_{-,0,k-1})^\perp$ with $(\ell^1)^{k-2}$
through the isomorphisms $L ^1_k \cong \ell ^1$.
 
The $GL^{\infty}_{+,k}$-action on the reduced manifold 
$GL^{\infty}_{+,k} \times_{GI^{\infty}_{+,0,k-1}} \left(
\ell^{\infty} \times \ell ^1 \times (I ^1_{-,0,k-1})^\perp
\right)$ is given by $g' \cdot [g, (\mathbf{q}, \mathbf{p}),
\gamma^\perp] = [g'\circ _k g, (\mathbf{q}, \mathbf{p}),
\gamma^\perp]$ for any $g', g \in GL^{\infty}_{+,k}$,
$(\mathbf{q}, \mathbf{p}) \in \ell^{\infty} \times \ell^1$,
and $\gamma^\perp \in (I^1_{-,0,k-1})^\perp$. 
 Via the globally trivializing
diffeomorphism
$\Phi$, the induced $GL^{\infty}_{+,k}$-action on
$(\ell^{\infty} \times \ell ^1) \times 
\left(\ell^{\infty} \right)^{k-2} \times
\left(\ell^1\right)^{k-2} $ has the expression
\begin{align}
&(g_0 + \dots + g_{k-1}S^{k-1}) \cdot
\left((\mathbf{q}, \mathbf{p}), \mathbf{q}_1,\dots,
\mathbf{q}_{k-2},  \mathbf{p}_1,\dots,
\mathbf{p}_{k-2} \right) \nonumber \\ 
& \qquad = \Phi^{-1} \left((g_0 + \dots +
g_{k-1}S^{k-1}) \cdot \Phi \left((\mathbf{q}, \mathbf{p}),
\mathbf{q}_1,\dots,
\mathbf{q}_{k-2},  \mathbf{p}_1,\dots,
\mathbf{p}_{k-2} \right) \right)\nonumber \\
& \qquad = \Phi^{-1}\left((g_0 + \dots + g_{k-1}S^{k-1})
\cdot \left[\left( \mathbb{I} +
\mathbf{q}_1 S + \dots + \mathbf{q}_{k-2}S^{k-2},
(\mathbf{q}, \mathbf{p}), \right. \right. \right. \nonumber
\\ & \qquad \qquad \qquad \qquad \qquad \qquad \qquad \qquad
\qquad \qquad \qquad 
\left. \left. \left. S^T \mathbf{p}_1 +
\dots + (S^T)^{k-2} \mathbf{p}_{k-2} \right)\right] \right)
\nonumber \\ 
& \qquad = \Phi^{-1}\left(\left[\left((g_0 + \dots +
g_{k-1}S^{k-1}) \circ _k (\mathbb{I} +
\mathbf{q}_1 S + \dots + \mathbf{q}_{k-2}S^{k-2}),
(\mathbf{q}, \mathbf{p}),  \right. \right. \right. \nonumber
\\ & \qquad \qquad \qquad \qquad \qquad \qquad \qquad \qquad
\qquad \qquad \qquad 
\left. \left. \left.
S^T \mathbf{p}_1 + \dots + (S^T)^{k-2} \mathbf{p}_{k-2}
\right)\right] \right)
\nonumber \\ 
& \qquad = \Phi^{-1}\left(\left[ \left( g_0 +
\sum_{l=1}^{k-2}  \left(
\sum_{i=0}^l g_{l-i} s^{l-i}( \mathbf{q}_i) \right)S ^l +
\left(\sum_{i=0}^{k-2} g_{k-1-i} s^{k-1-i}( \mathbf{q}_i)
\right)S^{k-1}, (\mathbf{q}, \mathbf{p}), 
 \right. \right. \right. \nonumber \\
& \qquad \qquad \qquad \qquad \qquad \qquad \qquad \qquad
\qquad \qquad \quad 
\left. \left. \left. \phantom{\sum_{i=0}^{k-2}}
S^T \mathbf{p}_1 +
\dots + (S^T)^{k-2} \mathbf{p}_{k-2} \right) \right]
\right) \nonumber \\
& \qquad = \left(\boldsymbol{\sigma}^{\nu_{k-1}}_{g_0 +
\left(\sum_{i=0}^{k-2} g_{k-1-i} s^{k-1-i}( \mathbf{q}_i)
\right)S^{k-1}} ( \mathbf{q}, \mathbf{p}),
\phantom{\sum_{i=0}^1}  \right. \nonumber \\
&\qquad \qquad \qquad \left. s(g_0 ^{-1}) \sum_{i=0}^1
g_{1-i} s^{1-i}(\mathbf{q}_i) ,
\dots , s(g_0 ^{-1}) \sum_{i=0}^{k-2} g_{k-2-i} s^{k-2-i}(
\mathbf{q}_i), \right.  \nonumber \\
& \qquad \qquad \; \;   \left. \phantom{\sum_{i=0}^1}
s(g_0) g_0^{-1}, \dots, s^{k-2}(g_0) g_0^{-1} \right),
\nonumber
\end{align}
where the equality in the last $k-2$ components follows from
\eqref{coadjoint k action on complement}. Let us summarize the
considerations above. Using 
\eqref{sigma action} and denoting 
$\left((\mathbf{q}', \mathbf{p}'), \mathbf{q}_1 ', \dots,
\mathbf{q}_{k-2}', \mathbf{p}_1 ', \dots,
\mathbf{p}_{k-2}'\right) : = (g_0 + \dots + g_{k-1}S^{k-1}) 
\cdot \left((\mathbf{q}, \mathbf{p}), \mathbf{q}_1,\dots,
\mathbf{q}_{k-2},  \mathbf{p}_1,\dots,
\mathbf{p}_{k-2} \right)$, we conclude that the
$GL^{\infty}_{+,k}$-action on the reduced manifold 
$(\ell^{\infty} \times \ell ^1) \times 
\left(\ell^{\infty} \right)^{k-2} \times
\left(\ell^1\right)^{k-2}$ is given by
\begin{align}
\label{k action q}
\mathbf{q}' &=  \mathbf{q} + \log g_0\\
\label{k action p}
\mathbf{p}' &= \mathbf{p}+
\left(\sum_{i=0}^{k-2} g_{k-1-i} s^{k-1-i}( \mathbf{q}_i)
\right) g_0 ^{-1}
\nu_{k-1} e^{s^{k-1}(\mathbf{q}) - \mathbf{q}} \nonumber \\
& \qquad  -
\tilde{s}^{k-1}
\left(\left(\sum_{i=0}^{k-2} g_{k-1-i} s^{k-1-i}( \mathbf{q}_i)
\right) g_0 ^{-1} \nu_{k-1} e^{s^{k-1}( \mathbf{q}) -
\mathbf{q}}\right) \\
\label{k action q prime}
\mathbf{q}_l' &= s(g_0 ^{-1}) \sum_{i=0}^{l} g_{l-i}
s^{l-i}( \mathbf{q}_i) \\
\label{k action p prime}
\mathbf{p}_l' &= s^{l}(g_0) g_0 ^{-1}\mathbf{p}_{l}, \qquad
l = 1, \dots, k-2.
\end{align}
All geometric objects described above satisfy the assumptions of
Propositions \ref{induction theorem} and 
\ref{induction theorem second} and thus one has the weak symplectic 
form $\Omega_k$ and the momentum map 
$\mathbf{J}_k: (\ell^{\infty} \times \ell^1) \times
\left(\ell ^{\infty}
\right)^{k-2} \times \left(\ell ^1 \right)^{k-2} \rightarrow
L ^1_{-,k}$ given by \eqref{quotient symplectic form} and
\eqref{induced momentum}, respectively. By
\eqref{group coordinate coadjoint minus k action}, $\mathbf{J}_k$
takes the form
\begin{align}
\label{momentum on concrete k quotient}
&\mathbf{J}_k\left((\mathbf{q}, \mathbf{p}), \mathbf{q}_1,
\dots
\mathbf{q}_{k-2}, \mathbf{p}_1, \dots
\mathbf{p}_{k-2} \right)   \nonumber\\ 
& \quad = \left(\operatorname{Ad}^{+,k}
\right)^\ast_{(\mathbb{I} +
\mathbf{q}_1S + \dots + \mathbf{q}_{k-2} S^{k-2})^{-1} }
\left(\mathcal{J}_{\nu_{k-1}}( \mathbf{q}, \mathbf{p}) +
S^T \mathbf{p}_1 + \dots + (S^T)^{k-2} \mathbf{p}_{k-2}
\right)  \nonumber\\
& \quad = \left(\operatorname{Ad}^{+,k}
\right)^\ast_{(\mathbb{I} +
\mathbf{q}_1S + \dots + \mathbf{q}_{k-2} S^{k-2})^{-1} }
\left(\mathbf{p} + S^T \mathbf{p}_1 + \dots + (S^T)^{k-2}
\mathbf{p}_{k-2}
\phantom{\nu_{k-1} e^{s^{k-1} ( \mathbf{q}) - \mathbf{q}} }
\right.  \nonumber \\ 
& \qquad \qquad \left.  + (S^T)^{k-1}
\nu_{k-1} e ^{s^{k-1} ( \mathbf{q}) - \mathbf{q}}
\right),
\end{align}
where the inverse $(\mathbb{I} +
\mathbf{q}_1S + \dots + \mathbf{q}_{k-2} S^{k-2})^{-1}$  is
given by \eqref{inverse k explicit}. We shall call
$\mathbf{J}_k $ the {\bfi generalized Flaschka map}. 

\medskip

In order to obtain the explicit expression of the weak symplectic
form $\Omega_k $ (see \eqref{concrete k symplectic form}) on the
induced symplectic  manifold
$(\ell ^{\infty} \times \ell ^1) \times \left(\ell ^{\infty}
\right)^{k-2} \times \left(\ell ^1 \right)^{k-2}$, let us
notice that the symplectic form $\omega+ \omega_L$ on $(\ell
^{\infty}\times \ell ^1) \times GL ^{\infty}_{+,k} \times L
^1_{-,k}$ is given by 
\begin{equation}
\label{symplectic form k concrete}
\omega + \omega_L = -\mathbf{d} \left(
\operatorname{Tr}(\mathbf{p}\mathbf{d}\mathbf{q})
+ \operatorname{Tr}(\rho g ^{-1}\circ_k\mathbf{d}g) \right),
\end{equation}
where $g ^{-1}\circ _k\mathbf{d}g$ is the left Maurer-Cartan
form on the Banach Lie group $GL ^{\infty}_{+,k}$. One has
the following decomposition
\begin{equation}
\label{theta formula}
\theta: = \operatorname{Tr}(\rho g ^{-1}\circ_k\mathbf{d}g)
= \operatorname{Tr}\left(\sum_{l=0}^{k-1} \rho_l
\theta_l \right)
\end{equation}
for $\rho= \rho_0 + S^T \rho_1 + \dots + (S^T)^{k-1}
\rho_{k-1} \in L ^1_{-,k}$ with 
\[
\theta_l = \sum_{i=0}^l h_i(g) s^i(\mathbf{d}g_{l-i}), \qquad
l=0,1, \dots, k-1.
\]
The diagonal operators $h_i$ are the components of $g ^{-1}
= h_0 +h_1 S + \dots + h_{k-1}S^{k-1}$ given by
\eqref{inverse k explicit}. Let $\tilde{\theta}$ be the pull
back of $\theta$ to the zero level set of the momentum map 
\eqref{big mom map}. Next, we pull back the form $\tilde{\theta}$ to 
$(\ell^{\infty} \times \ell ^1) \times \left(\ell^{\infty}
\right)^{k-2} \times \left(\ell ^1 \right)^{k-2}$ by
the global section $\Sigma : (\ell^{\infty} \times \ell^1)
\times \left(\ell^{\infty} \right)^{k-2} \times \left(\ell^1
\right)^{k-2} \rightarrow GL^{\infty}_{+,k} \times \left(
\ell^{\infty} \times \ell^1\right) \times (I^1_{-,0,k-1})^\perp$
defined by 
\begin{align*}
&\Sigma((\mathbf{q}, \mathbf{p}), \mathbf{q}_1,
\dots \mathbf{q}_{k-2}, \mathbf{p}_1, \dots,
\mathbf{p}_{k-2}) \\
& \qquad := \left(\mathbb{I} + \mathbf{q}_1 S, +
\dots + \mathbf{q}_{k-2}S^{k-2}, ( \mathbf{q}, \mathbf{p}),
S^T \mathbf{p}_1 + \dots + (S^T)^{k-2} \mathbf{p}_{k-2}
\right).
\end{align*}
Therefore, we get
\begin{align}
\label{concrete k one form}
\Sigma^\ast \tilde{\theta} := &\operatorname{Tr}(
\mathbf{p}\mathbf{d}\mathbf{q}) +
\operatorname{Tr} \left[\left(\mathcal{J}_{ \nu_{k-1}} (
\mathbf{q},\mathbf{p}) \right)_0 \theta_0 \right] +
\operatorname{Tr} \left(\left(\mathcal{J}_{ \nu_{k-1}} (
\mathbf{q},\mathbf{p}) \right)_{k-1} \theta_{k-1} \right)
+ \operatorname{Tr} \left(\sum_{l=1}^{k-2} \mathbf{p}_l
\theta_l \right) \nonumber\\
= &\operatorname{Tr}(\mathbf{p}\mathbf{d}\mathbf{q}) +
\operatorname{Tr} \left(\sum_{l=1}^{k-2} \mathbf{p}_l
\sum_{i=0}^{l-1}h_i(\mathbf{q}_1,\dots , \mathbf{q}_{i}) s^i
(\mathbf{d}\mathbf{q}_{l-i}) \right) \nonumber\\
& \qquad  + \operatorname{Tr}\left(
\nu_{k-1} e^{ s^{k-1}( \mathbf{q})- \mathbf{q}}
\sum_{i=1}^{k-2} h_i(\mathbf{q}_1,\dots , \mathbf{q}_{i})
s^i( \mathbf{d} \mathbf{q}_{k-1-i}) \right),
\end{align}
since $\theta_0 = 0 $, where $h_i( \mathbf{q}_1, \dots,
\mathbf{q}_{i})$ is given by \eqref{inverse k explicit}
with $g_0 = (1,1, \dots )$, $g_1 = \mathbf{q}_1, \dots$, $g_{k-2} =
\mathbf{q}_{k-2}$, $g_{k-1} = ( 0, 0, \dots )$. Since
$\operatorname{Tr}\delta  = \operatorname{Tr} \tilde{s}^j( \delta)$
for any $\delta\in L ^1_0$ and $j \in \mathbb{N}$, by \eqref{useful
s identities} the last summand in  \eqref{concrete k one form}
becomes 
\begin{align*}
&\sum_{i=1}^{k-2}\operatorname{Tr}\left[ \tilde{s}^i\left(\nu_{k-1}
e^{ s^{k-1} ( \mathbf{q})- \mathbf{q}} h_i(\mathbf{q}_1,\dots ,
\mathbf{q}_{i}) \right) \left(\mathbb{I} - \sum_{r=0}^{i-1}p_r
\right)\mathbf{d} \mathbf{q}_{k-1-i}\right] \\
& \qquad = \sum_{i=1}^{k-2} \operatorname{Tr} \left[
\tilde{s}^i\left(\nu_{k-1} e^{ s^{k-1} ( \mathbf{q})- \mathbf{q}}
h_i(\mathbf{q}_1,\dots ,
\mathbf{q}_{i})\right)\mathbf{d} \mathbf{q}_{k-1-i}\right]
\end{align*}
because 
\[
\tilde{s}^j( \delta)\sum_{r=0}^{j-1}p_r = 0 \qquad \text{for all}
\qquad \delta\in L ^1_0 \qquad \text{and} \qquad j \in \mathbb{N}.
\]
Similarly, the second  summand in \eqref{concrete k one form} equals
\[
\sum_{l=1}^{k-2} \sum _{i=0}^{l-1} \operatorname{Tr}
\left[\tilde{s}^i \big( \mathbf{p}_l
 h_i(\mathbf{q}_1,\dots , \mathbf{q}_{i}) \big) \mathbf{d}
\mathbf{q}_{l-i}
\right]
\]
so that \eqref{concrete k one form} becomes
\begin{align}
\label{concrete k one form shifted}
\Sigma^\ast \tilde{\theta} &= \operatorname{Tr}(
\mathbf{p}\mathbf{d}\mathbf{q}) + \sum_{l=1}^{k-2}
\operatorname{Tr} \left(\sum_{i=0}^{l-1} \tilde{s}^i 
\left(\mathbf{p}_l h_i(\mathbf{q}_1,\dots , \mathbf{q}_{i}) \right)
\mathbf{d} \mathbf{q}_{l-i} \right. \nonumber \\
& \qquad  \left. \phantom{\sum_{i=0}^{l-1}} +
\tilde{s}^l\left(\nu_{k-1} e^{ s^{k-1} ( \mathbf{q})- \mathbf{q}}
h_l(\mathbf{q}_1,\dots , \mathbf{q}_{l})\right)\mathbf{d}
\mathbf{q}_{k-1-l} \right) \nonumber\\
& = \operatorname{Tr}(\mathbf{p}\mathbf{d}\mathbf{q}) +
\sum_{l=1}^{k-2}\left[ \operatorname{Tr}\left(\sum_{i=0}^{k-2-l}
\tilde{s}^i  \left(\mathbf{p}_l h_i(\mathbf{q}_1,\dots ,
\mathbf{q}_{i}) \right)  \right. \right. \nonumber \\
& \qquad \left. \left.  \phantom{\sum_{i=0}^{l-1}} +
\tilde{s}^l\left(\nu_{k-1} e^{ s^{k-1} (
\mathbf{q})- \mathbf{q}} h_l(\mathbf{q}_1,\dots , \mathbf{q}_{l})
\right) \right) \mathbf{d} \mathbf{q}_l \right].
\end{align}
 
Then the reduced symplectic form is 
\begin{equation}
\label{concrete k symplectic form}
\Omega_k = - \mathbf{d}
\Sigma^\ast \tilde{\theta}.
\end{equation}
Indeed, a straightforward
verification shows that $- \mathbf{d} \Sigma^\ast
\tilde{\theta}$ satisfies the condition characterizing the
reduced symplectic form, so it must be equal to it. Note
that the one-form $\Sigma^\ast \tilde{\theta} $
depends on the chosen section $\Sigma$, but that if
$\tilde{\Sigma}$ is any other global section, then
$\mathbf{d} \Sigma^\ast \tilde{\theta} = \mathbf{d}
\tilde{\Sigma}^\ast\tilde{\theta} = \Omega_k$. 
In particular, the reduced symplectic form $\Omega_k$ is in
this case exact. Note also that the symplectic form $\Omega_k$ is
canonical only if $k=2$ and magnetic only if $k=3$, a case that  we
shall analyze in detail below. In general, if $k>3$, the weak
symplectic  form $\Omega_k$ is neither canonical nor magnetic due
to the presence of the $\mathbf{p}_j$-dependent coefficients of
$\mathbf{d}\mathbf{q}_l$ in the first sum of the second term. 

To deal with the Hamiltonian systems defined by the
function $I _l^{S,k}$ we need to regard the momentum map
$\mathbf{J}_k$ as having values in $L ^1_{S,k}$. This is achieved
by defining the equivariant momentum map $\mathbf{J}_k^S: =
\Phi_{S,-,k} \circ \mathbf{J}_k ^S: (\ell^{\infty} \times \ell^1)
\times \left(\ell ^{\infty} \right)^{k-2} \times \left(\ell ^1
\right)^{k-2} \rightarrow L ^1_{S,k}$, where $\Phi_{S,-,k}: 
\left(L^1_{-,k}, \{ \cdot , \cdot \}_{-,k} \right) \rightarrow
\left( L ^1_{S,k}, \{ \cdot , \cdot \}_{S,k} \right)$ is the
isomorphism of Banach Lie-Poisson spaces introduced at the end of
\S \ref{section: induced and coinduced from ell one}. Recall that
the effect of $\Phi_{S,-,k}$ on an element in $L ^1_{-,k}$ is
adding to it the transpose of its strictly lower triangular part.
Since $\mathbf{J}_k^S$ is a Poisson map and the functions
$I_l^{S,k}$ are in involution on $L ^1_{S,k}$, it follows that
$I_l^{S,k} \circ \mathbf{J}_k^S$ are also in involution on the weak
symplectic manifold $\left((\ell^{\infty} \times \ell^1)
\times \left(\ell ^{\infty} \right)^{k-2} \times \left(\ell ^1
\right)^{k-2}, \Omega_k\right)$ provided that these  functions
admit Hamiltonian vector fields.

\medskip

\noindent \textbf{The case $k=2$.} 
In this case we have $I ^1_{-,0,1} = L^1_{-,2}$ and $GI
^{\infty}_{+,0,1} = GL ^{\infty}_{+,2} $. As we discussed
earlier, the induction method yields in this case the
original weak symplectic manifold $\left(\ell ^{\infty}
\times \ell ^1, \omega \right)$. This is the case of the
standard semi-infinite Toda lattice.

\medskip

\noindent \textbf{The case $k=3$.} This is the first
situation that goes beyond  the Toda lattice. The 
Banach Lie group $G : = (GL^{\infty}_{+,3}, \circ_3)$
consists of bounded operators having only three upper
diagonals, while the operators in $GI^{\infty}_{+,0,2}$
have non-zero entries only on the main and the second
strictly upper diagonal. The induced space is now $(\ell
^{\infty} \times
\ell ^1)
\times (\ell ^{\infty}
\times \ell ^1)$. The $GL^{\infty}_{+,3}$-action on 
$(\ell^{\infty} \times \ell ^1) \times 
\left(\ell^{\infty} \times \ell^1\right)$ is given,
according to \eqref{k action q} - \eqref{k action p prime} by
\begin{align}
\label{3 action q}
\mathbf{q}' &=  \mathbf{q} + \log g_0\\
\label{3 action p}
\mathbf{p}' &= \mathbf{p}+ g_2 g_0 ^{-1}\nu_2
e^{s^2(\mathbf{q}) - \mathbf{q}} + g_1 s(\mathbf{q}_1)
g_0^{-1} \nu_2  e^{s(\mathbf{q}) - \mathbf{q}}  \nonumber\\
&  \qquad - \tilde{s}^{2} \left( g_2 g_0 ^{-1}\nu_2
e^{s^2(\mathbf{q}) - \mathbf{q}} + g_1 s(\mathbf{q}_1) g _0
^{-1} \nu_2 
 e^{s(\mathbf{q}) - \mathbf{q}}
\right)   \\
\label{3 action q prime}
\mathbf{q}_1' &= s(g_0 ^{-1})( g_1 + g_0\mathbf{q}_1) \\
\label{3 action p prime}
\mathbf{p}_1' &= s(g_0) g_0 ^{-1}\mathbf{p}_{1}, \qquad
l = 1, \dots, k-2.
\end{align}

The reduced symplectic form on 
$(\ell ^{\infty} \times \ell ^1) \times (\ell ^{\infty}
\times \ell ^1)$ is, according to \eqref{inverse k explicit},
\eqref{concrete k one form shifted}, and 
\eqref{concrete k symplectic form}, equal to
\begin{align}
\label{concrete 3 symplectic form}
\Omega_3 &= - \mathbf{d} \left[\operatorname{Tr}
\left(\mathbf{p} \mathbf{d}\mathbf{q}\right) +
\operatorname{Tr}\left(\mathbf{p}_1
\mathbf{d}\mathbf{q}_1\right) 
- \operatorname{Tr}\left(\nu_{2} e^{
s^{2}( \mathbf{q})- \mathbf{q}} \mathbf{q}_1s ( \mathbf{d}
\mathbf{q}_1)\right) \right] \nonumber\\
& = - \mathbf{d} \left[\operatorname{Tr}
\left(\mathbf{p} \mathbf{d}\mathbf{q}\right) +
\operatorname{Tr}\left(\left(\mathbf{p}_1 - \tilde{s} \left(
\nu_{2} e^{s^{2}( \mathbf{q})- \mathbf{q}} \mathbf{q}_1
\right)\right) \mathbf{d}\mathbf{q}_1\right) \right] 
\nonumber \\
& = - \mathbf{d} \left[\operatorname{Tr}
\left(\mathbf{p} \mathbf{d}\mathbf{q}\right) +
\operatorname{Tr}(\tilde{\mathbf{p}}_1
\mathbf{d}\mathbf{q}_1 \right],
\end{align}
where 
\begin{equation}
\label{tilde p one}
\tilde{\mathbf{p}}_1 := \mathbf{p}_1 - \tilde{s} \left(
\nu_{2} e^{s^{2}( \mathbf{q})- \mathbf{q}} \mathbf{q}_1
\right).
\end{equation}
We see here exactly the same phenomenon as in
classical electrodynamics, where a momentum shift by the
magnetic  potential transforms the non-canonical magnetic
symplectic form to the canonical one.

The equivariant momentum map 
\eqref{momentum on concrete k quotient} of this action is by
\eqref{group coordinate coadjoint minus k action} and 
\eqref{tilde p one} equal to
\begin{align}
\label{momentum on concrete 3 quotient}
&\mathbf{J}_3\left(\mathbf{q}, \mathbf{p}, \mathbf{q}_1,
\mathbf{p}_1\right)  
=  \left(\operatorname{Ad}^{+,3}
\right)^\ast_{(\mathbb{I} +
\mathbf{q}_1S)^{-1} }
\left(\mathbf{p} + S^T \mathbf{p}_1 + (S^T)^{2} \nu_{2}
e ^{s^{2} ( \mathbf{q}) - \mathbf{q}} \right) \nonumber\\
& \qquad = \mathbf{p}+ \mathbf{q}_1 \mathbf{p}_1 -
\tilde{s}\left(\mathbf{q}_1\mathbf{p}_1 + s( \mathbf{q}_1) \nu_2 
e^{s^2(\mathbf{q}) - \mathbf{q}} \mathbf{q}_1\right) + 
\tilde{s}^2 \left(\nu_2 e^{s^2(\mathbf{q}) - \mathbf{q}}
\mathbf{q}_1 s( \mathbf{q}_1) \right) \nonumber\\
& \qquad  \qquad  + S^T \left(\mathbf{p}_1 + s(
\mathbf{q}_1) \nu_2 e^{s^2(\mathbf{q}) - \mathbf{q}} -
\tilde{s}\left(\nu_2 e^{s^2(\mathbf{q}) - \mathbf{q}}
\mathbf{q}_1\right) \right) 
+ \left(S^T \right)^2 \nu_2e^{s^2(\mathbf{q}) - \mathbf{q}}
\nonumber \\
& \qquad = \mathbf{p} + \mathbf{q}_1 \mathbf{p}_1 -
\tilde{s} \left(\mathbf{q}_1 \mathbf{p}_1 \right) - \tilde{s}
\left(\nu_2 e^{s^2(\mathbf{q}) - \mathbf{q}} \mathbf{q}_1 \right)
\mathbf{q}_1 + \tilde{s}^2
\left(\nu_2 e^{s^2(\mathbf{q}) - \mathbf{q}} \mathbf{q}_1 \right)
\tilde{s}(\mathbf{q}_1)   \nonumber\\
& \qquad  \qquad + S^T \left(\mathbf{p}_1 + s(
\mathbf{q}_1) \nu_2 e^{s^2(\mathbf{q}) - \mathbf{q}} -
\tilde{s}\left(\nu_2 e^{s^2(\mathbf{q}) - \mathbf{q}}
\mathbf{q}_1\right) \right) 
+ \left(S^T \right)^2 \nu_2e^{s^2(\mathbf{q}) - \mathbf{q}}
\nonumber \\ 
& \qquad = \mathbf{p}+ \mathbf{q}_1 \tilde{\mathbf{p}}_1 -
\tilde{s}\left(\mathbf{q}_1 \tilde{\mathbf{p}}_1 \right) 
 + S^T \left(\tilde{\mathbf{p}}_1 + s(
\mathbf{q}_1) \nu_2 e^{s^2(\mathbf{q}) - \mathbf{q}}  \right) 
+ \left(S^T \right)^2 \nu_2e^{s^2(\mathbf{q}) - \mathbf{q}}
\end{align}
since the inverse of $\mathbb{I} + \mathbf{q}_1 S$ in the Banach
Lie group $GL ^{\infty}_{+,3}$ is equal to $(
\mathbb{I} +
\mathbf{q}_1S) ^{-1} =
\mathbb{I} - 
\mathbf{q}_1S + \mathbf{q}_1 s (\mathbf{q}_1) S^2 \in
GL^{\infty}_{+,3}$. 

The Hamiltonians $I^{S,3}_l$ given by \eqref{i l s k}
are in involution on $L ^1_{S,3}$ and hence the functions
$I_l^{S,3} \circ \mathbf{J}_3^S$ are in involution on 
$\left((\ell^{\infty} \times \ell ^1) \times
\left(\ell ^{\infty} \times \ell ^1 \right), \Omega_3 \right)$,
provided that they have Hamiltonian  vector  fields relative to the
weak symplectic  form $\Omega_3$. 

For  $l=1,2$, the Hamiltonians
$H_1 : = I^{S,3}_1 \circ \mathbf{J}_3^S$ and $H_2: = I^{S,3}_2
\circ \mathbf{J}_3^S$ have the expressions
\begin{equation}
\label{first degree 3 hamiltonian}
H_1( \mathbf{q}, \mathbf{p}, \mathbf{q}_1, \mathbf{p}_1) =
\operatorname{Tr}(\mathbf{p} )
\end{equation}
and 
\begin{align}
\label{second degree 3 hamiltonian}
H_2(\mathbf{q}, \mathbf{p}, \mathbf{q}_1, \mathbf{p}_1) 
& = \frac{1}{2} \operatorname{Tr} \left[\mathbf{p}+ \mathbf{q}_1 
\tilde{\mathbf{p}}_1 -
\tilde{s}\left(\mathbf{q}_1 \tilde{\mathbf{p}}_1 \right)
\right]^2 + \operatorname{Tr}\left(\tilde{\mathbf{p}}_1 + s(
\mathbf{q}_1) \nu_2 e^{s^2(\mathbf{q}) - \mathbf{q}}  \right)^2
\nonumber\\
& \qquad + \operatorname{Tr} \left(\nu_2  e^{s^2(\mathbf{q}) -
\mathbf{q}}\right)^2.
\end{align}
The Hamiltonian system defined by $H_2$ describes a semi-infinite
family of particles in an external field (given by the magnetic term
of the symplectic  form \eqref{concrete k symplectic form})
and where the interaction is between every second neighbor.
In the case of the Toda lattice (obtained for $k=2$, as
discussed above), there is no external field and the
interaction  is between nearest neighbors. The solution of
the semi-infinite Toda lattice will be given in \S 
\ref{toda example}. For arbitrary
$k$ there is an external field and interaction of 
particles is between every $(k-1)$st neighbor. 

We have given here only the first two Hamiltonians of an
infinite family of functions in involution. Involutivity follows
because they are obtained from a family of integrals in involution,
namely the $I_k^{S,3}$ by pull back with the Poisson map
$\mathbf{J}^S_3$.

\section{The semi-infinite Toda lattice}
\label{toda example}

In this section we illustrate the theory of the $k $-diagonal
Hamiltonian systems by the detailed investigation of the
semi-infinite Toda lattice which is an example of a bidiagonal
system (see Remark (v) at the end of \S 
\ref{sect: bidiagonal case}). We shall follow the method of
orthogonal polynomials first proposed in \cite{be}, as far as we know. We shall extend below the results in \cite{M} for the finite Toda lattice by explicitly solving the the semi-infinite Toda lattice both in action-angle variables as well as giving all the flows of the full hierarchy in the original variables.

The family of Hamiltonians $I_l^{S,2} \in C ^{\infty}(
L^1_{S,2})$, $l \in \mathbb{N}$, leads to the chain of
Hamilton equations
\begin{equation}
\label{toda hierarchy}
\frac{\partial}{\partial t_l} \boldsymbol{\rho} = \left[\boldsymbol{\rho},
B_l\right], \quad
\text{where} \quad B_l : =
P_-^{\infty}(\boldsymbol{\rho}^{l}) -
\left(P_-^{\infty}(\boldsymbol{\rho}^{l})\right)^T,
\end{equation}
on the Banach Lie-Poisson space $\left(L ^1_{S,2}, \{ \cdot ,
\cdot \}_{S,2}\right)$ (or on the space
$(L^1_{-,2}, \{ \cdot , \cdot \}_2)$ isomorphic to it)
induced from
\eqref{l s one i l 1} by the inclusion $\iota_{S,2}: L
^1_{S, 2}
\hookrightarrow L^1_S$.

The selfadjoint trace class operator $\boldsymbol{\rho}\in L
^1_{S, 2}$ acts on the orthonormal basis $\{|k \rangle\}_{k =
0} ^{\infty}$ of
$\mathcal{H}$ as follows:
\begin{equation}
\label{rho action}
\boldsymbol{\rho} |k \rangle = \rho_{k-1,k} |k-1\rangle +
\rho_{kk}|k \rangle + \rho_{k, k+1} |k+1 \rangle,
\end{equation}
where $k \in \mathbb{N} \cup \{0\}$ and we set $\rho_{-1,0} =
0$.

Note that if $\boldsymbol{\rho}$ is replaced by $\boldsymbol{\rho}' : = c\boldsymbol{\rho} + b \mathbb{I}$, where $b, c \in \mathbb{R}$, $c \neq 0$, 
then the equations \eqref{toda hierarchy} remain unchanged by rescaling the time $t'_ l : = c ^{-l} t_l $. As will be explained later,
the norm $\|\boldsymbol{\rho}\|_\infty$ and the positivity $\boldsymbol{\rho} \geq 0 $ are preserved by the
evolution defined by \eqref{toda hierarchy}. Taking into account  the above facts, we can assume, without loss of generality, that $\|\boldsymbol{\rho}\|_\infty < 1 $ and $\boldsymbol{\rho} \geq 0 $. Consequently, from now on we shall work with generic initial conditions $\boldsymbol{\rho}(0)$ for the Hamiltonian system \eqref{toda hierarchy}, i.e., 
\begin{equation}
\label{simple_spectrum}
\lambda_m(0) \neq \lambda_n(0), \quad \text{for} \quad n \neq m
\end{equation}
\begin{equation}
\label{norm_less_one}
\lambda_m(0) > 0 \quad \text{and} \quad \operatorname{sup}_{m \in \mathbb{N}\cup \{0\}} \{\lambda_m(0)\} <1,
\end{equation}
where $\lambda_m(0)$ are the eigenvalues of $\boldsymbol{\rho}(0)$. This means that $\boldsymbol{\rho}(0)$ has simple spectrum,  $\boldsymbol{\rho}(0)\geq 0 $, and $\|\boldsymbol{\rho}(0)\|_\infty <1$. These hypotheses imply that $\rho_{k, k+1} (0) >0$ for all $k \in \mathbb{N} \cup \{0\}$ 
 and are
consistent with the
physical interpretation of the semi-infinite Toda system. Let us denote
by $\Omega ^1_{-,2} \subset L^1_{S,2}$ the open set consisting of operators satisfying \eqref{simple_spectrum} and \eqref{norm_less_one}. 

From \eqref{rho action}, it follows that 
\begin{equation}
\label{rho polynomial relation}
|k \rangle = P_k(\boldsymbol{\rho}) |0\rangle,
\end{equation}
where the the polynomials $P_k(\lambda) \in \mathbb{R}[
\lambda]$, $k \in
\mathbb{N}\cup \{0\}$, are obtained by solving the three
term recurrence equation
\begin{equation}
\label{recurrence}
\lambda P_k(\lambda) = \rho_{k-1, k}P_{k-1}(\lambda) +
\rho_{kk}P_k(\lambda) + \rho_{k, k+1}P_{k+1}(\lambda)
\end{equation}
with initial condition $P_0(\lambda) \equiv 1$. Note that the
degree of $P _k(\lambda)$ is $k $.

We show now that that the operator $\boldsymbol{ \rho} \in L ^1_{S,2}$ evolving according to \eqref{toda hierarchy} also has simple spectrum independent of all times $t _l$. To do this, we write the spectral resolution
\begin{equation}
\label{spectral resolution}
\boldsymbol{\rho} = \sum_{m=0} ^{\infty} \lambda_m
\mathbb{P}_m, \qquad
\mathbb{P}_m\mathbb{P}_n= \delta_{mn} \mathbb{P}_n, \qquad
\sum_{m=0} ^{\infty}\mathbb{P}_m= \mathbb{I},
\end{equation}
where
\begin{equation}
\label{projection definition}
\mathbb{P}_m: = \frac{| \lambda_m\rangle \langle \lambda
_m|}{\langle \lambda _m | \lambda_m \rangle} 
\end{equation}  
are the projectors on the one-dimensional eigenspaces spanned 
by the eigenvector $| \lambda_m\rangle $. From \eqref{toda hierarchy} 
one  obtains
\begin{equation}
\label{intermediate projection}
\left(\frac{\partial}{\partial t_l}\lambda_k\right)\mathbb{P}_n\mathbb{P}_k
+ \left(\lambda_n- \lambda_k\right) \left[\left( 
\frac{\partial}{ \partial t_l} \mathbb{P}_n\right) \mathbb{P}_k -
\mathbb{P}_nB _l \mathbb{P}_k \right]  = 0
\end{equation}
for any $n , k \in \mathbb{N}\cup \{0\}$ and $l \in \mathbb{N}$.
Putting $n = k $ in \eqref{intermediate projection} one finds
\begin{equation}
\label{eigenvalue conservation}
\frac{\partial}{\partial t_l}\lambda_n= 0
\end{equation}
for any $n \in \mathbb{N}\cup \{0\}$ and $l \in \mathbb{N}$. Thus $\lambda_m = \lambda_m(0) \neq \lambda_n$ for $n \neq m $ and we can conclude that the coefficients in 
\begin{equation}
\label{m eigenvector}
|\lambda _m\rangle = \sum_{l=0} ^{\infty} P_l (\lambda_m )
|l \rangle
\end{equation}
are the values $P_l(\lambda_m)$
at the eigenvalue $\lambda_m$ of the polynomials
$P_l(\lambda)$ which are orthogonal relative to the $L ^2$-scalar
product given by the measure $\sigma$ in \eqref{def of sigma}. 

Taking $n \neq k$ in \eqref{intermediate projection} and using properties of orthogonal projectors one obtains
\begin{equation}
\label{time derivative of projectors}
\frac{\partial}{\partial t_l} \mathbb{P}_n = \left[\mathbb{P}_n, B _l
\right] \quad \text{for any} \quad n \in \mathbb{N}\cup \{0\}
\quad \text{and} \quad l \in \mathbb{N}.
\end{equation}

Similarly, for the resolvent
\begin{equation}
\label{resolvent definition}
R_\lambda : = (\boldsymbol{\rho} - \lambda\mathbb{I}) ^{-1}=
\sum_{m=0} ^{\infty} \frac{1}{\lambda_m- \lambda}
\mathbb{P}_m
\end{equation}
by \eqref{time derivative of projectors} one has 
\begin{equation}
\label{time derivative of resolvent}
\frac{\partial }{\partial t_l} R_\lambda =  \sum_{m=0} ^{\infty}
\frac{1}{\lambda_m- \lambda}
\left[\mathbb{P}_m, B _l\right]
= \left[R_\lambda , B_l\right].
\end{equation}

Note that \eqref{rho polynomial relation} implies that the
vector $|0 \rangle$ is cyclic for $\boldsymbol{\rho}$. Thus,
one has a unitary isomorphism of $\mathcal{H}$ with $L ^2(
\mathbb{R}, d\sigma)$, where the measure
\begin{equation}
\label{def of sigma}
d\sigma( \lambda) : = d\langle 0 | \mathbb{P}_ \lambda 0  
\rangle
= \sum_{m=0} ^{\infty}\mu_m \delta(\lambda- \lambda_m)
d\lambda,
\end{equation}
is given by the orthogonal resolution of the unity
$\mathbb{P} : \mathbb{R} \ni \lambda \mapsto \mathbb{P}_
\lambda \in L ^{\infty}( \mathcal{H})$ for 
\[
 \rho = \int _{ \mathbb{R}} \lambda d \mathbb{P}_ \lambda.
\]
The {\bfi masses} $\mu_m$ in \eqref{def of sigma} are given
by
\begin{equation}
\label{def of masses}
\mu_m^{-1} = \langle \lambda_m| \lambda_m\rangle =
\sum_{l=0} ^{\infty}\left(P_l( \lambda_m) \right)^2.
\end{equation}

Using  $\mathbb{P}_m|0\rangle = \mu_m| \lambda_m\rangle$  
and $\mu_m = \langle 0 | \mathbb{P}_m 0 \rangle$, one 
obtains from \eqref{time derivative of projectors} the 
 differential equation
\begin{equation}
\label{derivative of mu m}
\frac{\partial}{\partial t_l}\mu_m = 2\langle \lambda_m | B_l 0 \rangle
\mu_m 
= 2 \left( \lambda_m^l - \langle 0 | 
\boldsymbol{\rho}^l0 \rangle \right) \mu_m
\end{equation}
for any $l \in \mathbb{N}$  and $m \in \mathbb{N}\cup \{0\}$.
In order to prove the second equality in 
\eqref{derivative of mu m} we notice that 
\begin{equation}
\label{b l formula}
B _l = \boldsymbol{\rho}^l - P_0
^{\infty}(\boldsymbol{\rho}^l) - 2
\left[P^{\infty}_-(\boldsymbol{\rho}^l) \right]^T 
\end{equation}
\begin{equation}
\label{vanishing formula for p minus}
\left[P^{\infty}_-(\boldsymbol{\rho}^l) \right]^T |0 \rangle
= 0
\end{equation}
\begin{equation}
\label{p minus on zero}
P_0 ^{\infty}(\boldsymbol{\rho}^l) |0\rangle = \langle 0 |
P_0^{\infty}(\boldsymbol{\rho}^l)0 \rangle |0\rangle 
= \left\langle 0 |\boldsymbol{\rho}^l 0 \right\rangle |0\rangle 
\end{equation}
which implies
\begin{equation}
\label{inner product lambda m}
\left\langle \lambda_m| B _l0 \right\rangle
= \lambda_m^l - \left\langle 0| \boldsymbol{\rho}^l 0
\right\rangle.
\end{equation}

Using \eqref{derivative of mu m} and noticing that 
\begin{equation}
\label{sigma l formula}
\sigma_k = \left\langle 0| \boldsymbol{\rho}^k 0
\right\rangle 
\end{equation}
one obtains the system of equations
\begin{equation}
\label{moment equation k}
\frac{\partial}{\partial t_l}\sigma_k 
= 2 \left(\sigma_{k+l} -  \sigma_l \sigma_k\right), 
\end{equation}
where $\sigma_0 = 1$, $k \in \mathbb{N}\cup \{0\}$, $l \in \mathbb{N}$, for the moments   
\begin{equation}
\label{k moment  of sigma}
\sigma_k = \int_{\mathbb{R}} \lambda^k
d\sigma( \lambda) = \sum_{m=0} ^{\infty} \lambda_m^k
\mu_m
\end{equation}
of the measure \eqref{def of sigma}. Let
us remark here that in the considered case the moment
problem is determined, i.e., the moments $\sigma_k$ 
determine the measure \eqref{def of sigma} in a unique way
(see, e.g. \cite{Akh}).

Let us comment on the formulas obtained above. Introduce the
diagonal trace class operators $\boldsymbol{\lambda},
\boldsymbol{\mu}, \boldsymbol{\sigma}\in L ^1_0$ by defining
their $m^{\operatorname{th}}$ components to be the eigenvalues 
$\lambda_m$, the masses $\mu_m$, and the moments
$\sigma_m$, $m \in \mathbb{N}\cup \{0\}$, respectively. On the
open subset
$\Omega_{-,2} ^1$ one has three naturally defined
smooth coordinate systems: 
\begin{itemize}
\item[(i)] $\boldsymbol{\rho} \in \Omega_{-,2} ^1$,
\item [(ii)] $(\boldsymbol{\lambda}, \boldsymbol{\mu})
\in L ^1_0 \times L ^1_0$, where
$\operatorname{Tr}\boldsymbol{\mu}= 1$ and 
$\boldsymbol{\mu}>0$,
\item [(iii)] $\boldsymbol{\sigma} \in L ^1_0$ with first
component
$\sigma_0 = 1$, $\boldsymbol{\boldsymbol{\sigma}}>0$, and $\mathbf{d}_{0}>0$,
\end{itemize}
where $\mathbf{d}_0 : = \sum_{k=0}^{\infty}d_{0k}|k \rangle \langle k|$, and 
\begin{equation}
\label{delta zero n}
d_{0k} := \operatorname{det}
\left[
\begin{array}{cccccc}
\sigma_0& \sigma_1 & \sigma_2  &\sigma_3  & \dots &
\sigma_k\\ 
\sigma_1& \sigma_2 & \sigma_3 & \sigma_4  & \dots &
\sigma_{k+1}\\ 
\sigma_2  & \sigma_3 & \sigma_4 & \sigma_5& \dots &
\sigma_{k+2}\\ 
\sigma_3 & \sigma_4   & \sigma_5 & \sigma_6 & \dots &
\sigma_{k+3}\\
\vdots & \vdots & \vdots& \vdots & \vdots & \vdots\\
\sigma_k & \sigma_{k+1} & \sigma_{k+2} & \sigma_{k+3}&
\dots & \sigma_{2k}
\end{array}
\right] > 0, 
\end{equation}
with the convention that $d_{0, -1} = 1$. 
In order to see that $\boldsymbol{\sigma} \in L ^1_0$ we notice that
\[
\sum_{k=0} ^{\infty}\sigma_k 
= \sum_{k=0} ^{\infty}\langle 0| \boldsymbol{\rho} ^k 0 \rangle 
\leq \sum_{k=0} ^{\infty}\|\boldsymbol{\rho} ^k \| 
\leq \sum_{k=0} ^{\infty}\|\boldsymbol{\rho} \|^k
= \frac{1}{1-\|\boldsymbol{ \rho} \|_\infty} < + \infty. 
\]

We also define $\mathbf{d}_1: = 
\sum_{k=0} ^{\infty} d_{1k}|k \rangle \langle k|$, where 
\begin{equation}
\label{delta one n}
d_{1k} := \operatorname{det}
\left[
\begin{array}{ccccccc}
\sigma_0& \sigma_1 & \sigma_2 & \sigma_3  & \dots &
\sigma_{k-1} & \sigma_{k+1}\\ 
\sigma_1& \sigma_2 & \sigma_3 & \sigma_4  & \dots &
\sigma_{k} & \sigma_{k+2}\\ 
\sigma_2  & \sigma_3 & \sigma_4
& \sigma_5& \dots &
\sigma_{k+1} & \sigma_{k+3}\\ 
\sigma_3 & \sigma_4   & \sigma_5 & \sigma_6 & \dots &
\sigma_{k+2} & \sigma_{k+4}\\
\vdots & \vdots & \vdots& \vdots & \vdots & \vdots\\
\sigma_k & \sigma_{k+1} & \sigma_{k+2} & \sigma_{k+3}&
\dots & \sigma_{2k-1} & \sigma_{2k+1}
\end{array}
\right]
\end{equation}
for $n \in \mathbb{N} \cup \{0\}$.

The transformation from $\boldsymbol{\rho}$-coordinates to 
$\boldsymbol{\boldsymbol{\sigma}}$-coordinates is given by formula 
\eqref{sigma l formula}. The inverse transformation to 
\eqref{sigma l formula} has the form 
\begin{align}
\label{sigma to rho}
\boldsymbol{\rho} &= S^T \rho_1+ \rho_0 + \rho_1S \nonumber \\
&= S^T  \left[\tilde{s}(\mathbf{d}_0) s(\mathbf{d}_0) \right]^{1/2}
\mathbf{d}_0 ^{-1} + \mathbf{d}_0 ^{-1}
\mathbf{d}_1 - \tilde{s}(\mathbf{d}_0 ^{-1}\mathbf{d}_1) +
\left[\tilde{s}(\mathbf{d}_0) s(\mathbf{d}_0) \right]^{1/2}
\mathbf{d}_0 ^{-1} S, 
\end{align}
or, in components (see, e.g.,  \cite{Akh}),
\begin{equation}
\label{sigma to rho components}
\rho_{kk} = d_{0k}^{-1}d_{1k} - d_{0, k-1}^{-1} d_{1, k-1}
\qquad \text{and} \qquad \rho_{k, k+1} = \left(d_{0, k-1}
d_{0, k+1} \right)^{1/2} d_{0k} ^{-1} >0.
\end{equation}

Formula \eqref{k moment  of sigma} gives the transformation
from $(\boldsymbol{\lambda}, \boldsymbol{\mu}) $-coordinates
to $\boldsymbol{\boldsymbol{\sigma}}$-coordinates. The inverse
transformation to \eqref{k moment  of sigma} is obtained by expanding the
so-called Weyl function $\langle 0 | R_\lambda 0 \rangle$ in a Laurent series
\begin{equation}
\label{Laurent Weyl}
\langle 0 | R_\lambda 0 \rangle = \sum_{m=0} ^{\infty}
\frac{ \mu_m}{ \lambda_m- \lambda} = - \sum_{k=0} ^{\infty}
\frac{\sigma_k}{\lambda^{k+1}}
\end{equation}
for $|\lambda|> \operatorname{sup}_{m \in \mathbb{N}\cup
\{0\}} \{| \lambda_m |\} = \|\boldsymbol{\rho}\|_\infty$. So, one finds $(\boldsymbol{\lambda},
\boldsymbol{\mu}) $ by computing the Mittag-Leffler
decomposition of the left hand side of \eqref{Laurent Weyl}.

The passage from $\boldsymbol{\rho}$-coordinates to
$(\boldsymbol{\lambda}, \boldsymbol{\mu})$-coordinates is
obtained by composing the previously described
transformations. This can also be done directly constructing the
spectral resolution for $\boldsymbol{\rho}$.

After these remarks  we present Hamilton's  equations
\eqref{toda hierarchy} in the coordinates
$(\boldsymbol{\lambda}, \boldsymbol{\mu})$
\begin{align}
\label{lambda equ}
\frac{\partial}{\partial t_l} \boldsymbol{\lambda} &=
\{\boldsymbol{\lambda}, I _l^{S,2}\}_{S,2} = 0\\
\label{mu equ}
\frac{\partial}{\partial t_l} \boldsymbol{\mu} & =\{\boldsymbol{\mu}, 
I_l^{S,2}\}_{S,2} = 2 \left(\boldsymbol{\lambda}^l -
\operatorname{Tr}(\boldsymbol{\lambda}^l \boldsymbol{\mu}) 
\right) \boldsymbol{\mu}
\end{align}
or, in components,
\begin{equation}
\label{lambda equ components}
\frac{\partial}{\partial t_l} \lambda_m= 0 \quad \text{and} \quad
\frac{\partial}{\partial t_l}\mu_m = 2 \left(\lambda_m^l -
\sum_{n=0}^{\infty} \lambda_n^l \mu_n \right)\mu_m
\end{equation}
and in the coordinates $\boldsymbol{\boldsymbol{\sigma}}$
\begin{equation}
\label{sigma Hamiltonian equ}
\frac{\partial}{\partial t_l} \boldsymbol{\boldsymbol{\sigma}} =
\{\boldsymbol{\sigma},  I_l^{S,2} \} = 2 \left(s
^l(\boldsymbol{\sigma}) -
\sigma_l\boldsymbol{\sigma}\right)
\end{equation}
whose coordinate expression was already given in \eqref{moment equation k}.
In deducing equations \eqref{lambda equ}, \eqref{mu equ},
and \eqref{sigma Hamiltonian equ} we used 
\eqref{sigma l formula} and \eqref{k moment  of sigma}.

Let us observe now that \eqref{moment equation k} implies that
\begin{equation}
\label{c m l k}
\frac{\partial \sigma_k}{\partial t_l}= \frac{\partial
\sigma_l}{\partial t_k}
\end{equation}
for $k,l\in \mathbb{N}$. Thus there exists a function $\tau(t _1, t _2,
\ldots )$ of infinitely many variables $(t _1, t _2, \ldots) =: {\bf t} \in \ell^{\infty}$ such that 
\begin{equation}
\label{k der tau}
\sigma_k = \frac{1}{2} \frac{ \partial}{ \partial t _k}
\log \tau, \quad  k \in \mathbb{N}.
\end{equation}
In order to be consistent with the notation assumed in the theory of integrable systems (see, e.g. \cite{MJ, N}), we have called this function $\tau$-function. 

Substituting \eqref{k der tau} into 
\eqref{moment equation k} we obtain the system of linear partial
differential equations
\begin{equation}
\label{hirota}
\frac{ \partial^2 \tau }{\partial t _l \partial t _k} = 
2  \frac{\partial \tau}{ \partial t_{k+l}}, \quad  k,l \in \mathbb{N},
\end{equation}
on the $\tau$-function.

In order to find the explicit form of the $\tau$-function,
use \eqref{k moment  of sigma}, substitute 
\eqref{k der tau} into \eqref{lambda equ components}, and
integrate both sides of the resulting equation to get
\begin{align}
\label{mu m solution}
&\mu_m(t _1, t _2, \ldots, t_{l-1}, t _l, t _{l+1}, \ldots)
\nonumber \\
&  \quad = \mu_m(t _1, t _2, \ldots, t_{l-1}, 0, t _{l+1},
\ldots)
\frac{ \tau(t _1, t _2, \ldots, t_{l-1}, 0, t _{l+1},
\ldots)}{ \tau(t _1, t _2, \ldots, t_{l-1}, t _l, t _{l+1},
\ldots)} e^{2 \lambda_m^l t _l}.
\end{align}
Iterating \eqref{mu m solution} relative to $l \in
\mathbb{N}$ yields the final formula for $\mu_m(t _1, t _2,
\ldots)$, namely
\begin{equation}
\label{mu m final solution}
\mu_m(t _1, t _2, \ldots)
= \mu_m(0, 0 , \ldots) \frac{ \tau(0, 0, \ldots,)}{
\tau(t _1, t _2, \ldots)} e^{2 \sum_{l=1}^{\infty}\lambda_m^l
t _l}.
\end{equation}
Since  $\sum_{m=0} ^{\infty} \mu_m(t _1, t _2, \ldots) = 1$,
we get the following expression for the $\tau$-function
\begin{equation}
\label{tau function formula final}
\tau(t _1, t _2, \ldots) = \tau(0, 0, \ldots)
\sum_{m=0}^{\infty}
\mu_m(0,0, \ldots)e^{2 \sum_{l=1}
^{\infty}\lambda_m^l t _l}
\end{equation}

Let us show that the series in \eqref{tau function formula final} 
is convergent if $\boldsymbol{ \mu}(0) \in L ^1_0 \cong \ell ^1$ and 
${\bf t} \in \ell^\infty$. In order to do this we prove that the linear operator
defined by 
\[
\left(\Lambda {\bf t} \right)_m : = \sum_{l=1} ^{\infty} \lambda_m^l t _l
\]
is bounded on $\ell ^{\infty}$. This follows from
\begin{align*}
\| \Lambda{\bf t}\|_\infty & = 
\sup_{m \in \mathbb{N}} \left|\sum_{l=1} ^{\infty} \lambda_m^l t _l \right|
\leq \| {\bf t}\|_\infty \sup_{m \in \mathbb{N}} \left|\sum_{l=1} ^{\infty} \lambda_m^l  \right| \\
& = \| {\bf t}\|_\infty \sup_{m \in \mathbb{N}} \frac{ \lambda_m}{1-\lambda_m} 
= \| {\bf t}\|_\infty  \frac{\| \boldsymbol{\rho} \|_\infty}{1-\| \boldsymbol{\rho} \|_\infty} .
\end{align*}
Thus the sequence $\{e^{2 \sum_{l=1}
^{\infty}\lambda_m^l t _l}\}_{m \in \mathbb{N}} \in \ell ^{\infty}$.
Since $\{\mu_m(0, 0, \ldots) \}_{m \in \mathbb{N}} \in \ell ^1$, 
the series in \eqref{tau function formula final} converges.

Summarizing, we see that the substitution of 
\eqref{tau function formula final} into \eqref{k der tau} and
\eqref{mu m solution} gives the $\mathbf{t} : = ( t _1, t_2,
\ldots)$-dependence of the moments $\sigma_k(\mathbf{t})$
and the masses $\mu_m(\mathbf{t})$, respectively. The
dependence of $\rho_{kk}(\mathbf{t})$ and $\rho_{k, k+1} (
\mathbf{t}) $ on $\mathbf{t}$ is given by 
\eqref{sigma to rho}, \eqref{delta zero n}, and 
\eqref{delta one n} which express these quantities in terms
of $\sigma_m( \mathbf{t}) $. From the discussion above we see that the
conditions \eqref{simple_spectrum}, \eqref{norm_less_one}
are preserved by the
$\mathbf{t}$-evolution. 

Next, using
\eqref{k der tau}, \eqref{mu m final solution}, and the
formula
\begin{equation}
\label{polynomial sigma}
P_n(\lambda_m) = \frac{1}{\sqrt{d_{0,n-1}d_{0,n}}}
\operatorname{det}
\left[
\begin{array}{ccccc}
\sigma_0 & \sigma_1 & \sigma_2 & \dots & \sigma_n\\ 
\sigma_1 & \sigma_2 & \sigma_3 & \dots & \sigma_{n+1}\\ 
\sigma_2 & \sigma_3 & \sigma_4 & \dots & \sigma_{n+2}\\ 
\vdots & \vdots & \vdots& \vdots & \vdots \\
\sigma_{n-1} & \sigma_{n} & \sigma_{n+1} & \dots &
\sigma_{2n-1}\\ 1 & \lambda_m & \lambda_m^2 & \dots &
\lambda_m^n
\end{array}
\right]
\end{equation}
obtained by orthonormalizing the monomials $\lambda^n$, $n
\in\mathbb{N}\cup\{0\}$, with respect to the  measure $\sigma$ 
(see, e.g.,  \cite{Akh}), we
obtain from \eqref{m eigenvector} the
$\mathbf{t}$-dependence of the eigenvectors $|
\lambda_m(\mathbf{t}) \rangle$ and the corresponding
projectors $\mathbb{P}_m(\mathbf{t}) $, $m \in \mathbb{N}\cup
\{0\}$.

Formula \eqref{m eigenvector} defines the operator
$O: \mathcal{H} \rightarrow \mathcal{H}$ whose matrix in the
basis $\{|k\rangle\}_{k=0}^{\infty}$ is given by $O_{kl}(
\mathbf{t}) : = P_l( \mathbf{t})(\lambda_k)$. One has the
following identities
\begin{align}
\label{rho-o}
&\boldsymbol{\rho}( \mathbf{t}) O( \mathbf{t}) =
O (\mathbf{t}) \boldsymbol{\lambda}( \mathbf{t}) \\
\label{o-mu}
& O( \mathbf{t}) \boldsymbol{\mu}( \mathbf{t}) 
O(\mathbf{t})^T = \mathbb{I}
\end{align}
relating the operators $\boldsymbol{\rho}( \mathbf{t}) $,
$\boldsymbol{\lambda}( \mathbf{t}) $, $\boldsymbol{\mu}(
\mathbf{t}) $, and $O( \mathbf{t}) $ for any $\mathbf{t}$.
Since $\boldsymbol{\lambda}( \mathbf{t}) =
\boldsymbol{\lambda}(\boldsymbol{0})$, where $\boldsymbol{0}
:= (0, 0, \ldots)$, we obtain from
\eqref{rho-o} and \eqref{o-mu}
\begin{equation}
\label{rho decomposition}
\boldsymbol{\rho}( \mathbf{t}) = O( \mathbf{t})
O(\boldsymbol{0}) ^{-1}
\boldsymbol{\rho}(\boldsymbol{0}) \left(O( \mathbf{t})
O(\boldsymbol{0}) ^{-1}
\right)^{-1} = Z( \mathbf{t})^T
\boldsymbol{\rho}(\boldsymbol{0}) Z(\mathbf{t}),
\end{equation}
where $Z( \mathbf{t}) : = O(\boldsymbol{0}) \boldsymbol{\mu}(
\boldsymbol{0}) ^{1/2} \left(O( \mathbf{t}) \boldsymbol{\mu}(
\mathbf{t})^{1/2}\right)^T$ is an orthonormal operator, i.e.,
$Z( \mathbf{t}) ^T Z( \mathbf{t}) = \mathbb{I}$. As shown in
\S \ref{section: induced and coinduced from ell one} and \S 
\ref{section: dynamics generated by Casimirs of ell one}, one
can express the flows $\mathbf{t} \mapsto \boldsymbol{\rho}(
\mathbf{t})$ through the coadjoint action
$\left(\operatorname{Ad}^{S,2} \right)^\ast:
GL^{\infty}_{+,2} \rightarrow \operatorname{Aut}
\left(L^1_{S,2} \right)$ of the bidiagonal group
$GL^{\infty}_{+,2}$ on the Banach Lie-Poisson space $L^1_{S,2}
\cong L ^1_{-,2}$, i.e., 
\begin{align}
\boldsymbol{\rho}( \mathbf{t}) & = 
\left(\operatorname{Ad}^{S,2} \right)^\ast_{g( \mathbf{t})
^{-1}} \boldsymbol{\rho}( \boldsymbol{0}) \nonumber \\
&= S^T s(g_0( \mathbf{t})) g_0( \mathbf{t}) ^{-1} \rho_1(
\boldsymbol{0}) + \rho_0( \boldsymbol{0}) +  g_0(
\mathbf{t})^{-1} g_1( \mathbf{t}) \rho_1( \boldsymbol{0}) -
\tilde{s} \left(g_0( \mathbf{t})^{-1} g_1( \mathbf{t})
\rho_1( \boldsymbol{0}) \right) \nonumber\\
& \qquad + s(g_0(\mathbf{t})) g_0(\mathbf{t})^{-1}
\rho_1(\boldsymbol{0}) S \nonumber\\
&= \sum_{i=0}^\infty \rho_{i,i+1}( \boldsymbol{0})
\frac{g_{i+1, i+1}( \mathbf{t})}{g_{ii}( \mathbf{t})} |i + 1
\rangle \langle i | \nonumber\\
& \qquad + \sum_{i=0}^\infty
\left(\rho_{ii}( \boldsymbol{0}) +
\rho_{i,i+1}(\boldsymbol{0})
\frac{g_{i+1,i}( \mathbf{t}) }{g_{ii}( \mathbf{t})} -
\rho_{i,i+1}( \boldsymbol{0})
\frac{g_{i+1,i}( \mathbf{t}) }{g_{i+1,i+1}( \mathbf{t})}
\right)|i \rangle \langle i | \nonumber \\
& \qquad + \sum_{i=0}^\infty \rho_{i,i+1}( \boldsymbol{0})
\frac{g_{i+1, i+1}( \mathbf{t})}{g_{ii}( \mathbf{t})} |i
\rangle \langle i+1 |
\end{align}
(the symmetric version of 
\eqref{group k coadjoint action}), where $\rho_0: =
\operatorname{diag}( \rho_{00}$,
$\rho_{11}, \ldots), \rho_1 : = \operatorname{diag}(\rho_{01},
\rho_{12}, \ldots)$, $g_0 : = ( g_{00}, g_{11}, \ldots)$, and
$g_1 : = (g_{10}, g_{21}, \ldots) \in L ^1_0$.

In order to find the time dependence $\mathbf{t} \mapsto g(
\mathbf{t}) = g_0 ( \mathbf{t}) + g_1( \mathbf{t})S$ for $g(
\mathbf{t}) \in GL ^{\infty}_{+,2}$ let us note that from
\eqref{group k coadjoint action} and the three term
recurrence  relation \eqref{recurrence} it follows that
\begin{align}
\label{group solution diagonal}
g_{kk}( \mathbf{t}) &= g_{00}(\mathbf{t}) \frac{\rho_{00}(
\mathbf{t}) \cdots \rho_{k-1,k-1}(\mathbf{t})}{\rho_{00}(
\mathbf{0}) \cdots \rho_{k-1,k-1}(\mathbf{0})} \nonumber\\ 
&= g_{00}(\mathbf{t}) \frac{P_{kk}(
\mathbf{0})}{P_{kk}(\mathbf{t})}
= g_{00}(\mathbf{t}) \sqrt{\frac{d_{0, k-1}( \boldsymbol{0})
d_{0k}(\mathbf{t})}{d_{0k}(\boldsymbol{0})
d_{0,k-1}(\mathbf{t})} }
\end{align}
and 
\begin{align}
\label{group solution upper diagonal}
g_{k+1,k}( \mathbf{t}) &= g_{00}(\mathbf{t}) \left(\frac{
\rho_{00}( \mathbf{t}) \cdots
\rho_{k-1,k-1}(\mathbf{t})}{\rho_{00}(
\mathbf{0}) \cdots \rho_{k-1,k-1}(\mathbf{0})}\right)\left(
\frac{\rho_{00} (\mathbf{t}) + \dots + \rho_{kk}(
\mathbf{t}) - \rho_{00} (
\boldsymbol{0}) - \cdots - \rho_{kk}( \boldsymbol{0})}{
\rho_{kk}( \boldsymbol{0})}\right) \nonumber \\
& = g_{00}(\mathbf{t}) \frac{ P_{k+1, k}( \boldsymbol{0})
P_{k+1, k+1} ( \mathbf{t}) - P_{k+1, k}( \mathbf{t}) P_{k+1,
k+1}( \boldsymbol{0})}{P_{kk}( \mathbf{t}) P_{k+1, k+1}(
\mathbf{t}) } \nonumber \\
& = g_{00}(\mathbf{t}) \frac{d_{1k}(\mathbf{t}) 
\sqrt{ d_{0,k+1}(\boldsymbol{0})} - 
d_{1k}(\boldsymbol{0}) \sqrt{d_{0,k+1}(\mathbf{t})}}
{\sqrt{d_{0k}( \boldsymbol{0}) d_{0k}( \mathbf{t})
 d_{0,k-1}( \mathbf{t}) d_{0,k+1}(
\boldsymbol{0})}}\,,
\end{align}
where $P_{kl}( \mathbf{t}) $ are the coefficients of the
polynomial $P _n( \mathbf{t}) ( \lambda) = P_{nn}(
\mathbf{t}) \lambda^n + P_{n,n-1}( \mathbf{t}) \lambda^{n-1}
+ \dots + P_{n1}( \mathbf{t})  \lambda + P_{n0}( \mathbf{t})
$. The last equalities in \eqref{group solution diagonal} and
\eqref{group solution upper diagonal} are obtained using
\eqref{polynomial sigma},  \eqref{delta zero n}, and
\eqref{delta one n} to get the expressions
\begin{align*}
P_{kk}( \mathbf{t}) = \sqrt{\frac{d_{0,k-1}(
\mathbf{t})}{d_{0k}(
\mathbf{t})} } \qquad \text{and} \qquad
P_{k+1, k}( \mathbf{t}) =  \frac{-d_{1k}(\mathbf{t})}{
\sqrt{d_{0k}(
\mathbf{t}) d_{0, k+ 1}( \mathbf{t})}}.
\end{align*}
Recall that $d_0( \mathbf{t}) $ and $d_1( \mathbf{t}) $ are
given by \eqref{delta zero n} and \eqref{delta one n},
respectively.

Finally, taking in \eqref{sigma action} (for $k=2$) $g_0(
\mathbf{t}) $ and $g_1( \mathbf{t}) $ given by 
\eqref{group solution diagonal} and
\eqref{group solution upper diagonal}, we obtain the explicit
expression for the time evolution of the position
$\mathbf{q}( \mathbf{t}) $ and the momentum $\mathbf{p}(
\mathbf{t}) $ for all flows in the Toda hierarchy described
by the Hamiltonians 
\[
H _l( \mathbf{q}, \mathbf{p}) : =
\left(I _l^{S,2} \circ \mathcal{J}_{ \nu_1} \right)(
\mathbf{q}, \mathbf{p}),
\]
where $\mathcal{J}_{ \nu_1}: \ell ^1\times \ell ^{\infty} 
\rightarrow L ^1_{-,2} \cong L ^1_{S,2}$ is the Flaschka map
given by \eqref{flaschka map} for $k=2$ and $I _l^{S,2} =
I_l^S \circ \iota_{S,2} = I_l \circ \iota_S \circ
\iota_{S,2}$ are the restrictions to $L ^1_{S,2}$ of the
Casimir functions $I _l$ of $L ^1$ (see \eqref{i l s k}).

Note that the formulas giving the group element
$g(\mathbf{t})$ depend on $g_{00}( \mathbf{t})$. This first
component cannot be determined but it does not matter because
$g_{00}( \mathbf{t}) \mathbb{I}$ is in the center of $GL
^{\infty}_{+,2}$ and hence the coadjoint action defined by it
is trivial. Also, in terms of the variables $\mathbf{q}$ and
$\mathbf{p}$, the action of this group element is a
translation in $\mathbf{q}$ and has no effect on
$\mathbf{p}$. This corresponds to the flow of $I _1^{S,2}$.
\medskip

To solve the Toda system one takes an initial condition
$\boldsymbol{\rho}(\boldsymbol{0})$ which determines a
coadjoint orbit of $GL ^{\infty}_{+,2}$ in $L ^1_{S,2}$. These
coadjoint  orbits were studied in detail in \S 
\ref{sect: bidiagonal case}. In the generic case, when all
entries on the strictly upper (and hence also strictly lower)
diagonal of
$\boldsymbol{\rho}(\boldsymbol{0})$ are strictly positive, the
solution of the Toda lattice was given above. If some upper
diagonal entries of $\boldsymbol{\rho}(\boldsymbol{0})$
vanish, Remark (iv) at the end of \S 
\ref{sect: bidiagonal case} describes such orbits as  blocks,
some of them finite and at most one  infinite. Then the
Toda lattice equations decouple and we get a smaller Toda
system for each block. On the infinite block, the solution
is as above. On each finite block one obtains a finite
dimensional Toda lattice whose solution is known (see, e.g., 
\cite{K, M, N, S}). The method we used above for the
semi-infinite case can be also used in the finite
case; one works then with measures $\sigma$ having finite
support and uses finite orthogonal polynomials. If one 
implements the solution method described in this section to
this finite dimensional case the results in \cite{M} are
reproduced.

\bigskip

\addcontentsline{toc}{section}{Acknowledgments}
\noindent\textbf{Acknowledgments.} This work was begun
while both authors were at the  Erwin
Schr\"odinger International Institute for Mathematical
Physics in the Fall of 2003 during the program \textit{The
Geometry of the Moment Map\/} and hereby thank ESI for its
hospitality. Some of the work on this paper was done
during the program \textit{Geometric Mechanics\/} at the
Bernoulli Center of the EPFL in the Fall of 2004. A.O.
thanks the Bernoulli Center for its hospitality and
excellent working conditions during his extended stay
there. We are grateful to D. Belti\c{t}\u{a} and H. Flaschka for several
useful discussions that influenced our presentation. The
 authors thank  the Polish and Swiss National Science
Foundations  (Polish State Grant P03A 0001 29 and Swiss NSF Grant 200021-109111/1) for partial support.

\end{document}